\documentclass[12pt]{amsart}
\usepackage{geometry, amssymb, dsfont, amsaddr, textcomp}
\usepackage[all]{xy} 
\title[Iwasawa theory for totally real fields]{The main conjecture of Iwasawa theory for totally real fields}
\author{Mahesh Kakde}
\date{\today}
\address{University College London}
\email{kakdem@math.ucl.ac.uk}

\newtheorem{theorem}{Theorem}
\newtheorem{lemma}[theorem]{Lemma}
\newtheorem{proposition}[theorem]{Proposition}
\newtheorem{remark}[theorem]{Remark}
\newtheorem{notation}[theorem]{Notations}
\newtheorem{definition}[theorem]{Definition}
\newtheorem{corollary}[theorem]{Corollary}

\newcommand{\ilim}[1]{%
	\displaystyle{%
	\lim_{\genfrac{}{}{0pt}{}{\longleftarrow}{\scriptstyle #1}} }\;}

\begin{document}

\maketitle

\begin{abstract}  The purpose of this paper is to prove the main conjecture of noncommutative Iwasawa theory for $p$-adic Lie extensions for totally real fields, for an odd prime $p$, assuming Iwasawa's conjecture on vanishing of the cyclotomic $\mu$-invariant for certain CM extension of the base field. 
\end{abstract}{

\tableofcontents

\section{Introduction}Iwasawa theory studies the mysterious relationship  between purely arithmetic objects and special values of complex $L$-functions. A precise form of this relationship is usually called the ``main conjecture". Historically, the first version of this main conjecture was stated for the cyclotomic $\mathbb{Z}_p$-extension of totally real number fields. In the last decade this main conjecture has been vastly generalised to all $p$-adic Lie extensions of a totally real number field containing the cyclotomic $\mathbb{Z}_p$-extension of the base field. The aim of this paper is to prove the main conjecture of ``noncommutative Iwasawa theory" for $p$-adic Lie extensions, for an odd prime $p$, of totally real number fields assuming that the Iwasawa $\mu$ invariant of a certain CM extension of the base field vanishes (see definition \ref{hypothesisonmu}). Throughout this paper $p$ is a fixed odd prime number.

The details of contents in the paper are as follows: In section \ref{sectionstatement} we set up the notation, recall some definitions and results from algebraic $K$-theory and state our main theorem. This is the main conjecture of Iwasawa theory for totally real number fields. 

In section \ref{sectionclassicalmainconjecture} we recall the results proven by Deligne-Ribet \cite{DeligneRibet:1980} , Cassou-Nogues \cite{Cassou:1979}, and Barsky \cite{Barsky:1978} on the construction of $p$-adic zeta functions in commutative Iwasawa theory and the statement of the main conjecture of classical Iwasawa theory as proven by Wiles \cite{Wiles:1990}, Mazur-Wiles \cite{MazurWiles:1984} and Rubin \cite{Lang:1990}. In this section we also prove that under the assumption `hypothesis $\mu = 0$ (definition \ref{hypothesisonmu}) is true' the main theorem for extensions whose Galois group is of the form $\mathcal{G} \times \Delta$, where $\mathcal{G}$ is an abelian p-adic Lie group and $\Delta$ is a finite group of order prime to $p$, can be easily deduced from the result of Wiles, Mazur-Wiles and Rubin (see theorem \ref{primetopmainconjecture}). 

In section \ref{sectionreductions} we prove several algebraic results and obtain various reductions. Burns and Kato made a beautiful observation which is used in all our reduction steps. It can be described as follows: for an open subgroup $U$ of $\mathcal{G}$ and a closed normal subgroup $V$ of $U$, then we have a map
\[
\theta_{U,V} : K_1'(\Lambda(\mathcal{G})) \rightarrow K_1'(\Lambda(U)) \rightarrow K_1'(\Lambda(U/V)),
\]
which is a composition of the norm map and the map induced by natural surjection. Similarly, we have a map
\[
\theta_{U,V,S}: K_1'(\Lambda(\mathcal{G})_S) \rightarrow K_1'(\Lambda(U/V)_S),
\]
where $S$ is the canonical Ore set defined by Coates et.al. \cite{CFKSV:2005} (see subsection \ref{subsectionsetup}). If we wish to reduce the proof of main conjecture for an extension with Galois group $\mathcal{G}$ to a class of subquotients, say $Cl(\mathcal{G})$, we must describe the kernel and the image of the map
\[
(\theta_{U,V})_{U/V \in Cl(\mathcal{G})}: K_1'(\Lambda(\mathcal{G})) \rightarrow \prod_{Cl(\mathcal{G})}K_1'(\Lambda(U/V)),
\]
and prove that 
\[
Im((\theta_{U,V,S})_{Cl(\mathcal{G})}) \cap \prod_{Cl(\mathcal{G})} \iota(K_1'(\Lambda(U/V))) = \iota(Im((\theta_{U,V})_{Cl(\mathcal{G})})),
\]
where $\iota$ is the natural map $K_1'(\Lambda(\cdot)) \rightarrow K_1'(\Lambda(\cdot)_S)$. The above intersection takes place in $\prod_{Cl(\mathcal{G})} K_1'(\Lambda(U/V)_S)$.

 We first reduce the proof of the main conjecture for arbitrary $p$-adic Lie extensions to those of one dimensional $p$-adic Lie extensions (theorem \ref{reductiontodim1}). This reduction is a result of  Burns \cite{Burns:2010} presented differently here. Next we reduce the proof of the main conjecture of one dimensional $p$-adic Lie extensions to the case when the $p$-adic Lie group has a $\mathbb{Q}_p$-elementary quotient (see definition \ref{defnhyperelementary}) by an open central pro-cyclic subgroup (theorem \ref{reductiontohyperelementary}). This reduction uses induction theory of Dress \cite{Dress:1973}. The $l$-$\mathbb{Q}_p$-elementary case, for $l \neq p$, is unsurprisingly easier and is handled first by the algebraic result of proposition \ref{propprimetopk1} and the theorem \ref{primetopmainconjecture}. The $p$-$\mathbb{Q}_p$-elementary case is reduced to the case pro-$p \times \Delta$, where $\Delta$ is a finite cyclic group of order prime to $p$, using ideas from Oliver \cite{Oliver:1988} chapter 12. This is theorem \ref{theoremreductiontoprop}.

In section \ref{sectioncomputationofk1} we describe the Whitehead group of the Iwasawa algebra of a one dimensional pro-$p$, $p$-adic Lie group in terms of the Whitehead groups of Iwasawa algebras of abelian subquotients of the group. We also prove a result about Whitehead groups of a localisation of the Iwasawa algebra used in the formulation of the main conjecture (see definition \ref{defnofphi} and theorems \ref{theorem1} and \ref{theorem2}). Such theorems were first proven by Kato \cite{Kato:2005} for Iwasawa algebra of open subgroups of $\mathbb{Z}_p \rtimes \mathbb{Z}_p^{\times}$. The strategy we use for proving these theorems is a generalisation of the strategy Kato used in \cite{Kato:2006} for $p$-adic Heisenberg groups. It uses the integral logarithm of Oliver and Taylor. This strategy was used by the author to prove such algebraic theorem for pro-$p$ meta-abelian groups in \cite{Kakde:2010} and by Hara for the groups of the form $H \times \mathbb{Z}_p$, where $H$ is a finite group of exponent $p$, in \cite{Hara:2010}. The author believes that results in sections \ref{sectionreductions} and \ref{sectioncomputationofk1} will be useful in handling noncommutative main conjectures for other motives. Hence the results in these sections are presented in a self contained manner independent of other sections in the paper.

Thanks to the observation of Burns and Kato mentioned above and Wiles' theorem (the abelian main conjecture), the algebraic results in section \ref{sectioncomputationofk1} reduces the proof of the main conjecture to verifying certain congruences between the Deligne-Ribet, Cassou-Nogues, Barsky $p$-adic $L$-functions (here we are using the hypothesis $\mu = 0$). This is explained in section \ref{sectioncongruences}, proposition \ref{propburnsstrategy} (see also the paper of Burns \cite{Burns:2010}). The required congruences, theorem \ref{theoremrequiredcongruences}, are proven using the Deligne-Ribet $q$-expansion principle. The author learned about the use of $q$-expansion principle from the preprint of Kato \cite{Kato:2006}, however, here we follow the more elementary approach of Ritter-Weiss \cite{RitterWeiss:congruences}. The proof of the congruences is presented in section \ref{sectioncongruences} after a quick introduction to the theory of Hilbert modular forms in subsections \ref{subsectionhmf}-\ref{subsectionqexpansionprinciple} respectively. 

I have accumulated quite a debt of gratitude in writing this paper. Most of all to my teacher Professor John Coates for introducing me to this problem, many invaluable suggestions and discussions and for constant inspiration. To Professor Kazuya Kato for generously sharing his ideas on the main conjecture with me while I was a graduate student in Cambridge. To Professor David Burns for motivating discussions and much needed encouragement towards the end of this paper. I would like to thank Professor Peter Schneider and Professor Otmar Venjakob for carefully reading an earlier version of the manuscript and pointing out several errors. Much of this work was done while I was visiting Newton Institute for the programme on ``Non-Abelian Fundamental Groups in Arithmetic Geometry" and I  thank the organisers,  especially Professor Minhyong Kim, for inviting me and providing a very stimulating environment.

\section{Statement of the main conjecture} \label{sectionstatement}
In this section we formulate the main conjecture in noncommutative Iwasawa theory for totally real number fields following Coates et. al. \cite{CFKSV:2005}.

\subsection{The set up} \label{subsectionsetup} Let $F$ be a totally real number field of degree $r$ over $\mathbb{Q}$. Let $p$ be an odd prime. The field $F(\mu_{p^{\infty}}) = \cup_{ n \geq 1} F(\mu_{p^n})$, where $\mu_{p^n}$ denotes the group of $p^n$th roots of unity, contains a unique Galois extension of $F$ whose Galois group over $F$ is isomorphic to the additive group of $p$-adic integers $\mathbb{Z}_p$. This extension is called the cyclotomic $\mathbb{Z}_p$-extension of $F$ and we denote it by $F^{cyc}$. We write $\Gamma$ for the Galois group of $F^{cyc}$ over $F$.

\begin{definition} An admissible $p$-adic Lie extension $F_{\infty}$ of $F$ is a Galois extension $F_{\infty}$ of $F$ such that (i)$F_{\infty}/F$ is unramified outside a finite set of primes of $F$, (ii) $F_{\infty}$ is totally real, (iii) $\mathcal{G} := Gal(F_{\infty}/F)$ is a $p$-adic Lie group, and (iv) $F_{\infty}$ contains $F^{cyc}$.
\label{defnadmissible}
\end{definition}

From now on, $F_{\infty}/F$ will denote an admissible $p$-adic Lie extension, and we put 
\[
\mathcal{G} = Gal(F_{\infty}/F), \qquad H=Gal(F_{\infty}/F^{cyc}), \qquad \Gamma=\mathcal{G}/H.
\]
Let $\Sigma$ denote a finite set of finite primes of $F$ which will always be assumed to contain the primes of $F$ which ramify in $F_{\infty}$. For an extension $L$ of $F$ we denote the set of primes in $L$ lying above $\Sigma$ by $\Sigma_L$.

For any profinite group $P$, and the ring of integers $O$ in a finite extension of $\mathbb{Q}_p$, we define the Iwasawa algebra
\[
\Lambda_O(P) := \ilim{U} O[P/U],
\]
where $U$ runs through the open normal subgroups of $P$, and $O[P/U]$ denotes the group ring of $P/U$ with coefficients in $O$. When $O = \mathbb{Z}_p$, we write simply $\Lambda(P)$. Unless stated otherwise, we shall consider left modules over the Iwasawa algebras $\Lambda_O(P)$. Following Coates et. al. \cite{CFKSV:2005}, we define 
\[
S(\mathcal{G}, H) = \{s \in \Lambda_O(\mathcal{G}) : \Lambda_O(\mathcal{G})/\Lambda_O(\mathcal{G})s \text{ is a finitely generated } \Lambda_O(H)-\text{module}\}.
\]
We will usually write just $S$ for $S(\mathcal{G},H)$ since $\mathcal{G}$ and $H$ should be clear from context. It is proven in \emph{loc. cit.}, that $S$ is multiplicatively closed subset of nonzero divisors in $\Lambda_O(\mathcal{G})$, which is a left and right Ore set. Hence the localisation $\Lambda_O(\mathcal{G})_S$ of $\Lambda_O(\mathcal{G})$ exists, and the natural map from $\Lambda_O(\mathcal{G})$ to $\Lambda_O(\mathcal{G})_S$ is injective. A $\Lambda_O(\mathcal{G})$-module $M$ is called $S$-torsion if every element of $M$ is annihilated by some element in $S$. It is proven in \emph{loc. cit.} that a $\Lambda_O(\mathcal{G})$-module $M$ is $S$-torsion if and only if it is finitely generated as a $\Lambda_O(H)$-module.

\subsection{$K$-theory} In this subsection we recall a few notions from algebraic $K$-theory. Most of this part is based on the introductory section of Fukaya-Kato \cite{FukayaKato:2006}.

\begin{definition}For any ring $\Lambda$, the group $K_0(\Lambda)$ is an abelian group, whose group law we denote additively, defined by the following set of generators and relations. 
Generators : $[P]$, where $P$ is a finitely generated projective $\Lambda$-module. 
Relations: ($i$) $[P] = [Q]$ if $P$ is isomorphic to $Q$ as a $\Lambda$-module, and ($ii$) $[P\oplus Q] = [P] + [Q]$.
\label{k0}
\end{definition}

\noindent It is easily seen that any element of $K_0(\Lambda)$ can be written as $[P]-[Q]$ for finitely generated projective $\Lambda$-modules $P$ and $Q$. Moreover, $[P]-[Q] = [P^{\prime}] - [Q^{\prime}]$ in $K_0(\Lambda)$ if and only if there is a finitely generated projective $\Lambda$-module $R$ such that $P\oplus Q^{\prime} \oplus R$ is isomorphic to $P^{\prime} \oplus Q \oplus R$. 

\begin{definition} For any ring $\Lambda$, the group $K_1(\Lambda)$ is an abelian group, whose group law we denote multiplicatively, defined by the following generators and relations. 
Generators : $[P, \alpha]$, where $P$ is a finitely generated projective $\Lambda$-module and $\alpha$ is an automorphism of $P$. 
Relations : ($i$) $[P, \alpha] = [Q, \beta]$ if there is an isomorphism $f$ from $P$ to $Q$ such that $f \circ \alpha = \beta \circ f$, ($ii$) $[P, \alpha \circ \beta] = [P, \alpha] [P, \beta]$, and ($iii$) $[P\oplus Q, \alpha \oplus \beta] = [P, \alpha] [Q, \beta]$.
\label{k1}
\end{definition}

Here is an alternate description of $K_1(\Lambda)$. We have a canonical homomorphism $GL_n(\Lambda) \rightarrow K_1(\Lambda)$ defined by mapping $\alpha$ in $GL_n(\Lambda)$ to $[\Lambda^n, \alpha]$, where $\Lambda^n$ is regarded as a set of row vectors and $\alpha$ acts on them from the right. Now using the inclusion maps $GL_n(\Lambda) \hookrightarrow GL_{n+1}(\Lambda)$ given by $g \mapsto \left( \begin{array}{cc} g & 0 \\ 0 & 1\end{array} \right)$, we let 
\[
GL(\Lambda) = \cup_{n \geq 1} GL_n(\Lambda).
\]
Then the canonical homomorphisms $GL_n(\Lambda) \rightarrow K_1(\Lambda)$ induce an isomorphism (see for example Oliver \cite{Oliver:1988}, chapter 1)
\[
\frac{GL(\Lambda)}{[GL(\Lambda), GL(\Lambda)]} \xrightarrow{\sim} K_1(\Lambda),
\]
where $[GL(\Lambda), GL(\Lambda)]$ is the commutator subgroup of $GL(\Lambda)$. If $\Lambda$ is commutative, then the determinant maps, $GL_n(\Lambda) \rightarrow \Lambda^{\times}$, induce the determinant map
\[
det: K_1(\Lambda) \rightarrow \Lambda^{\times}, 
\]
via the above isomorphism. This gives a splitting of the canonical homomorphism $\Lambda^{\times} = GL_1(\Lambda) \rightarrow K_1(\Lambda)$. If $\Lambda$ is semilocal then Vaserstein (\cite{Vaserstein:1969} and \cite{Vaserstein:2005}) proves that the canonical homomorphism $\Lambda^{\times} =GL_1(\Lambda) \rightarrow K_1(\Lambda)$ is surjective. From these two facts we conclude that if $\Lambda$ is a semilocal commutative ring, then the determinant map induces a group isomorphism between $K_1(\Lambda)$ and $\Lambda^{\times}$.

{\bf Example:} If $P$ is a compact $p$-adic Lie group, then $\Lambda_O(P)$ is a semilocal ring. Hence if $P$ is an abelian compact $p$-adic Lie group then $K_1(\Lambda_O(P)) \cong \Lambda_O(P)^{\times}$.
 
Let $I \subset \Lambda$ be any ideal. Denote the group of invertible matrices in $GL(\Lambda)$ which are congruent to the identity modulo $I$ by $GL(\Lambda,I)$. Denote the smallest normal subgroup of $GL(\Lambda)$ containing all elementary matrices congruent to the identity modulo $I$ by $E(\Lambda, I)$. Finally, set $K_1(\Lambda,I) = GL(\Lambda,I)/E(\Lambda,I)$. A lemma of Whitehead says that $E(\Lambda,I) = [GL(\Lambda), GL(\Lambda,I)]$. Hence $K_1(\Lambda,I)$ is an abelian group. 

Let $f: \Lambda \rightarrow \Lambda'$ be a ring homomorphism. We consider the category $\mathcal{C}_f$ of all triplets $(P,\alpha, Q)$, where $P$ and $Q$ are finitely generated projective $\Lambda$-modules and $\alpha$ is an isomorphism between $\Lambda' \otimes_{\Lambda}P$ and $\Lambda' \otimes_{\Lambda} Q$ as $\Lambda'$-modules. A morphism between $(P, \alpha, Q)$ and $(P', \alpha', Q')$ is a pair of $\Lambda$-module morphisms $g: P \rightarrow P'$ and $h: Q \rightarrow Q'$ such that 
\[
\alpha'\circ (Id_{\Lambda'}\otimes g) = (Id_{\Lambda'}\otimes h)\circ \alpha .
\]
It is an isomorphism if both $g$ and $h$ are isomorphisms. A sequence of maps
\[
0 \rightarrow (P', \alpha', Q') \rightarrow (P, \alpha, Q) \rightarrow (P'', \alpha'', Q'') \rightarrow 0,
\]
is a short exact sequence if the underlying sequences
\[
0 \rightarrow P' \rightarrow P \rightarrow P'' \rightarrow 0, \ \text{and} 
\]
\[
0 \rightarrow Q' \rightarrow Q \rightarrow Q'' \rightarrow 0,
\]
are short exact sequences.

\begin{definition} For any ring homomorphism $f: \Lambda \rightarrow \Lambda'$, the group $K_0(f)$ is an abelian group, whose group law we denote additively, defined by the following generators and relations. Generators : $[(P, \alpha, Q)]$, where $(P, \alpha, Q)$ is an object in $\mathcal{C}_f$. Relations : ($i$) $[(P, \alpha, Q)] = [(P', \alpha', Q')]$ if $(P, \alpha, Q)$ is isomorphic to $(P', \alpha', Q')$, ($ii$) $[(P, \alpha, Q)] = [(P', \alpha', Q')] + [(P'', \alpha'', Q'')]$ for every short exact sequence as above and ($iii$) $[(P_1, \beta\circ\alpha, P_3)] = [(P_1, \alpha, P_2)] + [(P_2, \beta, P_3)]$.
\label{relativek0}
\end{definition} 

\noindent For the canonical injection, say $i$, of $\Lambda(\mathcal{G})$ in $\Lambda(\mathcal{G})_S$, we write $K_0(\Lambda(\mathcal{G}), \Lambda(\mathcal{G})_S)$ for $K_0(i)$ and call it the \emph{relative} $K_0$. We give two other descriptions of this group $K_0(\Lambda(\mathcal{G}), \Lambda(\mathcal{G})_S)$. For details see Weibel \cite{Weibel:2000}. Let $\mathcal{C}_S$ be the category of bounded complexes of finitely generated projective $\Lambda(\mathcal{G})$-modules whose cohomologies are $S$-torsion. Consider the abelian group $K_0(\mathcal{C}_S)$ with the following set of generators and relations. Generators : $[C]$, where $C$ is an object of $\mathcal{C}_S$. Relations : ($i$) $[C]=0$ if $C$ is acyclic, and ($ii$) $[C] = [C'] + [C'']$, for every short exact sequence 
\[
0 \rightarrow C' \rightarrow C \rightarrow C'' \rightarrow 0, 
\]
in $\mathcal{C}_S$. \\

\noindent Let $\mathcal{H}_S$ be the category of finitely generated $\Lambda(\mathcal{G})$-modules which are $S$-torsion and which have a finite resolution by finitely generated projective $\Lambda(\mathcal{G})$-modules. Let $K_0(\mathcal{H}_S)$ be the abelian group defined by the following set of generators and relations. Generators : $[M]$, where $M$ is an object of $\mathcal{H}_S$. Relations : ($i$) $[M] =[M']$ if $M$ is isomorphic to $M'$, and ($ii$) $[M] =[M'] + [M'']$ for every short exact sequence 
\[
0 \rightarrow M' \rightarrow M \rightarrow M'' \rightarrow 0, 
\]
of modules in $\mathcal{H}_S$. There are isomorphisms
\[
K_0(\Lambda(\mathcal{G}), \Lambda(\mathcal{G})_S) \xrightarrow{\sim} K_0(\mathcal{C}_S), \ \text{and}
\]
\[
K_0(\Lambda(\mathcal{G}), \Lambda(\mathcal{G})_S) \xrightarrow{\sim} K_0(\mathcal{H}_S),
\]
given as follows. First observe that every isomorphism $\alpha$ from $\Lambda(\mathcal{G})_S\otimes_{\Lambda(\mathcal{G})}P$ to $\Lambda(\mathcal{G})_S\otimes_{\Lambda(\mathcal{G})}Q$ is of the form $s^{-1}a$ with $a$ a $\Lambda(\mathcal{G})$-homomorphism from $P$ to $Q$ and $s$ an element of $S$. Then the above mentioned isomorphisms are respectively given by
\[
[(P, \alpha, Q)] \mapsto [P \xrightarrow{a} Q] - [Q \xrightarrow{s} Q], \ \text{and}
\]
\[
[(P, \alpha, Q)] \mapsto [Q/a(P)] - [Q/Qs].
\]
The relative $K_0$ fits into the localisation exact sequence of $K$-theory
\[
K_1(\Lambda(\mathcal{G})) \rightarrow K_1(\Lambda(\mathcal{G})_S) \xrightarrow{\partial} K_0(\Lambda(\mathcal{G}), \Lambda(\mathcal{G})_S) \xrightarrow{\eta} K_0(\Lambda(\mathcal{G})) \rightarrow K_0(\Lambda(\mathcal{G})_S).
\]
The homomorphism $\partial$ maps $\alpha$ in $K_1(\Lambda(\mathcal{G})_S)$ to $[(\Lambda(\mathcal{G})^n, \tilde{\alpha}, \Lambda(\mathcal{G})^n)]$, where $\tilde{\alpha}$ is any lift of $\alpha$ in $GL_{\infty}(\Lambda(\mathcal{G})_S)$ and $\tilde{\alpha}$ lies in $GL_n(\Lambda(\mathcal{G})_S)$. The homomorphism $\eta$ maps $[(P, \alpha, Q)]$ to $[P] -[Q]$ in $K_0(\Lambda(\mathcal{G}))$. The following lemma is essentially proven in Coates et.al. \cite{CFKSV:2005} under the assumption that $\mathcal{G}$ has no element of order $p$. The same technique gives this more general result as we now show.

\begin{lemma} The homomorphism $\partial$ is surjective.
\label{localisationsequence}
\end{lemma}
\noindent{\bf Proof:} We will show that the homomorphism $\eta$ is 0. Fix a pro-$p$ open normal subgroup $P$ of $\mathcal{G}$, and put $\Delta = \mathcal{G}/P$. We write $\mathcal{V} = \mathcal{V}(\Delta)$ for the set of irreducible representations of the finite group $\Delta$ over $\bar{\mathbb{Q}}_p$ and we take $L$ to be some fixed finite extension of $\mathbb{Q}_p$ such that all representations in $\mathcal{V}$ can be realised over $L$. Recall the construction in Coates et.al. \cite{CFKSV:2005} of the canonical homomorphism $\lambda$ from $K_0(\Lambda(\mathcal{G}))$ to $\prod_{\rho \in \mathcal{V}} K_0(L)$. It is the composition $\lambda = \lambda_4 \circ \lambda_3 \circ \lambda_2 \circ \lambda_1$ of four natural maps $\lambda_i$ ($i = 1,2,3,4$) as follows
\[
\lambda_1 : K_0(\Lambda(\mathcal{G})) \rightarrow K_0(\mathbb{Z}_p[\Delta]),
\]
\[
\lambda_2 : K_0(\mathbb{Z}_p[\Delta]) \rightarrow K_0(\mathbb{Q}_p[\Delta]),
\]
\[
\lambda_3 : K_0(\mathbb{Q}_p[\Delta]) \rightarrow K_0(L[\Delta]),
\]
\[
\lambda_4 : K_0(L[\Delta]) \xrightarrow{\sim} \prod_{\rho \in \mathcal{V}}K_0(M_{n_{\rho}}(L)) \xrightarrow{\sim} \prod_{\rho \in \mathcal{V}} K_0(L).
\]
It can be shown that $\lambda_1$, $\lambda_2$, and $\lambda_3$ are injective. Hence $\lambda$ is injective as well. We now recall an alternate description of $\lambda$ given in \emph{loc. cit.}. The augmentation map from $\Lambda_{\mathcal{O}}(\mathcal{G})$ to $\mathcal{O}$ induces a map from $K_0(\Lambda_{\mathcal{O}}(\mathcal{G}))$ to $K_0(\mathcal{O})$ which we denote by $\tau$. Let $U$ be a finitely generated $\Lambda_{\mathcal{O}}(\mathcal{G})$-module having a finite resolution by finitely generated projective $\Lambda_{\mathcal{O}}(\mathcal{G})$-modules. Then one can define the class of $U$, denoted by $[U]$, in $K_0(\Lambda_{\mathcal{O}}(G))$. Since $\mathcal{O}$ is a domain, $K_0(\mathcal{O})$ can be identified with the Grothendieck group of the category of all finitely generated $\mathcal{O}$-modules. Then $\tau$ is explicitly given by 
\[
\tau([U]) = \sum_{ i \geq 0} (-1)^i [H_i(\mathcal{G}, U)].
\]
$\tau$ factors through the map
\[
\epsilon : K_0(\Lambda_{\mathcal{O}}(\mathcal{G})) \rightarrow K_0(\Lambda_{\mathcal{O}}(\Gamma)),
\]
given by the natural surjection of $\Lambda_{\mathcal{O}}(\mathcal{G})$ on $\Lambda_{\mathcal{O}}(\Gamma)$. Moreover, $\epsilon$ is given explicitly by 
\[
\epsilon([U]) = \sum_{ i \geq 0} (-1)^i[H_i(H, U)].
\]
Now take $\mathcal{O}$ to be the ring of integers of $L$ and $j$ be the isomorphism from $K_0(\mathcal{O})$ to $K_0(L)$ induced by the inclusion of $\mathcal{O}$ in $L$. Let $tw_{\rho}(U)$ be defined by $U \otimes_{\mathbb{Z}_p} \mathcal{O}^{n_{\rho}}$, for any $\rho$ in $\mathcal{V}$. It can be made into a left $\mathcal{G}$ module by the diagonal action. This action extends to make $tw_{\rho}(U)$ a $\Lambda(\mathcal{G})$-module. It is proven in Venjakob \cite{Venjakob:2005} that $tw_{\rho}(U)$ is a finitely generated $S$-torsion $\Lambda(\mathcal{G})$-module if $U$ is. 

\noindent We now finish the proof. Take $[(P, \alpha, Q)]$ in $K_0(\Lambda(\mathcal{G}), \Lambda(\mathcal{G})_S)$. As remarked earlier, $\alpha$ is of the form $s^{-1}a$, with $a$ a $\Lambda(\mathcal{G})$-homomorphism from $P$ to  $Q$ and $s$ is an element of $S$. We will show that $[Q/a(P)]$ is 0 in $K_0(\Lambda(\mathcal{G}))$. Put $M= Q/a(P)$. Then $H_i(H, tw_{\rho}(M))$ is a $\Lambda_{\mathcal{O}}(\Gamma)$-torsion module for all $i \geq 0$ and thus its class in $K_0(\Lambda_{\mathcal{O}}(\Gamma))$ vanishes. Hence it follows that $\epsilon([tw_{\rho}(M)]) = 0$, whence $\tau([tw_{\rho}(M)]) = 0$ for all $\rho$ in $\mathcal{V}$. But this implies that $\lambda([M])=0$ and so $[M]=0$ in $K_0(\Lambda(\mathcal{G}))$. This completes the proof. \qed \\

It will be more convenient to work with a slight modification of the $K_1$ group which we now define. Lets put $\widehat{\Lambda_O(\mathcal{G})_S}$ as the $p$-adic completion of the ring $\Lambda_O(\mathcal{G})_S$. 

\begin{definition} Let $O$ be the ring of integers in a finite extension $L$ of $\mathbb{Q}_p$. Let $P$ be a finite group. We define $SK_1(O[P]) = ker(K_1(O[P]) \rightarrow K_1(L[P]))$. For a pro-finite group $\mathcal{G} = \ilim{U} \mathcal{G}/U$, we define $SK_1(\Lambda_O(\mathcal{G})) := \ilim{U} SK_1(O[\mathcal{G}/U])$. We define $SK_1(\Lambda_O(\mathcal{G})_S)$ and $SK_1(\widehat{\Lambda_O(\mathcal{G})_S})$ to be the image of $SK_1(\Lambda_O(\mathcal{G}))$ in $K_1(\Lambda(\mathcal{G})_S)$ and $K_1(\widehat{\Lambda_O(\mathcal{G})_S})$ respectively under the natural maps
\[
K_1(\Lambda_O(\mathcal{G})) \rightarrow K_1(\Lambda_O(\mathcal{G})_S),
\]
and
\[
K_1(\Lambda_O(\mathcal{G})) \rightarrow K_1(\widehat{\Lambda_O(\mathcal{G})_S}).
\]
\label{defnsk1}
 \end{definition}

\begin{definition} Let $R$ be either $\Lambda_O(\mathcal{G})$ or $\Lambda_O(\mathcal{G})_S$ or $\widehat{\Lambda_O(\mathcal{G})_S}$, then we define
\[
K_1'(R) = K_1(R)/SK_1(R).
\]
\label{defnkprime}
\end{definition}

\subsection{The arithmetic side of the main conjecture} We now explain the basic arithmetic objects attached to $F_{\infty}/F$ needed to formulate the main conjecture.

Let $M_{\infty}$ be the maximal abelian $p$-extension of $F_{\infty}$, which is unramified outside the set of primes above $\Sigma$. As usual $\mathcal{G}$ acts on $Gal(M_{\infty}/F_{\infty})$ by $g\cdot x = \tilde{g}x\tilde{g}^{-1}$, where $g$ is in $\mathcal{G}$, and $\tilde{g}$ is any lifting of $g$ to $Gal(M_{\infty}/F)$. This action extends to a left action of $\Lambda(\mathcal{G})$. 

\begin{definition} We say that $F_{\infty}/F$ satisfies the hypothesis $\mu=0$ if there exists a pro-$p$ open subgroup $H'$ of $H$ such that the Galois group over $F_{\infty}^{H'}$ of the maximal abelian $p$-extension of $F_{\infty}^{H'}$ unramified outside $\Sigma$ is a finitely generated $\mathbb{Z}_p$-module. 
\label{hypothesisonmu}
\end{definition}

\noindent Let $L$ be a finite extension of $F$ such that $L^{cyc} = F_{\infty}^{H'}$. Then the above hypothesis is equivalent to vanishing of the Iwasawa $\mu$-invariant of $L(\mu_p)$. This vanishing of the $\mu$-invariant is a special case of a general conjecture of Iwasawa asserting that, for every finite extension $K$ of $\mathbb{Q}$, the Galois group over $K^{cyc}$ of the maximal unramified abelian $p$-extension of $K^{cyc}$ is a finitely generated $\mathbb{Z}_p$-module. When $K$ is an abelian extension of $\mathbb{Q}$, this conjecture is proven by Ferrero-Washington \cite{FerreroWashington:1979} (A different proof is given by Sinnott \cite{Sinnott:1984}). \\

\begin{lemma} $Gal(M_{\infty}/F_{\infty})$ is finitely generated over $\Lambda(H)$ if and only if $F_{\infty}/F$ satisfies the hypothesis $\mu = 0$.
\label{xistorsion}
\end{lemma}
\noindent{\bf Proof:} Put $X = Gal(M_{\infty}/F_{\infty})$, and let $H'$ be any pro-$p$ open subgroup of $H$. Thus $\Lambda(H')$ is a local ring. It follows from Nakayama's lemma that $X$ is finitely generated over $\Lambda(H)$ if and only if $X_{H'}$ is a finitely generated $\mathbb{Z}_p$-module. Let $F_{\Sigma}$ denote the maximal pro-$p$ extension of $F_{\infty}$ which is unramified outside primes above $\Sigma$. Then $F_{\Sigma}$ is Galois over $F$. Let $K_{\infty} = F_{\infty}^{H'}$. Then we have the inflation-restriction exact sequence

\begin{align*}
0 \rightarrow H^1(H', \mathbb{Q}_p/\mathbb{Z}_p)  \rightarrow  H^1(Gal(F_{\Sigma}/K_{\infty}), \mathbb{Q}_p/\mathbb{Z}_p) & \\ \rightarrow H^1(Gal(F_{\Sigma}/F_{\infty}), \mathbb{Q}_p/\mathbb{Z}_p)^{H'} \rightarrow H^2(H', \mathbb{Q}_p/\mathbb{Z}_p).
\end{align*}
As $H'$ is a $p$-adic Lie group, $H^i(H', \mathbb{Q}_p/\mathbb{Z}_p)$ are cofinitely generated $\mathbb{Z}_p$-modules for all $i \geq 0$. Moreover, since $Gal(F_{\Sigma}/F)$ acts trivially on $\mathbb{Q}_p/\mathbb{Z}_p$, we have 
\[
H^1(Gal(F_{\Sigma}/L), \mathbb{Q}_p/\mathbb{Z}_p) = Hom(Gal(M_L/L), \mathbb{Q}_p/\mathbb{Z}_p),
\]
for every intermediate field $L$ with $F_{\Sigma} \supset L \supset F$; here $M_L$ denotes the maximal abelian $p$-extension of $L$ which is unramified outside $\Sigma_L$. We conclude from the above exact sequence that $X_{H'}$ is a finitely generated $\mathbb{Z}_p$-module if and only if $Gal(M_{K_{\infty}}/K_{\infty})$ is a finitely generated $\mathbb{Z}_p$-module. The conclusion of the lemma is now plain since $Gal(M_{K_{\infty}}/K_{\infty})$ being a finitely generated $\mathbb{Z}_p$-module is precisely the hypothesis $\mu=0$ for some open subgroup $H'$.
\qed 

In the view of the lemma, we shall henceforth make the following \\
{\bf Assumption:} Our admissible $p$-adic Lie extension $F_{\infty}/F$ satisfies the hypothesis $\mu = 0$. 

\begin{remark} When $F_{\infty} = F^{cyc}$, Iwasawa \cite{Iwasawa:1973} proves that $X = Gal(M_{\infty}/F_{\infty})$ is a finitely generated $\Lambda(\Gamma)$-torsion module by showing that the `Leopoldt defect' is bounded in the cyclotomic extension. When $\mathcal{G}$ is a $p$-adic Lie group of dimension greater than 1, Ochi and Venjakob \cite{OchiVenjakob:2003} show that $X$ is finitely generated $\Lambda(\mathcal{G})$-torsion module.  
\end{remark}

The group $X =Gal(M_{\infty}/F_{\infty})$ is a fundamental arithmetic object which is studied through the main conjecture. The above lemma shows that $X$ is a $S$-torsion module, however, if $\mathcal{G}$ has elements of order $p$, then $X$ may not have a finite resolution by finitely generated projective $\Lambda(\mathcal{G})$-modules. In our approach to the proof of the main conjecture it is necessary to consider $\mathcal{G}$ having elements of order $p$. So we consider a complex of $\Lambda(\mathcal{G})$-modules which is closely related to $X$, and is quasi-isomorphic to a bounded complex of finitely generated projective $\Lambda(\mathcal{G})$-modules. We put $C(F_{\infty}/F)$ to be the complex
\[
C(F_{\infty}/F) = RHom(R\Gamma_{\acute{e}t}(Spec(\mathcal{O}_{F_{\infty}}[\frac{1}{\Sigma}]), \mathbb{Q}_p/\mathbb{Z}_p), \mathbb{Q}_p/\mathbb{Z}_p).
\]
Here $\mathbb{Q}_p/\mathbb{Z}_p$ is the constant sheaf corresponding to the abelian group $\mathbb{Q}_p/\mathbb{Z}_p$ on the \'{e}tale site of $Spec(\mathcal{O}_{F_{\infty}}[\frac{1}{\Sigma}])$. Since $\mathbb{Q}_p/\mathbb{Z}_p$ is a direct limit of finite abelian groups of $p$-power order, we have an isomorphism with Galois cohomology
\[
R\Gamma_{\acute{e}t}(Spec(\mathcal{O}_{F_{\infty}}[\frac{1}{\Sigma}]), \mathbb{Q}_p.\mathbb{Z}_p) \xrightarrow{\sim} R\Gamma(Gal(F_{\Sigma}/F_{\infty}), \mathbb{Q}_p/\mathbb{Z}_p).
\]
Here $F_{\Sigma}$ is the maximal $p$-extension of $F_{\infty}$ unramified outside primes above $\Sigma$. Then 
\[
H^i(C(F_{\infty}/F)) = \left\{ \begin{array}{l l}
0 & \text{if $i \neq 0,-1$} \\
\mathbb{Z}_p & \text{if $i=0$} \\
Gal(M_{\infty}/F_{\infty}) & \text{if $i=-1$}.
\end{array} \right.
\]
The following results are proven in Fukaya-Kato \cite{FukayaKato:2006}. \\
($i$) The complex $C(F_{\infty}/F)$ is quasi-isomorphic to a bounded complex of finitely generated projective $\Lambda(\mathcal{G})$-modules. By passing to the derived category we identify $C(F_{\infty}/F)$ with a quasi-isomorphic bounded complex of finitely generated $\Lambda(\mathcal{G})$-modules. \\
($ii$) If $F \subset K \subset F_{\infty}$ is any extension of $F$, then 
\begin{equation}
\Lambda(Gal(K/F)) \otimes^{L}_{\Lambda(\mathcal{G})}C(F_{\infty}/F) \xrightarrow{\sim} C(K/F), 
\label{decent}
\end{equation}
where $\Lambda(\mathcal{G})$ acts on the right on $\Lambda(Gal(K/F))$ through the natural surjection of $\Lambda(\mathcal{G})$ on $\Lambda(Gal(K/F))$ and $C(K/F)$ is the complex
\[
C(K/F) = RHom(R\Gamma_{\acute{e}t}(Spec(O_K[\frac{1}{\Sigma}]), \mathbb{Q}_p/\mathbb{Z}_p),\mathbb{Q}_p/\mathbb{Z}_p).
\]

By lemma \ref{xistorsion}, we know that $C(F_{\infty}/F)$ is $S$-torsion i.e. $\Lambda(\mathcal{G})_S \otimes^{L}_{\Lambda(\mathcal{G})} C(F_{\infty}/F)$ is acyclic. Hence we can talk about the class of $C(F_{\infty}/F)$, denoted by $[C(F_{\infty}/F)]$, in the group $K_0(\Lambda(\mathcal{G}), \Lambda(\mathcal{G})_S)$. \\

\subsection{The analytic side of the main conjecture} We now explain the analytic objects in the Iwasawa theory of $F_{\infty}/F$ needed for formulating the main conjecture. 

Let $\rho$ be an Artin representation (i.e. kernel of $\rho$ is open) of $Gal(\overline{F}/F)$ on a finite dimensional vector space over $\overline{\mathbb{Q}}_p$, factoring through $\mathcal{G}$. Let $\alpha$ be an embedding of $\overline{\mathbb{Q}}_p$ in $\mathbb{C}$. The resulting complex Artin representation $\alpha \circ \rho$ gives the complex Artin $L$-function $L(\alpha \circ \rho, s)$. A famous result of Klingen \cite{Klingen:1962} and Siegel \cite{Seigel:1970} says that $L(\alpha \circ \rho, -n)$ is an algebraic number for any odd positive integer $n$. This number depends on the choice of $\alpha$. If $\alpha$ is replaced by another embedding then we get a conjugate of the algebraic number $L(\alpha \circ \rho, -n)$ over $\mathbb{Q}$. However, the beauty of the result of Klingen and Siegel is that it makes it possible to choose a canonical conjugate. Hence we obtain an algebraic number $L(\rho, -n)$, \emph{``the value of complex $L$-function associated to $\rho$ at $-n$"}. For details see section 1.2 in Coates-Lichtenbaum \cite{CoatesLichtenbaum:1973} (when $\rho$ is one dimensional see equation (\ref{lvalue}) below in section \ref{subsectionlvalues}. For $\rho$ of higher dimension one has to use Brauer's induction theorem). Similarly, we can define the value of complex $L$-function associated to $\rho$ at $-n$ with Euler factors at primes in $\Sigma$ removed. We denote it by $L_{\Sigma}(\rho, -n)$. This value obviously depends only on the character associated to $\rho$. We use the same letter to denote Artin representation and its character.  

Let $L$ be a finite extension of $\mathbb{Q}_p$ with ring of integers $O$. Let $\rho$ be a continuous homomorphism from $\mathcal{G}$ into $GL_n(O)$. It induces a ring homomorphism from $\Lambda(\mathcal{G})$ into $M_n(\Lambda_{O}(\Gamma))$. This homomorphism is given on elements of $\mathcal{G}$ by mapping $\sigma$ in $\mathcal{G}$ to $\rho(\sigma)\overline{\sigma}$, where $\overline{\sigma}$ is the image of $\sigma$ in $\Gamma$. We write $Q_{O}(\Gamma)$ for the field of fraction of $\Lambda_{O}(\Gamma)$. It is proven in \cite{CFKSV:2005} that this homomorphism extends to a homomorphism 
\[
\Phi_{\rho} : \Lambda(\mathcal{G})_S \rightarrow M_n(Q_{O}(\Gamma)).
\]
$\Phi_{\rho}$ induces a homomorphism 
\[
\Phi'_{\rho} : K_1'(\Lambda(\mathcal{G})_S) \rightarrow K_1(M_n(Q_{O}(\Gamma))) = Q_{O}(\Gamma)^{\times}.
\]
Now let $\varphi$ be the augmentation map from $\Lambda_{O}(\Gamma)$ to $O$. We denote its kernel by $\mathfrak{p}$. If we write $\Lambda_{O}(\Gamma)_{\mathfrak{p}}$ for the localisation of $\Lambda_{O}(\Gamma)$ at the prime ideal $\mathfrak{p}$, then $\varphi$ extends to a homomorphism 
\[
\varphi : \Lambda_{O}(\Gamma)_{\mathfrak{p}} \rightarrow L.
\]
We extend this map to a map $\varphi'$ from $Q_{O}(\Gamma)$ to $L \cup \{\infty\}$ by mapping any $x$ in $Q_{O}(\Gamma) - \Lambda_{O}(\Gamma)_{\mathfrak{p}}$ to $\infty$. The composition of $\Phi'_{\rho}$ with $\varphi'$ gives a map
\begin{center}
$K_1'(\Lambda(\mathcal{G})_S) \rightarrow L \cup \{\infty\}$ \\
$ x \mapsto x(\rho) $.
\end{center}
In the classical Iwasawa theory $x(\rho)$ would be written as $ \int_{\mathcal{G}} \rho dx$. \\
This map has the following properties: \\
($i$) Let $\mathcal{G}'$ be an open subgroup of $\mathcal{G}$. Let $\chi$ be an one dimensional representation of $\mathcal{G}'$ and $\rho = Ind_{\mathcal{G}'}^{\mathcal{G}}(\chi)$. If $N$ is the norm map from $K_1'(\Lambda(\mathcal{G})_S)$ to $K_1'(\Lambda(\mathcal{G}')_S)$, then for any $x$ in $K_1'(\Lambda(\mathcal{G})_S)$, we have $x(\rho) = (Nx)(\chi)$. Note that $\Lambda(\mathcal{G})_S$ is a free $\Lambda(\mathcal{G}')_S$-module of finite rank as proven in \cite{SchneiderVenjakob:2010a}, proposition 4.5 (i). \\
($ii$) Let $\rho_i$ be continuous homomorphisms from $\mathcal{G}$ into $GL_{n_i}(L)$, for $i=1,2$. We then get a continuous homomorphism $\rho_1\oplus \rho_2$ from $\mathcal{G}$ into $GL_{n_1+n_2}(L)$. Then for any $x$ in $K_1'(\Lambda(\mathcal{G})_S)$, $x(\rho_1\oplus \rho_2) = x(\rho_1)x(\rho_2)$. \\
($iii$) Let $U$ be a subgroup of $H$, normal in $\mathcal{G}$. Let $\pi$ be the homomorphism from $K_1'(\Lambda(\mathcal{G})_S)$ to $K_1'(\Lambda(\mathcal{G}/U)_S)$ induced by the natural surjection of $\Lambda(\mathcal{G})_S$ onto $\Lambda(\mathcal{G}/U)_S$. Let $\rho$ be a continuous homomorphism of $\mathcal{G}/U$ into $GL_n(L)$. We write $inf(\rho)$ for the composition of the natural surjection from $\mathcal{G}$ onto $\mathcal{G}/U$ with $\rho$. Then for any $x$ in $K_1(\Lambda(\mathcal{G})_S)$, we have $x(inf(\rho)) = \pi(x)(\rho)$. \\

\subsection{The Statement} 
{\bf Notation:} Let $\kappa_F$ be the $p$-adic cyclotomic character 
\[
Gal(F(\mu_{p^{\infty}})/F) \rightarrow \mathbb{Z}_p^{\times},
\]
defined by 
\[
\sigma(\zeta) = \zeta^{\kappa_{F}(\sigma)},
\]
for any $\sigma$ in $Gal(F(\mu_{p^{\infty}})/F)$ and any $p$-power root of unity $\zeta$. \\

\noindent We are now ready to state the main conjecture. Recall the localisation sequence of $K$-theory
\[
K_1'(\Lambda(\mathcal{G})) \rightarrow K_1'(\Lambda(\mathcal{G})_S) \xrightarrow{\partial} K_0(\Lambda(\mathcal{G}), \Lambda(\mathcal{G})_S) \rightarrow 0.
\]

\begin{theorem} ({\bf Main Conjecture}) Let $F_{\infty}/F$ be an admissible $p$-adic Lie extension satisfying the hypothesis $\mu=0$. Then there is a unique $\zeta(F_{\infty}/F)$ in $K_1'(\Lambda(\mathcal{G})_S)$ such that $\partial(\zeta(F_{\infty}/F)) = -[C(F_{\infty}/F)]$ and for any Artin character $\rho$ of $\mathcal{G}$ and any positive integer $r$ divisible by $[F_{\infty}(\mu_p):F_{\infty}]$,
\[
\zeta(F_{\infty}/F)(\rho\kappa_F^r) = L_{\Sigma}(\rho, 1-r).
\]
\label{mainconjecture}
\end{theorem}

\begin{remark} The Main Conjecture in this form was first formulated by Kato \cite{Kato:1993}, and Fontaine and Perrin-Riou \cite{FontaineRiou:1991} in the case when $\mathcal{G}$ is abelian. This was generalised to include noncommutative groups $\mathcal{G}$ by Burns and Flach \cite{BurnsFlach:2001}, Huber and Kings \cite{HuberKings:2002}, Coates et.al. \cite{CFKSV:2005}, and Fukaya and Kato \cite{FukayaKato:2006}. The case of one dimensional $p$-adic Lie groups was considered by Ritter and Weiss in \cite{RitterWeiss:2} and proved very recently in \cite{RitterWeiss:2010}. Using their result D. Burns \cite{Burns:2010} proved the general case. We must warn that in all above sources the $p$-adic zeta function is conjectured to live in $K_1$ as opposed to $K_1'$ and hence uniqueness may not hold. However, the two versions are equivalent because $SK_1(\Lambda(\mathcal{G})_S)= Ker(\partial) \cap (\cap_{\rho} Ker(\varphi' \circ \Phi_{\rho}'))$, where $\rho$ runs through all Artin representations of $\mathcal{G}$, as proven by Burns \cite{Burns:2010}.
\end{remark}

\begin{remark} $\zeta(F_{\infty}/F)$ is called a $p$-adic zeta function for the extension $F_{\infty}/F$. It depends on $\Sigma$ but we suppress this fact in the notation. If $\mathcal{G}$ is abelian then existence and uniqueness of the $p$-adic zeta function is well known, as we will soon see.
\end{remark}

\begin{remark} One can show that the validity of the Main Conjecture is independent of $\Sigma$ as long as it contains all primes of $F$ which ramify in $F_{\infty}$. 
\end{remark}

\noindent{\bf Example:} Let $F= \mathbb{Q}(\mu_p)^+$. We take $F_{\infty}$ to be the maximal abelian $p$-extension of $F^{cyc}$ unramified outside the unique prime above $p$. Then $F_{\infty}$ is Galois over $F$. If $p$ is an irregular prime then dimension of $\mathcal{G} = Gal(F_{\infty}/F)$ is more than 1. Our main theorem applies unconditionally to $F_{\infty}/F$ since validity of the hypothesis $\mu=0$ for $F_{\infty}/F$ is a theorem of Ferrero-Washington \cite{FerreroWashington:1979}. \\

\noindent{\bf Example:} Ravi Ramakrishna, in corollary 1 of \cite{Ramakrishna:2002}, constructs infinitely many surjective representation $G_{\mathbb{Q}} \rightarrow SL_2(\mathbb{Z}_7)$ each of which is unramified outside a finite set of primes. If we take any one of this representation, say $\rho : G_{\mathbb{Q}} \rightarrow SL_2(\mathbb{Z}_7)$, then $\overline{\mathbb{Q}}^{ker(\rho)}$ is either a CM field or a totally real field. Let $F_{\infty}$ be the compositum of $\mathbb{Q}^{cyc}$, the cyclotomic $\mathbb{Z}_7$-extension, and the maximal totally real subfield of $\overline{\mathbb{Q}}^{ker(\rho)}$. Then $F_{\infty}/\mathbb{Q}$ is an admissible $7$-adic Lie extension. The author learned about this example from H. Hida.

\section{Classical Iwasawa main conjecture} \label{sectionclassicalmainconjecture}
 In this subsection we recall the classical main conjecture (formulated for $F=\mathbb{Q}$ by Iwasawa and for totally real number fields by Coates \cite{Coates:1977} and Greenberg \cite{Greenberg:1977}) proven by Mazur-Wiles \cite{MazurWiles:1984} in the case $F=\mathbb{Q}$ (and also by Rubin by a different method) and in general by Wiles \cite{Wiles:1990}.

Let $\psi$ be an Artin representation of $Gal(\overline{F}/F)$. Let $O_{\psi}$ be the ring obtained by adjoining values of $\psi$ to $\mathbb{Z}_p$. Let $F_{\psi}$ be the field defined by $\psi$ i.e. the fixed field of kernel of $\psi$. Assume $F_{\psi} \cap F^{cyc} = F$. Recall that we have fixed a finite set $\Sigma$ of finite primes of $F$ and we assume that it contains all primes which ramify in $F_{\psi}F^{cyc}/F$. Let $\gamma$ be a fixed topological generator of $\Gamma = Gal(F^{cyc}/F)$ and let $u \in \mathbb{Z}_p^{\times}$ be the image of $\gamma$ under the $p$-adic cyclotomic character. If $\psi$ is one dimensional then Deligne-Ribet \cite{DeligneRibet:1980} and Cassou-Nogu\'{e}s \cite{Cassou:1979} and Barsky \cite{Barsky:1978} proved that there is a power $G_{\psi}(T) \in O_{\psi}[[T]]$ such that for any non-trivial $\psi$ and any positive integer $n$ divisible by $[F_{\psi}F^{cyc}(\mu_p):F_{\psi}F^{cyc}]$
\[
G_{\psi}(u^n-1) = L_{\Sigma}(\psi, 1-n),
\]
and if $\psi$ is the trivial character $\mathds{1}$ and $n$ is a positive integer divisible by $[F^{cyc}(\mu_p):F^{cyc}]$,
\[
\frac{G_{\mathds{1}}(u^n-1)}{u^n-1} = L_{\Sigma}(\mathds{1}, 1-n).
\]

Now let $M_{\psi}$ be the maximal abelian $p$-extension of $F_{\psi}F^{cyc}$ unramified outside prime above $\Sigma$. Let $X_{\psi}$ be the Galois group of $M_{\psi}$ over $F_{\psi}F^{cyc}$. Let $h_{\psi}(T)$ be the characteristic polynomial of $\gamma$ acting on the $e_{\psi}(X_{\psi} \otimes_{\mathbb{Z}_p} \overline{\mathbb{Q}}_p)$, where $e_{\psi}$ is the idempotent in $\overline{\mathbb{Q}}_p[H' \times H_p]$ corresponding to $\psi$. We write $G_{\psi}(T)$ as $\pi^{\mu_{\psi}}G_{\psi}^*(T)$, with $\pi$ a uniformiser of $O_{\psi}$ and $G_{\psi}^*(T)$ not divisible by $\pi$. Then Wiles shows that $h_{\psi}(T)$ and $G_{\psi}^*(T)$ generate the same ideal in $O_{\psi}[[T]]$ (theorem 1.3 \cite{Wiles:1990}). Wiles also shows that if $p$ does not divide the order of $Gal(F_{\psi}/F)$, then in fact $\mu_{\psi} = 0$ (theorem 1.4 \cite{Wiles:1990}. Note that we are using the assumption $\mu=0$). 

Now let us assume that $\mathcal{G} = H \times \Gamma$, where $H$ is a finite abelian group. Serre's account \cite{Serre:1978} of the work of Deligne-Ribet shows that there is a pseudomeasure $\zeta \in Q(\Lambda(\mathcal{G}))$ such that for any nontrivial one dimensional character $\psi$ of $H$, we get $\psi(\zeta) = G_{\psi}(T)$ and $\mathds{1}(\zeta)=G_{\mathds{1}}(T)/T$. Here for a one dimensional character $\psi$ of $H$, we denote the induced map
\[
Q(\Lambda(\mathcal{G})) \rightarrow Q(\Lambda_{O_{\psi}}(\Gamma)) \cong Q(O_{\psi}[[T]]),
\]
by the same symbol $\psi$. The last isomorphism is the one obtained by mapping $\gamma$ to $1+T$.

\begin{lemma} The $p$-adic zeta function $\zeta$ is a unit in $\Lambda(\mathcal{G})_S$.
\end{lemma}

\noindent{\bf Proof:} We write $H$ as $H' \times H_p$, where $H_p$ is the $p$-Sylow subgroup of $H$ and $H'$ is a finite group whose order is prime to $p$. Let $\widehat{H'}$ be the set of all irreducible characters if $H'$. Let $O_{\psi}$ be the finite extension of $\mathbb{Z}_p$ obtained by adjoining all values of $\psi$ and let $\pi_{\psi}$ be a uniformiser of $O_{\psi}$. We have the following decomposition
\[
\Lambda(\mathcal{G})_S \xrightarrow{\sim} \oplus_{\psi \in \hat{H'}} \Lambda_{O_{\psi}}(H_p \times \Gamma)_{S},
\]
which maps $g = (h',\overline{g}) \in \mathcal{G}$, with $h' \in H'$ and $\overline{g} \in H_p \times \Gamma$,  to $(\psi(h'))\overline{g})_{\psi \in \hat{H'}}$. Let $\mathds{1}$ be the trivial character of $H_p$. Each of the summands $\Lambda_{O_{\psi}}(H_p \times \Gamma)_S$ is a local ring with maximal ideal 
\[
\mathfrak{m_{\psi}} = \{x \in \Lambda_{O_{\psi}}(H_p \times \Gamma)_S : \mathds{1}(x) \in \pi_{\psi}\Lambda_{O_{\psi}}(\Gamma)_S\}.
\]
For each $\psi \in \hat{H'}$, the element $(\psi \times \mathds{1})(\zeta) \in \Lambda_{O_{\psi}}(\Gamma)_S$ does not belong to $\pi_{\psi}\Lambda_{O_{\psi}}(\Gamma)_S$ because of our assumption $\mu = 0$ and the result of Wiles mentioned above. Hence $\zeta$ is a unit in $\Lambda(\mathcal{G})_S$. 
\qed

\begin{theorem} Under the boundary map $K_1(\Lambda(\mathcal{G})_S) \xrightarrow{\partial} K_0(\Lambda(\mathcal{G}), \Lambda(\mathcal{G})_S)$,
\[
\zeta \mapsto -[C(F_{\infty}/F)].
\]
Note that $K_1'(\Lambda(\mathcal{G})_S) = K_1(\Lambda(\mathcal{G})_S)$ as $SK_1(\Lambda(\mathcal{G})) = \{1\}$.
\end{theorem}

\noindent{\bf Proof:} Put $C=C(F_{\infty}/F)$. Let $\hat{H}$ be the set of all irreducible character of $H$. Consider the following diagram with exact rows
\[
\xymatrix{1  \ar[r] & \Lambda(\mathcal{G})^{\times} \ar[d]_{\theta} \ar[r] & \Lambda(\mathcal{G})_S^{\times} \ar[d]_{\theta_S} \ar[r]^(.35){\partial} & K_0(\Lambda(\mathcal{G}), \Lambda(\mathcal{G})_S) \ar[r] \ar[d]_{\theta_0} & 1 \\
1 \ar[r] & \prod \Lambda_{O_{\psi}}(\Gamma)^{\times} \ar[r] & \prod \Lambda_{O_{\psi}}(\Gamma)_S^{\times} \ar[r]^(.35){\partial} & \prod K_0(\Lambda_{O_{\psi}}(\Gamma), \Lambda_{O_{\psi}}(\Gamma)_S) \ar[r] & 1,}
\]
where the product runs through all $\psi \in \hat{H}$. The maps $\theta$ and $\theta_S$ are injective. The map $\theta_0$ is
\[
\theta_0(P^{\cdot}) = (\Lambda_{O_{\psi}}(\Gamma) \otimes_{\Lambda(\mathcal{G})}^{L} P^{\cdot}) ,
\]
for any complex $P^{\cdot}$ of $\Lambda(\mathcal{G})$ modules with $S$-torsion cohomologies. Let $\theta_0(C(F_{\infty}/F)) = (C^{\cdot}_{\psi})$. The cohomologies of the complex $C^{\cdot}_{\psi}$ are 
\[
H^i(C^{\cdot}_{\psi}) = \left\{ 
\begin{array}{l l}
0 & \text{if $i \neq 0,-1$} \\
M_{\psi} & \text{if $i=-1$} \\
\mathbb{Z}_p & \text{if $i=0$ and $\psi =\mathds{1}$} \\
0 & \text{if $i=0$ and $\psi \neq \mathds{1}$} 
\end{array} \right.
\] 
Hence the classical main conjecture proven by Wiles says that $\partial(\psi(\zeta)) = -[C^{\cdot}_{\psi}]$ for all $\psi \in \hat{H}$. Since $\theta_0(\partial(\zeta)) = \partial(\theta_S(\zeta))$, the result will follow if we show that $\theta_0$ is injective. The injectivity of $\theta_0$ will follow from snake lemma if we can show that the cokernel of $\theta$ injects in cokernel of $\theta_S$. Let $x \in \Lambda(\mathcal{G})_S^{\times}$ such that $\theta_S(x) = (x_{\psi}) \in \prod \Lambda_{O_{\psi}}(\Gamma)^{\times}$. Then 
\[
x = \frac{1}{|H|} \sum_{h \in H} h \Big( \sum_{\psi\in \hat{H}} x_{\psi} \psi(h^{-1}) \Big).
\]
Hence $x \in \Lambda(\mathcal{G})[\frac{1}{p}] \cap \Lambda(\mathcal{G})_S = \Lambda(\mathcal{G})$. Similarly, $x^{-1} \in \Lambda(\mathcal{G})$. Hence $x \in \Lambda(\mathcal{G})^{\times}$. Hence cokernel of $\theta$ injects in the cokernel of $\theta_S$ and $\theta_0$ is injective.
\qed

Next we show how our formulation of the main conjecture holds for admissible $p$-adic Lie extension whose Galois group is of the form $\Delta \times \mathcal{G}$, with $\Delta$ a finite group whose order is prime to $p$ and $\mathcal{G}$ is abelian. Let $\psi$ be an arbitrary Artin representation of $Gal(\overline{F}/F)$ such that $F_{\psi} \cap F_{\infty} =F$. Put $\Delta = Gal(F_{\psi}/F)$. If $\psi$ is one dimensional, then $Gal(F_{\psi}F_{\infty}/F) \cong \Delta \times \mathcal{G}$ is abelian. We have the Deligne-Ribet $p$-adic zeta function 
\[
\zeta \in Q(\Lambda(\Delta \times \mathcal{G})).
\]
We define the $p$-adic $L$-function $L_p(\psi) \in Q(\Lambda_{\psi}(\mathcal{G}))$ as the image of $\zeta$ under the map
\[
Q(\Lambda(\Delta \times \mathcal{G})) \rightarrow Q(\Lambda_{O_{\psi}}(\mathcal{G})),
\]
induced by $\psi$.

Now assume that $\psi$ is an irreducible Artin character of degree greater than 1. Then $\Delta$ is nonabelian. By Brauer induction theorem there are subgroups $\Delta_i$ of $\Delta$ and degree one character $\psi_i$ of $\Delta_i$ such that
\[
\psi = \sum a_i Ind^{\Delta}_{\Delta_i} \psi_i, \hspace{1cm} a_i \in \mathbb{Z}.
\]
Let $L_i = F_{\psi}^{\Delta_i}$. We view $\mathcal{G}$ as a Galois group of $L_iF_{\infty}/L_i$ and get the $p$-adic $L$-function $L_p(\psi_i) \in \Lambda_{O_{\psi_i}}(\mathcal{G})$. For the character $\psi$ we define the $p$-adic $L$ function $L_p(\psi)$ by
\[
L_p(\psi) = \prod_{i} L_p(\psi_i)^{a_i} \in Q(\Lambda_{O_{\psi}}(\mathcal{G})).
\]
It follows from the classical main conjecture and the assumption $\mu =0$ that if $p \nmid |\Delta|$, then $L_p(\psi) \in \Lambda_{O_{\psi}}(\mathcal{G})_S^{\times}$ (This is the $p$-adic Artin conjecture formulated by Greenberg \cite{Greenberg:1983}). If $p \nmid |\Delta|$, then we define $\zeta$ as 
\[
\zeta = \sum_{\chi \in R(\Delta)} e_{\chi} L_p(\chi),
\]
where $R(\Delta)$ is the set of irreducible characters of $\Delta$. 

\begin{lemma} $\zeta$ belongs to $K_1(\Lambda(\mathcal{G} \times \Delta)_S)$ and is the $p$-adic zeta function for the extension $F_{\psi}F_{\infty}/F$. Note that $K_1'(\Lambda(\mathcal{G} \times \Delta)_S) = K_1(\Lambda(\mathcal{G} \times \Delta)_S)$.
\label{primetopzeta}
\end{lemma}

\noindent{\bf Proof:} Since the order of $\Delta$ is prime to $p$, we have the isomorphism
\[
K_1(\Lambda(\mathcal{G} \times \Delta)_S) \xrightarrow{\sim} \oplus_{\chi \in R(\Delta)} \Lambda_{O_{\chi}}(\mathcal{G})_S^{\times}.
\]
Since the image of $\zeta$ under this isomorphism is $(L_p(\chi))$, we have the lemma.
\qed

\begin{theorem} Let $F_{\infty}/F$ be an admissible $p$-adic Lie extension satisfying hypothesis $\mu=0$. Assume that $Gal(F_{\infty}/F) = \mathcal{G} \times \Delta$, with $\mathcal{G}$ abelian and $\Delta$ a finite group of order prime to $p$. Then the main conjecture is true for the extension $F_{\infty}/F$. Moreover, we have an injection 
\[
K_1(\Lambda(\mathcal{G} \times \Delta)) \rightarrow K_1(\Lambda(\mathcal{G} \times \Delta)_S).
\]
\label{primetopmainconjecture}
\end{theorem}
\noindent{\bf Proof:} Let $C^{\cdot}_{\chi} = (\Lambda_{O_{\chi}}(\mathcal{G})_S \otimes^{L}_{\Lambda(\mathcal{G} \times \Delta)_S} C(F_{\infty}/F))$ for all $\chi \in R(\Delta)$. Then 
\[
\partial(L_p(\chi)) = -[C^{\cdot}_{\chi}] \in K_0(\Lambda_{O_{\chi}}(\mathcal{G}), \Lambda_{O_{\chi}}(\mathcal{G})_S).
\]
The theorem follows from the following isomorphism
\[
K_0(\Lambda(\mathcal{G} \times \Delta), \Lambda(\mathcal{G} \times \Delta)_S) \rightarrow \oplus_{\chi \in R(\Delta)} K_0(\Lambda_{O_{\chi}}(\mathcal{G}), \Lambda_{O_{\chi}}(\mathcal{G})_S).
\]
The injection of $K_1(\Lambda(\mathcal{G} \times \Delta))$ in $K_1(\Lambda(\mathcal{G} \times \Delta)_S)$ follows from the fact that 
\[
K_1(\Lambda(\mathcal{G} \times \Delta)) \cong \oplus_{\chi \in R(\Delta)} \Lambda_{O_{\chi}}(\mathcal{G})^{\times}
\]
 and 
 \[
 K_1(\Lambda(\mathcal{G} \times \Delta)_S) \cong \oplus_{\chi \in R(\Delta)} \Lambda_{O_{\chi}}(\mathcal{G})_S^{\times}.
 \]

\qed

\section{Various reductions} \label{sectionreductions} In this section we prove several algebraic results and obtain various reductions. In this section $O$ is the ring of integers in a finite extension of $\mathbb{Q}_p$.

\subsection{Reduction to one dimensional case} In this section we show that the validity of the main conjecture for one dimensional $p$-adic Lie group implies the main conjecture for $p$-adic Lie extensions of arbitrary dimension. This is proven by D. Burns \cite{Burns:2010}. 

Let $Q_1(\mathcal{G}) = \{ \mathcal{G}/U : U \text{ is an open pro-$p$ subgroup of } H \text{ and is normal in } \mathcal{G}\}$. Using a result of Fukaya-Kato \cite{FukayaKato:2006} we have

\begin{lemma} Let $P$ be a compact $p$-adic Lie group. Then we have an isomorphism
\[
K_1'(\Lambda_O(P)) \rightarrow \ilim{\Delta} K_1'(O[\Delta]),
\]
where $\Delta$ runs through all finite quotients of $P$.
\label{inverselimitofk1}
\end{lemma}
\noindent{\bf Proof:} For a compact $p$-adic Lie group $P$, let $J_P$ be the Jacobson radical of $\Lambda(P)$. Then Fukaya-Kato (proposition 1.5.1 \cite{FukayaKato:2006}) show that 
\[
K_1(\Lambda(P)) \xrightarrow{\sim} \ilim{n} K_1(\Lambda(P)/J_{P}^n)
\]
In the rest of the proof $\Delta$ varies over all finite quotients of $P$. For a non-negative integer $n$ define
\[
I_{\Delta,n} := ker(\Lambda_O(P) \rightarrow O[\Delta]/J_{\Delta}^n).
\]
Clearly, $J_P^n \subset \cap_{\Delta} I_{\Delta,n}$. Consider the exact sequence 
\[
K_1(\Lambda_O(P)/J_P^n,I_{\Delta,n}/J_P^n) \rightarrow K_1(\Lambda_O(P)/J_P^n) \rightarrow K_1(O[\Delta]/J_{\Delta}^n) \rightarrow 1.
\]
As all the groups in the above short exact sequence are finite, we have an exact sequence
\[
\ilim{\Delta}K_1(\frac{\Lambda_O(P)}{J_P^n},\frac{I_{\Delta,n}}{J_P^n}) \rightarrow K_1(\frac{\Lambda_O(P)}{J_P^n}) \rightarrow \ilim{\Delta} K_1(\frac{O[\Delta]}{J_{\Delta}^n}) \rightarrow 1.
\]
We have a surjection 
\[
1+ I_{\Delta,n}/J_P^n \rightarrow K_1(\Lambda_O(P)/J_P^n,I_{\Delta,n}/J_P^n) \rightarrow 1.
\]
Since $I_{\Delta,n}/J_P^n$ is finite we get the induced surjection
\[
\ilim{\Delta}\Big(1+I_{\Delta,n}/J_P^n \Big) = 1+ (\cap_{\Delta}I_{\Delta,n})/J_P^n \rightarrow \ilim{\Delta} K_1(\Lambda_O(P)/J_P^n,I_{\Delta,n}/J_P^n) \rightarrow 1.
\]
In particular, we deduce that $\ilim{\Delta} K_1(\Lambda_O(P)/J_P^n,I_{\Delta,n}/J_P^n)$ is finite and 
\[
\ilim{n}\Big(1+ (\cap_{\Delta}I_{\Delta,n})/J_P^n\Big) \rightarrow \ilim{n} \ilim{\Delta} K_1(\Lambda_O(P)/J_P^n,I_{\Delta,n}/J_P^n) \rightarrow 1
\]
is exact. But $\ilim{n}\Big(1+ (\cap_{\Delta}I_{\Delta,n})/J_P^n\Big) = \{1\}$ as $\cap_{n}\cap_{\Delta} I_{\Delta,n} = \{0\}$. Hence $\ilim{n} \ilim{\Delta} K_1(\Lambda_O(P)/J_P^n,I_{\Delta,n}/J_P^n)=\{1\}$ and from the above exact sequence we get
\[
\ilim{n} K_1(\Lambda_O(P)/J_P^n) \xrightarrow{\cong} \ilim{n}\ilim{\Delta} K_1(O[\Delta]/J_{\Delta}^n).
\]
Hence
\begin{align*}
K_1(\Lambda_O(P)) & \cong \ilim{n} \ilim{\Delta} K_1(O[\Delta]/J_{\Delta}^n) \\
& \cong \ilim{\Delta} \ilim{n} K_1(O[\Delta]/J_{\Delta}^n) \\
& \cong \ilim{\Delta} K_1(O[\Delta]) 
\end{align*}
The interchange of inverse limits is justified because the groups $K_1(O[\Delta]/J_{\Delta}^n)$ are finite for all $\Delta$ and all $n$. The first and third isomorphism use the result of Fukaya-Kato mentioned above. Hence by definition of $SK_1(\Lambda_O(P))$, we have
\[
K_1'(\Lambda_O(P)) \cong \ilim{\Delta}K_1(O[\Delta])/\ilim{\Delta}SK_1(O[\Delta]).
\]
But since $SK_1(O[\Delta])$ are finite for all $\Delta$ (by a result of Wall,  theorem 2.5 in \cite{Oliver:1988}), we have
\[
\ilim{\Delta}K_1(O[\Delta])/\ilim{\Delta}SK_1(O[\Delta]) \cong \ilim{\Delta} K_1'(O[\Delta]).
\]
\qed

\begin{corollary} The following natural map is an isomorphism
\[
K_1'(\Lambda(\mathcal{G})) \xrightarrow{\sim} \ilim{G \in Q_1(\mathcal{G})} K_1'(\Lambda(G)).
\]
\end{corollary}

\noindent{\bf Proof:} 
\begin{align*} 
\ilim{G \in Q_1(\mathcal{G})} K_1'(\Lambda(G)) & \cong \ilim{G} \ilim{\Delta_G} K_1'(\mathbb{Z}_p[\Delta_G]) \\
& \cong \ilim{\Delta} K_1'(\mathbb{Z}_p[\Delta]) \\
& \cong K_1'(\Lambda(\mathcal{G})),
\end{align*}
where $\Delta_G$ runs through finite quotients of $G$ and $\Delta$ runs through all finite quotients of $\mathcal{G}$.
\qed

\begin{theorem} Assume $F_{\infty}/F$ is an admissible $p$-adic Lie extension satisfying $\mu=0$ hypothesis. Then the main conjecture is true for $F_{\infty}/F$ if it is true for each of the extensions $F_{\infty}^U/F$ for all open pro-$p$ subgroups $U \leq H$ such that $U$ is normal in $\mathcal{G}$ and if for each $G = \mathcal{G}/U$ the group $K_1'(\Lambda(G))$ injects into $K_1'(\Lambda(G)_S)$.
\label{reductiontodim1}
\end{theorem}

\noindent{\bf Proof:} Note that for each $U$ the extension $F_{\infty}^U/F$ is an admissible $p$-adic Lie extension satisfying $\mu=0$ hypothesis. We consider the following diagram
\[
\xymatrix{ & K_1'(\Lambda(\mathcal{G})) \ar[r] \ar[d]_{\sim} & K_1'(\Lambda(\mathcal{G})_S) \ar[r]^{\partial} \ar[d] & K_0(\Lambda(\mathcal{G}), \Lambda(\mathcal{G})_S) \ar[r] \ar[d] & 0 \\
1 \ar[r] & \ilim{G}K_1'(\Lambda(G)) \ar[r] & \ilim{G} K_1'(\Lambda(G)_S) \ar[r] & \ilim{G} K_0(\Lambda(G), \Lambda(G)_S) & }
\]
We assume that the main conjecture is true for each quotient $G = \mathcal{G}/U$ in $Q_1(\mathcal{G})$. Then for each $G \in Q_1(\mathcal{G})$ we have a unique $\zeta_G \in K_1'(\Lambda(G)_S)$ satisfying the main conjecture. Uniqueness of $\zeta_G$ implies that $(\zeta_G)_G \in \ilim{G} K_1'(\Lambda(G)_S)$. Let $f^{-1} \in K_1'(\Lambda(\mathcal{G})_S)$ be a characteristic element of $C(F_{\infty}/F)$. Let $(f_G)_{G} \in \ilim{G} K_1'(\Lambda(G)_S)$ be the image of $f$. Then $\partial(f_G) = -[C(F_{\infty}^U/F)]$ by isomorphism (\ref{decent}). Put $u_G = \zeta_G f_G^{-1} \in K_1'(\Lambda(G)_S)$. Then $(u_G)_{G \in Q_1(\mathcal{G})} \in \ilim{G}K_1'(\Lambda(G))$ as $\partial(\zeta_G) = -[C(F_{\infty}^U/F)]$. Hence there is $u \in K_1'(\Lambda(\mathcal{G}))$ which maps to $(u_G)_{G}$. We claim that $\zeta = uf$ is the $p$-adic zeta function we seek. Uniqueness is clear and so is the fact that $\partial(\zeta) = -[C(F_{\infty}/F)]$. Now we show that $\zeta$ satisfies the required interpolation property. Let $\rho$ be an Artin representation of $\mathcal{G}$, then there is a $G \in Q_1(\mathcal{G})$, such that $\rho$ factors through $G$. Since $\zeta$ maps to $\zeta_G \in K_1(\Lambda(G)_S)$ we have for any positive integer $r$ divisible by $[F_{\infty}(\mu_p):F_{\infty}]$
\[
\zeta(\rho\kappa_F^r) = \zeta_G(\rho\kappa_F^r) = L_{\Sigma}(\rho, 1-r).
\]
Note that since $U$ is pro-$p$, we have $[F_{\infty}(\mu_p):F_{\infty}] = [F_{\infty}^U(\mu_p):F_{\infty}^U]$.
\qed

\begin{remark} Since we interpolate Artin $L$-values the above result is not surprising though the proof is delicate. It is a typical example of the strategy of Burns and Kato (explained in the introduction) and will appear in all our reduction steps. 
\end{remark}

\subsection{Reduction to $\mathbb{Q}_p$-elementary extensions} Now we consider one dimensional $\mathcal{G}$. We pick and fix a lift of $\Gamma$ in $\mathcal{G}$ and hence get an isomorphism $\mathcal{G} \cong H \rtimes \Gamma$. Let us fix $\Gamma^{p^e}$, an open subgroup of $\Gamma$ acting trivially on $H$. Let $G = \mathcal{G}/\Gamma^{p^e}$. 

\begin{definition} Let $l$ be a prime. A finite group $P$ is called $l$-hyperelementary if $P$ is of the form $C_n \rtimes P_1$, with $C_n$ a cyclic group of order $n$ and $P_1$ a finite $l$-group and $l \nmid n$. Let $K$ be a field of characteristic 0. A $l$-hyperelementary group $C_n \rtimes P_1$ is called $l$-$K$-elementary if 
\[
Im[P_1 \rightarrow Aut(C_n) \cong (\mathbb{Z}/n\mathbb{Z})^{\times}] \subset Gal(K(\mu_n)/K).
\]
A hyperelementary group is one which is $l$-hyperelementary for some prime $l$. A $K$-elementary group is one which is $l$-$K$-elementary for some prime $l$.
\label{defnhyperelementary}
\end{definition} 

\begin{definition} Let $l$ be a prime. A $p$-adic Lie group is called $l$-$K$-elementary if it is of the form $P \rtimes \Gamma$ for a finite group $P$ and there is a central open subgroup $\Gamma^{p^r}$ of $P \rtimes \Gamma$ such that $(P \rtimes \Gamma)/\Gamma^{p^r}$ is $l$-$K$-elementary finite group. A $p$-adic Lie group is $K$-elementary if it is $l$-$K$-elementary for some prime $l$. 
\end{definition}

The induction theory of A. Dress \cite{Dress:1973} gives the following theorem (see theorem 11.2 in R. Oliver \cite{Oliver:1988} and also theorem 5.4 CTC Wall \cite{Wall:1974}). 

\begin{theorem} (Dress, Wall) Let $P$ be a finite group. Let $K$ be the field of fractions of $O$. Then we have the following isomorphism
\[
K_1'(O[P]) \xrightarrow{\sim} \ilim{\pi} K_1'(O[\pi]),
\]
where $\pi$ runs through all $K$-elementary subgroups of $P$. The inverse limit is with respect to norm maps and the maps induced by conjugation. In other words, $(x_{\pi}) \in \prod_{\pi} K_1'(O[\pi])$ lies in $\ilim{\pi} K_1'(O[\pi])$ if and only if \\
(i) for all $g \in P$,  $gx_{\pi}g^{-1} = x_{g\pi g^{-1}}$, and \\
(ii) for $\pi \leq \pi' \leq P$ the norm homomorphism $K_1'(O[\pi']) \rightarrow K_1'(O[\pi])$ maps $x_{\pi'}$ to $x_{\pi}$.

\label{inductiontheorem}
\end{theorem}

\begin{lemma} For a fixed $n \geq 0$, let $G_n = \mathcal{G}/\Gamma^{p^{e+n}}$. We have an isomorphism 
\[
K_1'(\Lambda(\mathcal{G})) \xrightarrow{\sim} \ilim{P} K_1'(\Lambda(U_P)),
\]
where the inverse limit ranges over all $\mathbb{Q}_p$-elementary subgroups $P$ of $G_n$ and it is with respect to the maps induced by conjugation and the norm maps. $U_P$ denotes the inverse image of $P$ in $\mathcal{G}$ under the surjection $\mathcal{G} \rightarrow G_n$.
\label{hyperelementaryinduction}
\end{lemma}

\noindent{\bf Proof:}  In this proof $P$ runs though all $\mathbb{Q}_p$-elementary subgroups of $G_n$ and all inverse limits are with respect to the maps induced by conjugation and the norm maps. We have an isomorphism 
\[
K_1'(\mathbb{Z}_p[G_n]) \xrightarrow{\sim} \ilim{P} K_1'(\mathbb{Z}_p[P])
\]
by theorem \ref{inductiontheorem}. We claim that
\[
K_1'(\mathbb{Z}_p[G_{n+i}]) \xrightarrow{\sim} \ilim{P} K_1'(\mathbb{Z}_p[U_P/\Gamma^{p^{e+n+i}}]).
\]
Note that $U_P/\Gamma^{p^{e+n+i}}$ is the inverse image of $P$ in the group $G_{n+i}$ under the surjection $G_{n+i} \rightarrow G_n$. Any $\mathbb{Q}_p$-elementary subgroup of $G_{n+i}$ is contained in $U_P/\Gamma^{p^{e+n+i}}$ for some $P$. Hence 
\begin{align*}
\ilim{P}K_1'(\mathbb{Z}_p[U_P/\Gamma^{p^{e+n+i}}]) 
& \xrightarrow{\sim} \ilim{P} \ilim{C_P} K_1'(\mathbb{Z}_p[C_P]) \\
& \xrightarrow{\sim} \ilim{Q} K_1'(\mathbb{Z}_p[Q])  \\
& \xrightarrow{\sim} K_1'(\mathbb{Z}_p[G_{n+i}]).
\end{align*}
Here $C_P$ runs through all $\mathbb{Q}_p$-elementary subgroups of $U_P/\Gamma^{p^{e+n+i}}$ and $Q$ runs through all $\mathbb{Q}_p$-elementary subgroups of $G_{n+i}$. Hence we get the claim. Passing to the inverse limit we get the result by lemma \ref{inverselimitofk1}. 
\qed

\begin{theorem} Assume that the main conjecture is valid for all admissible $p$-adic Lie extensions satisfying $\mu=0$ hypothesis whose Galois group is $\mathbb{Q}_p$-elementary. Also assume that for all $\mathbb{Q}_p$-elementary $p$-adic Lie groups $U$ the group $K_1'(\Lambda(U))$ injects in $K_1'(\Lambda(U)_S)$. Then the main conjecture is valid for all one dimensional admissible $p$-adic Lie extensions $F_{\infty}/F$ satisfying $\mu=0$ hypothesis. Moreover, if $\mathcal{G} = Gal(F_{\infty}/F)$, then $K_1'(\Lambda(\mathcal{G}))$ injects in $K_1'(\Lambda(\mathcal{G})_S)$.
\label{reductiontohyperelementary} 
\end{theorem}

\noindent{\bf Proof:} Let $F_{\infty}/F$ be an admissible $p$-adic Lie extension of dimension one satisfying $\mu=0$ hypothesis. Let $\mathcal{G} =Gal(F_{\infty}/F)$. Then all its admissible $p$-adic Lie subextensions satisfy $\mu=0$ hypothesis. Assume that the main conjecture is valid for all admissible $p$-adic Lie subextensions of $F_{\infty}/F$ whose Galois group is $\mathbb{Q}_p$-elementary. As before, for a fixed but arbitrary $n \geq 0$, let $G_n = \mathcal{G}/\Gamma^{p^{e+n}}$. Consider the following diagram 
\[
\xymatrix{ & K_1'(\Lambda(\mathcal{G})) \ar[r] \ar[d]_{\sim} & K_1'(\Lambda(\mathcal{G})_S) \ar[r] \ar[d] & K_0(\Lambda(\mathcal{G}), \Lambda(\mathcal{G})_S) \ar[r] \ar[d] & 0 \\
1 \ar[r] & \ilim{P} K_1'(\Lambda(U_P)) \ar[r] & \ilim{P} K_1'(\Lambda(U_P)_S) \ar[r] & \ilim{P} K_0(\Lambda(U_P), \Lambda(U_P)_S) &}
\]
Here $P$ runs through all $\mathbb{Q}_p$-elementary subgroups of $G_n$ and $U_P$ denotes the inverse image of $P$ in $\mathcal{G}$.

Take $f \in K_1'(\Lambda(\mathcal{G})_S)$ such that $\partial(f)= -[C(F_{\infty}/F)]$. Let $(f_P)$ be the image of $f$ in $\ilim{P} K_1'(\Lambda(U_P)_S)$. Then $\partial(f_P)=-[C(F_{\infty}/F_{\infty}^{U_P})]$. For each $P$ let $\zeta_P$ be the $p$-adic zeta function satisfying the main conjecture for $F_{\infty}/F_{\infty}^{U_P}$. We claim that uniqueness implies that $(\zeta_P) \in \ilim{P} K_1'(\Lambda(U_P)_S)$. We must show that for any $g \in \mathcal{G}$, we have $g\zeta_Pg^{-1} = \zeta_{gPg^{-1}}$ and if $P \leq P' \leq G_n$, then $norm(\zeta_{P'}) = \zeta_P$. Note that 
\[
\partial(g\zeta_Pg^{-1}) = \partial(\zeta_{gPg^{-1}}) = -[C(F_{\infty}/F_{\infty}^{U_{gPg^{-1}}}].
\]
If $\rho$ is an Artin Character of $U_{gPg^{-1}}$, then we denote by $\rho_g$ the character of $U_P$ given by $\rho_g(h) = \rho(ghg^{-1})$ for any $h \in U_P$. Then 
\begin{align*}
g\zeta_Pg^{-1} (\rho\kappa^r) & = \zeta_P(\rho_g \kappa^r) \\
& = L_{\Sigma}(\rho_g, 1-r) \\
& = L_{\Sigma}(\rho, 1-r) \\
& = \zeta_{gPg^{-1}}(\rho\kappa^r)
\end{align*}
Hence $g\zeta g^{-1} = \zeta_{gPg^{-1}}$.

We have $\partial(norm(\zeta_{P'})) = -[C(F_{\infty}/F_{\infty}^{U_P}] = \partial(\zeta_P)$. On the other hand, if $\rho$ is an Artin Character of $U_P$, then 
\begin{align*}
norm(\zeta_{P'})(\rho\kappa^r) & = \zeta_{P'}(Ind^{U_{P'}}_{U_P}(\rho)) \\
& = L_{\Sigma}(Ind^{U_{P'}}_{U_P}(\rho), 1-r) \\
& = L_{\Sigma}(\rho, 1-r) \\
& = \zeta_P(\rho \kappa^r)
\end{align*}
Hence $norm(\zeta_{P'}) = \zeta_P$ and our claim holds.

Put $u_P = \zeta_Pf_P^{-1} \in K_1'(\Lambda(U_P)_S)$. As $\partial(\zeta_P) = \partial(f_P) = -[C(F_{\infty}/F_{\infty}^{U_P})]$, we have $u_P \in K_1'(\Lambda(U_P))$. Moreover, $(u_P)_P \in \ilim{P}K_1'(\Lambda(U_P))$. Then there is a $u \in K_1'(\Lambda(\mathcal{G}))$ mapping to $(u_P)_P$.  Let $\zeta = uf$. It is clear that $\partial (\zeta) = -[C(F_{\infty}/F)]$. This is the $p$-adic zeta function we seek. First note that it is the only element in $K_1(\mathcal{G})_S)$ such that $\partial(\zeta) = -[C(F_{\infty}/F)]$ and whose image in $K_1'(\Lambda(U_P)_S)$ is $\zeta_P$. Hence it is independent of the choice of $n$. 

 Let $\rho$ be an Artin character of $\mathcal{G}$. It factors through $G_n$ for some $n$. Then by Brauer's induction theorem (theorem 19, \cite{Serre:representationtheory}) there are Artin representations $\rho_P$ of $U_P$ such that
\[
\rho = \sum_{P} n_P Ind^{\mathcal{G}}_{U_P}\rho_P.
\]
Note that an elementary subgroup (as defined in \cite{Serre:representationtheory}) is $\mathbb{Q}_p$-elementary. Though for Brauer's induction theorem we only need elementary subgroups we need to work with $\mathbb{Q}_p$-elementary subgroups for the isomorphism in the previous lemma. Also note than $\rho_P$ may not be one dimensional. Then for any positive integer $r$ divisible by $[F_{\infty}(\mu_p):F_{\infty}]$, we have
\begin{align*}
\zeta(\rho\kappa_F^r) & = \prod_{P} \zeta(Ind^{\mathcal{G}}_{U_P}\rho_P\kappa_F^r)^{n_P}  \\
& = \prod_{P} \zeta_{P}(\rho_P\kappa_{F_P}^r)^{n_P} \\
&= \prod_{P} L_{\Sigma}(\rho_P, 1-r)^{n_P} \\
& = L_{\Sigma}(\rho, 1-r).
\end{align*}
Hence $\zeta$ satisfies the required interpolation property. Uniqueness of $\zeta$ and assertion about the $K_1'$ groups follows by an easy diagram chase. 
\qed

Hence we have reduced the main conjecture to extensions whose Galois group is a $\mathbb{Q}_p$-elementary $p$-adic Lie group.

\subsection{The case of $l$-$\mathbb{Q}_p$-elementary for $l \neq p$} We now assume that $\mathcal{G}$ is $l$-$\mathbb{Q}_p$-elementary for some prime $l \neq p$. 

\begin{lemma} If a $p$-adic Lie group is $l$-$\mathbb{Q}_p$-elementary for some prime $l \neq p$, then it is isomorphic to $\Gamma^{p^e} \times G$ for some finite $l$-$\mathbb{Q}_p$-elementary group $G$ and some integer $e \geq 0$.
\label{hyperelementarystructure}
\end{lemma}

\noindent{\bf Proof:} Let $\mathcal{G} = H \rtimes \Gamma$. Let $\Gamma^{p^r}$ be a central open subgroup of $\mathcal{G}$ such that $\mathcal{G}/\Gamma^{p^r}$ is $l$-$\mathbb{Q}_p$-elementary. Then for any $s \leq r$ such that $\Gamma^{p^s}$ is a central subgroup of $\mathcal{G}$, the quotient $\mathcal{G}/\Gamma^{p^s}$ is $l$-$\mathbb{Q}_p$-elementary. Let $e \geq 0$ be the smallest integer such that $\Gamma^{p^e}$ is an open central subgroup of $\mathcal{G}$. Let $G = \mathcal{G}/\Gamma^{p^e}$. It is a $l$-$\mathbb{Q}_p$-elementary finite group. We claim that $\mathcal{G} = \Gamma^{p^e} \times G$. Let $G = C \rtimes P$, where $C$ is a cyclic subgroup of order prime to $l$ and $P$ is a $l$-group. Note that $C \rtimes P = G = H \rtimes \Gamma/\Gamma^{p^e}$. Since order of $P$ is prime to $p$ it is a subgroup of $H$. Hence $P$ is also a subgroup of $\mathcal{G}$. Let $U_C$ be the inverse image of $C$ in $\mathcal{G}$. Since $\Gamma^{p^e}$ is central and $C$ is cyclic, $U_C$ is an abelian group. Write $U_C = Q \times D$, where $Q$ is a pro-$p$ pro-cyclic subgroup of $U_C$ and $D$ is a torsion subgroup. Then $\mathcal{G} = (Q \times D) \rtimes P$. The action of $P$ on $Q$ has to be trivial since it is trivial on an open subgroup of $Q$. Hence $Q$ is a central pro-$p$ pro-cyclic subgroup of $\mathcal{G}$ which implies by minimality of $e$ that $Q = \Gamma^{p^e}$ and $D= C$. 
\qed

Hence we may assume that the $l$-$\mathbb{Q}_p$-elementary group $\mathcal{G}$ is of the form $\Gamma \times H$, with $H$ a $l$-$\mathbb{Q}_p$-elementary finite group. Let $H = C \rtimes P$, with $C$ a cyclic group of order prime to $l$ and $P$ a finite $l$-group. If $p$ does not divide the order of $C$, then the main conjecture for $F_{\infty}/F$ is valid by theorem \ref{primetopmainconjecture}. 

Write $C$ as $C_p \times C_{p'}$, where $C_p$ is a cyclic group of $p$-power order and $C_{p'}$ is a cyclic group of order prime to $p$. Note that both $C_p$ and $C_{p'}$ are normal subgroups of $H$. Hence, after replacing $C$ by $C_p$ and $P$ by $C_{p'} \rtimes P$, we assume that $H = C \rtimes P$, with $P$ a finite group of order prime to $p$ and $C$ a cyclic group of $p$-power order. By proposition 12.7 in Oliver \cite{Oliver:1988}, $SK_1(\Lambda(\mathcal{G})) =1$. Put $\overline{C} = H_0(P,C)$, $\overline{H} = \overline{C} \times P$ and $\overline{\mathcal{G}} = \Gamma \times \overline{H}$. Put $P_1 = Ker(P \rightarrow Aut(C))$, $H_1= C \times P_1$ and $\mathcal{G}_1 = \Gamma \times H_1$. Consider the maps
\[
\theta: K_1(\Lambda(\mathcal{G})) \rightarrow K_1(\Lambda(\overline{\mathcal{G}})) \times K_1(\Lambda(\mathcal{G}_1)),
\]
\[
\theta_S : K_1(\Lambda(\mathcal{G})_S) \rightarrow K_1(\Lambda(\overline{\mathcal{G}})_S) \times K_1(\Lambda(\mathcal{G}_1)_S).
\]
In both cases the map in the first component is the one induced by natural surjection and the map in the second component is the norm map.

\begin{definition} Let $\Phi$ (resp. $\Phi_S$) be the set of all pairs $(x_0,x_1) \in K_1(\Lambda(\overline{\mathcal{G}})) \times K_1(\Lambda(\mathcal{G}_1))$ (resp. $K_1(\Lambda(\overline{\mathcal{G}})_S) \times K_1(\Lambda(\mathcal{G}_1)_S)$) such that \\
1. Under the map
\[
K_1(\Lambda(\overline{\mathcal{G}})) \xrightarrow{norm} K_1(\Lambda(\Gamma \times \overline{C} \times P_1)) \xleftarrow{can} K_1(\Lambda(\mathcal{G}_1))
\]

\[
(\text{resp. }K_1(\Lambda(\overline{\mathcal{G}})_S) \xrightarrow{norm} K_1(\Lambda(\Gamma \times \overline{C} \times P_1)_S) \xleftarrow{can} K_1(\Lambda(\mathcal{G}_1)_S)),
\]
we have $norm(x_0) = can(x_1)$. \\
2. $x_1$ is fixed under the conjugation action by every element of $P$.
\label{definitionofphi}
\end{definition}

\begin{proposition} The map $\theta$ induces an isomorphism between $K_1(\Lambda(\mathcal{G}))$ and $\Phi$. The image of $\theta_S$ is contained in $\Phi_S$. In particular 
\[
Im(\theta_S) \cap (K_1(\Lambda(\overline{\mathcal{G}})) \times K_1(\Lambda(\mathcal{G}_1))) = Im(\theta).
\]
Note that $K_1(\Lambda(\overline{\mathcal{G}}))$ and $K_1(\Lambda(\mathcal{G}_1))$ inject in $K_1(\Lambda(\overline{\mathcal{G}})_S)$ and $K_1(\Lambda(\mathcal{G}_1)_S)$ respectively and hence the above intersection makes sense in $K_1(\Lambda(\overline{\mathcal{G}})_S) \times K_1(\Lambda(\mathcal{G}_1)_S)$.   
\label{propprimetopk1}
\end{proposition}

\noindent{\bf Proof:} We first verify that the image of $\theta$ and $\theta_S$ is contained in $\Phi$ and $\Phi_S$ respectively. We do this only for $\theta$ and the proof of $\theta_S$ is similar. Let $B$ be a set of right coset representatives of $H_1$ in $H$. Then $B$ is a basis of the $\Lambda(\mathcal{G}_1)$-module $\Lambda(\mathcal{G})$. Let $x \in K_1(\Lambda(\mathcal{G}))$. The image of $x$ in $K_1(\Lambda(\mathcal{G}_1))$ under the norm map is explicitly described as follows. Let $\tilde{x}$ be a lift of $x$ in $\Lambda(\mathcal{G})^{\times}$. Multiplication on the right by $\tilde{x}$ gives a $\Lambda(\mathcal{G}_1)$-linear map on $\Lambda(\mathcal{G})$. Let $A(B,\tilde{x})$ be the matrix of this map with respect to the basis $B$. Then norm of $x$ is the class of this matrix in $K_1(\Lambda(\mathcal{G}_1))$. It is independent of the choice of $\tilde{x}$ and the basis $B$.

Let $g \in P$. Then $gBg^{-1}$ is also a $\Lambda(\mathcal{G}_1)$-basis of $\Lambda(\mathcal{G})$. Then $gA(B,\tilde{x})g^{-1} = A(gBg^{-1},g\tilde{x}g^{-1})$. Since $g\tilde{x}g^{-1}$ is also a lift of $x$, the class of $A(gBg^{-1},g\tilde{x}g^{-1})$ in $K_1(\Lambda(\mathcal{G}_1))$ is the same as that of $A(B,\tilde{x})$. Hence $\theta(x)$ satisfies the second condition in the definition of $\Phi$.

Since we can choose the image of $B$ in $\Lambda(\overline{\mathcal{G}})$ as a basis of the $\Lambda(\Gamma \times \overline{C} \times P_1)$-module $\Lambda(\overline{\mathcal{G}})$, the following diagram commutes.
\[
\xymatrix{ K_1(\Lambda(\mathcal{G})) \ar[rr]^{norm} \ar[d]_{can} & & K_1(\Lambda(\mathcal{G}_1)) \ar[d]^{can} \\
K_1(\Lambda(\overline{\mathcal{G}})) \ar[rr]_{norm} & & K_1(\Lambda(\Gamma \times \overline{C} \times P_1))}
\]
Hence $\theta(x)$ satisfies the first condition in the definition of $\Phi$. 

Next we show that $\theta$ surjects on $\Phi$. Let $(x_0,x_1) \in \Phi$. The natural map from $K_1(\Lambda(\mathcal{G}))$ to $K_1(\Lambda(\overline{\mathcal{G}}))$ is a surjection (since $\Lambda(\mathcal{G})^{\times}$ surjects on $K_1(\Lambda(\mathcal{G}))$ and also on $\Lambda(\overline{\mathcal{G}})^{\times}$. Moreover, $\Lambda(\overline{\mathcal{G}})^{\times}$ surjects on $K_1(\Lambda(\overline{\mathcal{G}}))$), hence we may and do assume that $x_0 = 1$. Then the first condition in the definition of $\Phi$ implies that 
\[
x_1 \in J:= Ker(can: K_1(\Lambda(\mathcal{G}_1)) \rightarrow K_1(\Lambda(\Gamma \times \overline{C} \times P_1))). 
\]
By theorem 2.10 (ii) Oliver \cite{Oliver:1988} and proposition 1.5.3 Fukaya-Kato \cite{FukayaKato:2006} the subgroup $J$ is pro-$p$; hence so is the subgroup $J^P$ which is the subgroup of point-wise fixed elements of $J$ under the conjugation action of $P$. Also $x_1 \in J^P$. Let $n$ be the order of $P/P_1$. Let $z \in J^P$ be such that $z^n = x_1$. Denote by $\tilde{z}$ the image of $z$ in $K_1(\Lambda(\mathcal{G}))$ under the natural map $K_1(\Lambda(\mathcal{G}_1)) \rightarrow K_1(\Lambda(\mathcal{G}))$. Let $\theta(\tilde{z}) = (z_0, z_1)$. By construction $z_1 = x_1$. Moreover, $norm(z_0) = z_0^n = 1$ in $K_1(\Lambda(\Gamma \times \overline{C} \times P_1))$. But $z_0$ lies  in a pro-$p$ subgroup and hence $z_0=1$. This proves the surjection of $\theta$ on $\Phi$. 

Lastly, we show that $\theta$ is injective. We need a lemma

\begin{lemma} Let $\widehat{C}$ be the set of irreducible characters of $C$. Then $P$ acts on them by
\[
(g\cdot\chi)(h) = \chi(ghg^{-1}), \hspace{1cm} \text{for all $g \in P$ and $h \in C$}
\]
Under this action satbiliser of $\chi$ is $P_1$ if $\chi \neq 1$. Hence any irreducible representation of $H$ is obtained either by inflating an irreducible representation of $\overline{C} \times P$ or by inducing an irreducible representation of $H_1$. 
\label{serreproposition}
\end{lemma}
\noindent{\bf Proof:} Let $\chi \in \widehat{C}$ and $\chi \neq 1$. Let $g \in P$. If $g\cdot \chi =\chi$, then $\chi(ghg^{-1}) = \chi(h)$ for all $h \in C$. Hence $ghg^{-1}h^{-1} \in ker(\chi)$ for all $h \in C$, and in particular for any generator $c$ of $C$. As $\chi \neq 1$, $ker(\chi)$ is a proper subgroup of $C$. Hence if $gcg^{-1} = c^a$, then $a \equiv 1 (\text{mod }p)$. On the other hand as order of $P$ is prime to $p$, $a^{p-1} \equiv 1 (\text{mod }|C|)$. Hence $a \equiv 1 (\text{mod } |C|)$, i.e. $g \in P_1$. The second assertion now follows directly from proposition 25 in Serre \cite{Serre:representationtheory}.  \qed

For a finite group $P$ we denote by $Conj(P)$ the set of conjugacy classes of $P$. Let $n$ be a non-negative integer. We define a map
\[
\beta: \mathbb{Q}_p[Conj(\mathcal{G}/\Gamma^{p^n})] \rightarrow \mathbb{Q}_p[Conj(\overline{\mathcal{G}}/\Gamma^{p^n})] \times \mathbb{Q}_p[\mathcal{G}_1/\Gamma^{p^n}],
\]
where the map into the first component is induced by natural surjection and the map into the second component is defined as follows: Let $B$ be a set of left coset representatives of $H_1$ in $H$ and let $g \in \mathcal{G}/\Gamma^{p^n}$. Then the map is the $\mathbb{Q}_p$-linear map induced by
\[
g \mapsto \sum_{x \in B} \{ x^{-1}gx : x^{-1}gx \in \mathcal{G}_1/\Gamma^{p^n} \}.
\]
For a finite group $P$ let $R(P)$ denote the ring of virtual characters of $P$. Then $id_{\overline{\mathbb{Q}}_p} \otimes \beta$ is dual to the following map 
\[
\overline{\mathbb{Q}}_p \otimes_{\mathbb{Z}} R(\overline{\mathcal{G}}/\Gamma^{p^n}) \times \overline{\mathbb{Q}}_p \otimes_{\mathbb{Z}} R(\mathcal{G}_1/\Gamma^{p^n}) \rightarrow \overline{\mathbb{Q}}_p \otimes_{\mathbb{Z}} R(\mathcal{G}/ \Gamma^{p^n}),
\]
\[
(\chi, \rho) \mapsto Inf(\chi) + Ind(\rho).
\]
It follows from the above lemma that this map is surjective. Hence $\beta$ is injective. It induces an injection
\[
\beta : \ilim{n} \mathbb{Q}_p[Conj(\mathcal{G}/\Gamma^{p^n})] \rightarrow \ilim{n} \mathbb{Q}_p[Conj(\overline{\mathcal{G}}/\Gamma^{p^n})] \times \ilim{n} \mathbb{Q}_p[\mathcal{G}_1/\Gamma^{p^n}].
\]

We need the $log$ map on $K_1$-groups defined by Oliver and Taylor. We will discuss it in more detail in the next section. Now we just refer to Chapter 2 of Oliver \cite{Oliver:1988}. For a finite group $G$ there is a group homomorphism
\[
log: K_1(\mathbb{Z}_p[G]) \rightarrow \mathbb{Q}_p[Conj(G)].
\]
Theorem 2.9 in Oliver \cite{Oliver:1988} says that the kernel of this homomorphism is the torsion subgroup. However, by a theorem of Wall (see \emph{loc. cit.} theorem 7.4) the torsion subgroup of $K_1(\mathbb{Z}_p[G])$ is $\mu_{p-1} \times G^{ab} \times SK_1(\mathbb{Z}_p[G])$. But $SK_1(\mathbb{Z}_p[\mathcal{G}/\Gamma^{p^n}]) = 1$ by 12.7 in Oliver \cite{Oliver:1988}. Consider the exact sequences
\[
1 \rightarrow \mu_{p-1} \times (\mathcal{G}/\Gamma^{p^n})^{ab} \rightarrow K_1(\mathbb{Z}_p[\mathcal{G}/\Gamma^{p^n}]) \xrightarrow{log} \mathbb{Q}_p[Conj(\mathcal{G}/\Gamma^{p^n})].
\]
They induce an exact sequence
\[
1 \rightarrow \mu_{p-1} \times \mathcal{G}^{ab} \rightarrow K_1(\Lambda(\mathcal{G})) \xrightarrow{log} \ilim{n} \mathbb{Q}_p[Conj(\mathcal{G}/\Gamma^{p^n})]
\]
We also have a commutative diagram (see proof of theorem 6.8 in Oliver \cite{Oliver:1988})
\[
\xymatrix{ K_1(\Lambda(\mathcal{G})) \ar[rr]^{log} \ar[d]_{\theta} & & \ilim{n} \mathbb{Q}_p[Conj(\mathcal{G}/\Gamma^{p^n})] \ar[d]^{\beta} \\
K_1(\Lambda(\overline{\mathcal{G}})) \times K_1(\Lambda(\mathcal{G}_1)) \ar[rr]_{log} & & \ilim{n} \mathbb{Q}_p[Conj(\overline{\mathcal{G}}/\Gamma^{p^n})] \times \ilim{n} \mathbb{Q}_p[\mathcal{G}_1/\Gamma^{p^n}] }
\]

We are now ready to show injectivity of $\theta$. Let $x \in Ker(\theta)$. Because $\beta$ is injective, $x \in Ker(log)$. Hence $x \in \mu_{p-1} \times \mathcal{G}^{ab}$. But under the natural surjection
\[
K_1(\Lambda(\mathcal{G})) \rightarrow K_1(\Lambda(\overline{\mathcal{G}})),
\] 
(which is the first component of the map $\theta$) $\mu_{p-1} \times \mathcal{G}^{ab}$ maps identically on $\mu_{p-1} \times \overline{\mathcal{G}}^{ab}$ (note that $\mathcal{G}^{ab} = \overline{\mathcal{G}}^{ab}$). Hence $x=1$. 
\qed

\begin{remark} The above proposition is a prototype of the main results in section \ref{sectioncomputationofk1} (theorems \ref{theorem1} and \ref{theorem2}). The congruence condition in those theorems do not appear here because norm is taken in subgroups whose index in $\mathcal{G}$ is prime to $p$.
\label{remarkaboutphi}
\end{remark}

\begin{theorem} Let $F_{\infty}/F$ be an admissible extension satisfying the hypothesis $\mu=0$ and let $\mathcal{G} = Gal(F_{\infty}/F)$. Assume that $\mathcal{G}$ is $l$-$\mathbb{Q}_p$-elementary for some prime $l \neq p$. Then the main conjecture for $F_{\infty}/F$ is valid. Moreover, $K_1(\Lambda(\mathcal{G}))$ injects in $K_1(\Lambda(\mathcal{G})_S)$.
\label{lhyperelementarymainconjecture}
\end{theorem}

\noindent{\bf Proof:} We use the above notation. Let $f \in K_1(\Lambda(\mathcal{G})_S)$ be any element such that $\partial(f) = -[C(F_{\infty}/F)]$. Let $\theta_S(f) = (f_0,f_1) \in \Phi_S$. Let $L=F_{\infty}^{\mathcal{G}_1}$. Let $F_{\infty}'$ be the subfield of $F_{\infty}$ such that $Gal(F_{\infty}'/F) = \overline{\mathcal{G}}$. Then the main conjecture is true for subextensions of $F_{\infty}/L$ and the extension $F_{\infty}'/F$ by theorem \ref{primetopmainconjecture}. Let $\zeta_1$ and $\zeta_0$ be the corresponding $p$-adic zeta functions satisfying the main conjecture. By the interpolation property and uniqueness we verify that $(\zeta_0, \zeta_1) \in \Phi_S$. Let $u_i = \zeta_if_i^{-1}$ (for $i=0,1$). Then $(u_0,u_1) \in \Phi$. Let $u \in K_1(\Lambda(\mathcal{G}))$ be the unique element such that $\theta(u) =(u_0,u_1)$. Then $\zeta = uf$ is the $p$-adic zeta function we seek. It is clear that $\partial(\zeta) = -[C(F_{\infty}/F)]$. We have to show that $\zeta$ satisfies the required interpolation property. Let $\rho$ be an irreducible Artin representation of $\mathcal{G}$. Then by lemma \ref{serreproposition} $\rho$ is either obtained by inflating a representation $\rho_0$ of $\overline{\mathcal{G}}$ or by inducing a representation $\rho_1$ of $\mathcal{G}_1$.  Hence for any positive integer $r$ divisible by $[F_{\infty}(\mu_p):F_{\infty}]$, we have, for $i=0 \text{ or } 1$ 
\[
\zeta(\rho \kappa_F^r) = \zeta_i(\rho_i \kappa_i^r) = L_{\Sigma}(\rho_i, 1-r) = L_{\Sigma}(\rho, 1-r),
\]
where $\kappa_0 = \kappa_F$ and $\kappa_1 = \kappa_L$.

The injection of $K_1(\Lambda(\mathcal{G}))$ in $K_1(\Lambda(\mathcal{G})_S)$ is clear from the following diagram
\[
\xymatrix{ & K_1(\Lambda(\mathcal{G})) \ar[r] \ar[d]_{\sim} & K_1(\Lambda(\mathcal{G})_S) \ar[d] \\
0 \ar[r] & \Phi \ar[r] & \Phi_S}
\] 
The injectivity of the bottom arrow comes from theorem (\ref{primetopmainconjecture}).
\qed

\subsection{Reduction to $p$-elementary extensions} In this case we assume that $G = \mathcal{G}/\Gamma^{p^e}$ is $p$-$\mathbb{Q}_p$-elementary. Say $G = C_n \rtimes H_p$, where $H_p$ is a $p$-group and $C_n$ is a cyclic group of order $n$ prime to $p$. The Galois group $Gal(\mathbb{Q}_p(\mu_n)/\mathbb{Q}_p)$ acts on the set of one dimensional characters $\hat{C}_n$ of $C_n$. We let $C$ be the orbit set of $\hat{C}_n$ under this action. By abuse of notation we also use $C$ for a set containing exactly one representative from each orbit. Let $L_{\chi} = \mathbb{Q}_p(\chi(c)|c \in C_n)$ be a finite field extension of $\mathbb{Q}_p$. Let $O_{\chi}$ be the ring of integers of $L_{\chi}$. Then the ring $\mathbb{Z}_p[C_n]$ decomposes as
\[
\mathbb{Z}_p[C_n] \cong \oplus_{\chi \in C} O_{\chi} 
\]
The action of $H_p$ on $C_n$ induces an action on $\oplus_{\chi \in C}O_{\chi}$. But $Im(H_p \rightarrow Aut(C_n)) \subset Gal(\mathbb{Q}_p(\mu_n)/\mathbb{Q}_p)$ as $G$ is $\mathbb{Q}_p$-elementary. As $Gal(\mathbb{Q}_p(\mu_n)/\mathbb{Q}_p)$ fixes each $O_{\chi}$, the group $H_p$ acts on each $O_{\chi}$ and we get a homomorphism (for details see proposition 11.6 in Oliver  \cite{Oliver:1988})
\[
H_p \rightarrow Gal(L_{\chi}/\mathbb{Q}_p).
\]
Recall that $U_{H_p}$ denotes the inverse image of $H_p$ in $\mathcal{G}$. For this section we simply denote it by $U$. Then $\mathcal{G} = C_n \rtimes U$. Let $t_{\chi}$ denote the composition
\[
t_{\chi} : U \rightarrow H_p \rightarrow Gal(L_{\chi}/\mathbb{Q}_p).
\]
Let $U_{\chi}$ denote $ker(t_{\chi})$.

\begin{remark} We have the map
\begin{align*}
K_1'(\Lambda(\mathcal{G})) & \rightarrow K_1'(\Lambda(\Gamma^{p^e})[G]^{\tau}) \\
& \rightarrow K_1'(\Lambda(\Gamma^{p^e})[C_n][H_p]^{\tau} \\
& \rightarrow \oplus_{\chi \in C} K_1'(\Lambda_{O_{\chi}}(\Gamma^{p^e})[H_p]^{\tau}) \\ 
&\rightarrow \oplus_{\chi \in C} K_1'(\Lambda_{O_{\chi}}(U)^{\tau}) \\
& \xrightarrow{norm} \oplus_{\chi \in C} K_1'(\Lambda_{O_{\chi}}(U_{\chi})) 
\end{align*}
Here $\tau$ in the second and third line is denoting the usual twisting map and the action of $H_p$ on $C_n$ and $O_{\chi}$ respectively. The superscript $\tau$ in the fourth line denotes the action of $U$ on $O_{\chi}$. It is easy to see that the image of $K_1'(\Lambda(\mathcal{G}))$ lands in $ \oplus_{\chi \in C} K_1'(\Lambda_{O_{\chi}}(U_{\chi}))^{U/U_{\chi}}$. The argument is exactly the same as in the proof of M2) in lemma \ref{multiplicativecontain} and will not be repeated here. Similarly, we have a map 
\[
K_1'(\Lambda(\mathcal{G})_S) \rightarrow \oplus_{\chi \in C} K_1'(\Lambda_{O_{\chi}}(\Lambda(U_{\chi}))_S)^{U/U_{\chi}}.
\]
\end{remark}

\begin{proposition} We have an isomorphism
\[
K_1'(\Lambda(\mathcal{G})) \rightarrow \oplus_{\chi \in C} K_1'(\Lambda_{O_{\chi}}(U_{\chi}))^{U/U_{\chi}},
\]
\label{k1ofphyperelementary}
\end{proposition}

\noindent{\bf Proof:} For every $n \geq 0$, let $G_n = \mathcal{G}/\Gamma^{p^{e+n}}$. Let $U_{\chi, n}$ be the kernel of the map
\[
U/\Gamma^{p^{e+n}} \rightarrow H_p \rightarrow Gal(L_{\chi}/\mathbb{Q}_p).
\]
Then by proposition 11.6 and theorem 12.3 (4) in Oliver \cite{Oliver:1988}, we have
\[
K_1'(\mathbb{Z}_p[G_n]) \xrightarrow{\sim} \oplus_{\chi \in C} K_1'(O_{\chi}[U_{\chi,n}])^{U/U_{\chi,n}}.
\]
By passing to the inverse limit gives the required result by lemma \ref{inverselimitofk1}.
\qed

Define $L:=F(\mu_p) \cap F_{\infty}$. Then $L$ is a totally real number field and $[L:F]$ divides $p-1$. Hence $L \cap F_{\infty}^{C_n} = F$ and $Gal(F_{\infty}^{C_n}L/F) \cong Gal(L/F) \times U$. Then $D := Gal(F_{\infty}/F_{\infty}^{C_n}L) \subset C_n$ is a cyclic group. For every $\chi \in \hat{C}_n$ we define $C_{\chi}  := ker(\chi) \cap D$. Then $Gal(F_{\infty}^{C_{\chi}}/F_{\infty}^{C_n \rtimes U_{\chi}}) \cong C_n/C_{\chi} \times U_{\chi}$ (the only statement to be explained here is the triviality of action of $U_{\chi}$ on $C_n/C_{\chi}$. For a subgroup $J$ of $C_n$, the action of $U_{\chi}$ is trivial on $C_n/J$ if and only if for any $g \in U_{\chi}$ and any $h \in C_n$ the element $hgh^{-1}g^{-1} \in J$. Since $U_{\chi}$ acts trivially on $C_n/D$ and $C_n/ker(\chi)$, it acts trivially on $C_n/C_{\chi}$). We then have 
\[
[F_{\infty}(\mu_p):F_{\infty}] = [F_{\infty}^{C_{\chi}}(\mu_p):F_{\infty}^{C_{\chi}}].
\]

\begin{theorem} Let $F_{\infty}/F$ be an addmissible $p$-adic Lie extension satisfying $\mu=0$ hypothesis. Assume that $\mathcal{G} = Gal(F_{\infty}/F)$ is $p$-$\mathbb{Q}_p$-elementary group. With the notation as above, assume that the main conjecture is true for $F_{\infty}^{C_{\chi}}/F_{\infty}^{C_n\rtimes U_{\chi}}$ for each $\chi \in \hat{C}_n$ and that $K_1'(\Lambda_{O_{\chi}}(U_{\chi}))$ injects in $K_1'(\Lambda_{O_{\chi}}(U_{\chi})_S)$. Then the main conjecture is true for $F_{\infty}/F$ and $K_1'(\Lambda(\mathcal{G}))$ injects in $K_1'(\Lambda(\mathcal{G})_S)$. 
\label{theoremreductiontoprop}
\end{theorem}

\noindent{\bf Proof:} Assume that the main conjecture is valid for each of the extensions $F_{\infty}^{C_{\chi}}/F_{\infty}^{C_n \rtimes U_{\chi}}$. Let $\zeta_{\chi} \in K_1'(\Lambda(C_n/C_{\chi} \times U_{\chi})_S)$ be the $p$-adic zeta function in the main conjecture. Let $L_p(\chi)$ be the image of $\zeta_{\chi}$ under the natural map 
\[
K_1'(\Lambda(C_n/C_{\chi} \times U_{\chi})_S) \rightarrow K_1'(\Lambda_{O_{\chi}}(U_{\chi})_S),
\]
induced by the natural surjection $\mathbb{Z}_p[C_n/C_{\chi}] \rightarrow O_{\chi}$. Consider the following commutative diagram
\[
\xymatrix{ K_1'(\Lambda(\mathcal{G})) \ar[r] \ar[d]_{\sim} & K_1'(\Lambda(\mathcal{G})_S) \ar[d] \ar[r]^{\partial} & K_0(\Lambda(\mathcal{G}), \Lambda(\mathcal{G})_S)  \ar[d]  \\
\oplus_{\chi} K_1'(\Lambda_{O_{\chi}}(U_{\chi}))^{U/U_{\chi}} \ar[r] & \oplus_{\chi} K_1'(\Lambda_{O_{\chi}}(U_{\chi})_S)^{U/U_{\chi}} \ar[r]^{\partial} & \oplus_{\chi}K_0(\Lambda_{O_{\chi}}(U_{\chi}), \Lambda_{O_{\chi}}(U_{\chi})_S) }
\]
By uniqueness of the $p$-adic zeta function satisfying the main conjecture for the extension $F_{\infty}^{C_{\chi}}/F_{\infty}^{C_n \rtimes U_{\chi}}$, we get that $\zeta_{\chi} \in K_1'(\Lambda(C_n/C_{\chi} \times U_{\chi})_S)^{U/U_{\chi}}$. Hence $L_p(\chi) \in K_1'(\Lambda_{O_{\chi}}(U_{\chi})_S)^{U/U_{\chi}}$. Let $f \in K_1'(\Lambda(\mathcal{G})_S)$ be such that $\partial(f) = -[C(F_{\infty}/F)]$. Let the image of $f$ under the middle vertical arrow be $(f_{\chi})$. Let $u_{\chi} = L_p(\chi) f_{\chi}^{-1}$. Then $(u_{\chi})_{\chi} \in \oplus_{\chi} K_1'(\Lambda_{O_{\chi}}(U_{\chi}))^{U/U_{\chi}}$. Hence there is a unique $u \in K_1(\Lambda(\mathcal{G}))$ mapping to $(u_{\chi})_{\chi}$ under the first vertical arrow. We claim that $\zeta = uf$ is the $p$-adic zeta function we seek. It is clear that $\partial(\zeta) = -[C(F_{\infty}/F)]$. We now show the interpolation property. Let $\rho$ be an irreducible Artin representation of $\mathcal{G}$. Then by proposition 25 in Serre \cite{Serre:representationtheory} there is a $\chi$ and an Artin representation $\rho_{\chi}$ of $U_{\chi}$ such that 
\[
\rho = Ind^{\mathcal{G}}_{C_n \rtimes U_{\chi}} (\chi\rho_{\chi}).
\]
Hence for any positive integer $r$ divisible by $[F_{\infty}(\mu_p):F_{\infty}]$, we get
\[
\zeta(\rho\kappa_F^r) = \zeta_{\chi}(\chi\rho_{\chi}\kappa^r) = L_{\Sigma}(\rho, 1-r),
\]
where $\kappa$ is the $p$-adic cyclotomic character of $F_{\infty}^{C_n \rtimes U_{\chi}}$.

Uniqueness of the $p$-adic zeta function and the statement about $K_1'$ groups follows from an easy diagram chase in the above diagram.
\qed

Hence we reduce the proof of the main conjecture to the case when $\mathcal{G} = \Delta \times \mathcal{G}_p$, where $\Delta$ is a finite cyclic group of order prime to $p$ and $\mathcal{G}_p$ is a pro-$p$ $p$-adic Lie group of dimension 1.

\section{The main algebraic results}
\label{sectioncomputationofk1}
We obtain a description of $K_1$ groups of Iwasawa algebras of one dimensional compact pro-$p$ $p$-adic Lie groups and its localisation. This result will predict certain congruences between abelian $p$-adic zeta functions which we prove in the next section.

\subsection{Notations} The purpose of this subsection is to introduce the notations to be used throughout this section. Let $\mathcal{G}$ be a one dimensional pro-$p$ $p$-adic Lie group with a closed normal subgroup $H$ such that $\mathcal{G}/H = \Gamma$ is isomorphic to the additive group of $p$-adic integers. We fix a lift of $\Gamma$ in $\mathcal{G}$. This gives an isomorphism $\mathcal{G} \cong H \rtimes \Gamma$ which we take as an identification. We write $\Gamma$ multiplicatively. Let $\Gamma^{p^e}$ be a fixed open subgroup of $\Gamma$ acting trivially on $H$. Put $Z= \Gamma^{p^e}$ and $G = \mathcal{G}/Z$. Let $O$ be the ring of integers in a finite unramified extension of $\mathbb{Q}_p$ or a finite direct sum of such rings. An important example for us is $O = \mathbb{Z}_p[\Delta]$, where $\Delta$ is a finite cyclic group of order prime to $p$. The group $Gal(\overline{\mathbb{Q}}_p/\mathbb{Q}_p)$ acts on $\hat{\Delta}$, the group of characters of $\Delta$. Then $\mathbb{Q}_p(\sigma \cdot \chi) = \mathbb{Q}_p(\chi)$, for all $\sigma \in Gal(\overline{\mathbb{Q}}_p/\mathbb{Q}_p)$ and all $\chi \in \hat{\Delta}$. Let $C=H_0(Gal(\overline{\mathbb{Q}}_p/\mathbb{Q}_p),\hat{\Delta})$ be the set of orbits of the action of $Gal(\overline{\mathbb{Q}}_p/\mathbb{Q}_p)$ on $\hat{\Delta}$. Then there is an isomorphism $\mathbb{Z}_p[\Delta] \cong \oplus_{\chi \in C} O_{\chi}$, where $O_{\chi}$ is the ring of integers in $\mathbb{Q}_p(\chi)$, given by $x \mapsto (\chi(x))_{\chi \in C}$. Under this isomorphism the Frobenius automorphism of $\oplus_{\chi \in C} O_{\chi}$ corresponds to the $\mathbb{Z}_p$-linear automorphism of $\mathbb{Z}_p[\Delta]$ given by $\delta \mapsto \delta^p$ for all $\delta \in \Delta$.

The Iwasawa algebra of $\mathcal{G}$ with coefficients in $O$ is defined as
\[
\Lambda_O(\mathcal{G}) = \ilim{U} O[\mathcal{G}/U],
\]
where the inverse limit is over all open normal subgroups $U$ of $\mathcal{G}$. Clearly if $O = O_1 \oplus O_2$, then $\Lambda_O(\mathcal{G}) = \Lambda_{O_1}(\mathcal{G}) \oplus \Lambda_{O_2}(\mathcal{G})$.

\subsubsection{Twisted group rings} Recall the definition of twisted group rings. Let $R$ be a ring and $P$ be any group. Let 
\[
\tau: P \times P \rightarrow R^{\times},
\]
be a two cocycle. Then the twisted group ring, denoted by $R[P]^{\tau}$, is a free $R$-module generated by $P$. We denote the image of $h \in P$ in $R[P]^{\tau}$ by $\overline{h}$. Hence, every element of $R[P]^{\tau}$ is a finite sum $\sum_{h \in P} r_h \overline{h}$ and the addition is component wise. The multiplication has the following twist: 
\[
\overline{h}\cdot \overline{h'} = \tau(h,h') \overline{hh'}.
\]

\subsubsection{The Iwasawa algebra as a twisted group ring} Fix a topological generator $\gamma$ of $\Gamma$. Then $\gamma^{p^e}$ is a topological generator of $Z$. The Iwasawa algebra $\Lambda_{O}(\mathcal{G})$ is a twisted group ring
\[
\Lambda_{O}(\mathcal{G}) = \Lambda_{O}(Z)[G]^{\tau},
\]
where $\tau$ is the twisting map given by 
\[
\tau(h_1\gamma^{a_1}, h_2\gamma^{a_2}) = \gamma^{[\frac{a_1+a_2}{p^e}]p^e} \in \Lambda_{O}(Z)^{\times}.
\]
Here $[x]$ denotes the greatest integer less than or equal to $x$. The twisting map $\tau$ does not depend on the choice of $\gamma$. Note that for any $g, g' \in G$, 
\[
\tau(g, g') = \tau(g',g) \qquad \text{i.e. $\tau$ is \emph{symmetric}},
\]
\[
\tau(g,g^{-1})=1=\tau(g,1).
\]

\subsubsection{Definition of the Ore set $S$} \label{setS} Recall the Ore set $S$, due to Coates et. al. \cite{CFKSV:2005}, mentioned in the introduction:
\[
S(\mathcal{G},H) = \{ f \in \Lambda_O(\mathcal{G}) | \Lambda_O(\mathcal{G})/\Lambda_O(\mathcal{G})f \text{ is a finitely generated } O\text{-module}\}
\]
We will usually denote $S(\mathcal{G},H)$ by $S$ as $\mathcal{G}$ and $H$ should be clear from context. It is proven in \emph{loc. cit.} that $S$ is a multiplicatively closed left and right Ore set and does not contain any zero divisors. Hence we may localise to get the ring $\Lambda_O(\mathcal{G})_S$ which contains $\Lambda_O(\mathcal{G})$. 

\subsubsection{Some useful lemmas}

\begin{lemma} The subset of $S$ defined by $T=\Lambda_O(Z) - p\Lambda_O(Z)$ is a multiplicatively closed left and right Ore subset of $\Lambda_O(\mathcal{G})$ (note that $p$ is a uniformiser in each direct summand of $O$). The inclusion of rings $\Lambda_O(\mathcal{G})_T \rightarrow \Lambda_O(\mathcal{G})_S$ is an isomorphism.
\label{lemma1}
\end{lemma}
\noindent{\bf Proof:} It is enough to prove the result when $O$ is the ring of integers in a finite unramified extension of $\mathbb{Q}_p$. Since the group $Z$ is central in $\mathcal{G}$, it is clear that $T$ is a left and right Ore set. Since $\Lambda_O(Z)$ is a domain, the set $T$ does not contain any zero divisors. Hence, the map $\Lambda_O(\mathcal{G})_T \rightarrow \Lambda_O(\mathcal{G})_S$, induced by the inclusion $T \rightarrow S$, is an inclusion. We now prove that it is surjective.

Note that $\Lambda_O(\mathcal{G})_{T} = \Lambda_O(Z)_T \otimes_{\Lambda_O(Z)} \Lambda_O(\mathcal{G})$. We first show that 
\[
Q(\Lambda_O(\mathcal{G})) = Q(\Lambda_O(Z)) \otimes_{\Lambda_O(Z)} \Lambda_O(\mathcal{G}),
\] 
where $Q(R)$ denotes the total ring of fractions of a ring $R$. Note that we have an injective map
\[
Q(\Lambda_O(Z)) \otimes_{\Lambda_O(Z)} \Lambda_O(\mathcal{G}) \hookrightarrow Q(\Lambda_O(\mathcal{G})).
\]
As $Q(\Lambda_O(Z))$ is a field and $\Lambda_O(\mathcal{G})$ is a free $\Lambda_O(Z)$-module of finite rank, the ring $Q(\Lambda_O(Z)) \otimes_{\Lambda_O(Z)} \Lambda_O(\mathcal{G})$ is an Artinian ring. Hence every regular element is invertible. The ring $\Lambda_O(\mathcal{G})$ is contained in $Q(\Lambda_O(Z)) \otimes_{\Lambda_O(Z)} \Lambda_O(\mathcal{G})$ and every regular element of $\Lambda_O(\mathcal{G})$ is invertible in $Q(\Lambda_O(Z)) \otimes_{\Lambda_O(Z)} \Lambda_O(\mathcal{G})$, hence the injection $Q(\Lambda_O(Z)) \otimes_{\Lambda_O(Z)} \Lambda_O(\mathcal{G}) \hookrightarrow Q(\Lambda_O(\mathcal{G}))$ must be surjective. Any element $x \in \Lambda_O(\mathcal{G})_S \subset Q(\Lambda_O(\mathcal{G}))$ can be written as $\frac{a}{t}$, with $a \in \Lambda_O(\mathcal{G})$ and $ t $ a non-zero element of $\Lambda_O(Z)$. If $t \in p^n\Lambda_O(Z)$, then $tx = a \in p^n\Lambda_O(\mathcal{G})_S$. On the other hand, $a$ also lies in $\Lambda_O(\mathcal{G})$, hence $a \in p^n\Lambda_O(\mathcal{G})$. Then we can divide the largest possible power of $p$ from $t$ and the same power of $p$ from $a$, and $x$ can be presented as $\frac{a}{t}$, with $a \in \Lambda_O(\mathcal{G})$ and $ t \in T$. \qed

\begin{remark} We remark that $\Lambda_O(\mathcal{G})_T = \Lambda_O(Z)_T[G]^{\tau}$, for the same twisting map $\tau$ as above. We also study the $p$-adic completion $\widehat{\Lambda_{O}(\mathcal{G})_{T}} = \widehat{\Lambda_O(Z)_T}[G]^{\tau}$ of $\Lambda_O(\mathcal{G})_T$. Note that $\Lambda_O(Z)_T = \Lambda_O(Z)_{(p)}$.
\end{remark}

\begin{notation} In the rest of this section $R$ denotes either $\Lambda_O(Z)$ or $\widehat{\Lambda_O(Z)_{(p)}}$. For any subgroup $P$ of $G$, we denote the inverse image of $P$ in $\mathcal{G}$ by $U_P$. Let $N_GP$ be the normaliser of $P$ in $G$. We put $W_GP= N_GP/P$, the Weyl group. We denote the set of all cyclic subgroups of $G$ by $C(G)$.
\end{notation}

\begin{notation} Let $P \leq P' \leq G$. Then we denote by $ver^{P'}_P$ the transfer homomorphism from $U_{P'}^{ab} \rightarrow U_{P}^{ab}$. The homomorphism $ver^{P'}_P$ induces maps on various Iwasawa algebras by acting through Frobenious on the coefficients. We denote this induced map again by $ver^{P'}_P$. 

Let $P$ be a cyclic subgroup of $G$. We denote the transfer map from $U_{P}$ to $U_{P^p}$ simply by $\varphi$. If $P$ is a cyclic subgroup of $G$ and $P'$ is any subgroup of $G$ containing $P^p$, then we denote by $\varphi$ the composition
\[
U_P \xrightarrow{\varphi} U_{P^p} \hookrightarrow U_{P'}.
\]
$\varphi$ induces maps on various Iwasawa algebras by acting through Frobenius on the coeffiecients. We denote this induced maps again by $\varphi$.
\end{notation}

\begin{definition} Let $I_1$ and $I_2$ be ideals in $R[G]^{\tau}$. Then we define $[I_1, I_2]$ to be the additive subgroup of $R[G]^{\tau}$ generated by elements of the form $ab-ba$, where $a \in I_1$ and $b \in I_2$.
\end{definition} 

\begin{remark} $R[Conj(G)]^{\tau}$ is just a $R$-module and multiplication is not defined. However, taking powers of elements in $Conj(G)$ is well defined and while doing so we must remember the twisting map $\tau$. Hence we put it in the notation.
\end{remark}

\begin{lemma} We have the isomorphism of $R$-modules
\[
R[G]^{\tau}/[R[G]^{\tau}, R[G]^{\tau}] \rightarrow R[Conj(G)]^{\tau}.
\]
\end{lemma}

\begin{notation} Throughout this paper we usually denote the class of $g$ in $Conj(G)$ by $[g]$.
\end{notation}

\noindent{\bf Proof:} Consider the $R$-module homomorphism 
\[
R[G]^{\tau} \xrightarrow{\phi} R[Conj(G)]^{\tau},
\]
\[
\sum r_g \overline{g} \mapsto \sum r_g [\overline{g}].
\]
This map $\phi$ is surjective and since $\tau$ is symmetric, the kernel of $\phi$ contains $[R[G]^{\tau}, R[G]^{\tau}]$. We must show that it is equal to $[R[G]^{\tau}, R[G]^{\tau}]$. Let $\sum r_g\overline{g} \in ker(\phi)$. Let $C_g$ denote the centraliser of $g$ in $G$.  Then for each $g \in G$, we get
\[
\sum_{x \in G/C_g} r_{xgx^{-1}} = 0.
\]
Consider
\begin{align*}
\sum r_{xgx^{-1}} \overline{xgx^{-1}} - \sum r_{xgx^{-1}} \overline{g} & = \sum r_{xgx^{-1}}(\overline{xgx^{-1}} - \overline{g}) \\
& = \sum r_{xgx^{-1}} (\tau(xg, x^{-1})^{-1} \overline{xg} \overline{x^{-1}} - \tau(x^{-1}, xg)^{-1} \overline{x^{-1}}\overline{xg}) \\
& = \sum r_{xgx^{-1}} \tau(xg, x^{-1})^{-1} ( \overline{xg}\overline{x^{-1}} - \overline{x^{-1}} \overline{xg}),
\end{align*}
which clearly belongs to $[R[G]^{\tau}, R[G]^{\tau}]$. All the sums are over $x \in G/C_g$. \qed

\begin{lemma} For any subgroup $P \leq G$, we have
\[
\Lambda_O(U_P) \cong \Lambda_O(Z)[P]^{\tau},
\]
and
\[
\widehat{\Lambda_O(U_P)_S} \cong \widehat{\Lambda_O(Z)_{(p)}}[P]^{\tau},
\]
for the same twisting map $\tau$ as above.
\end{lemma}

\begin{lemma} If $P \leq G$ is a cyclic subgroup, then $U_P$ is abelian (though $U_P$ is not necessarily a direct product of $Z$ and $P$). 
\end{lemma}

\begin{definition} For any $P \in C(G)$ with $P \neq \{1\}$, we choose and fix a non-trivial character $\omega_P$ of $P$ of order $p$. We put $\omega_1 := \omega_{\{1\}}$ to be a non-trivial character of $Z$ whose restriction to $Z^p$ is trivial. Then $\omega_P$ induces a map on $\Lambda_{O[\mu_p]}(U_P)^{\times}$, on $\Lambda_{O[\mu_p]}(U_P)_S^{\times}$, and on $\widehat{\Lambda_{O[\mu_p]}(U_P)_S}^{\times}$ given by mapping $g \in U_P$ to $\omega_P(g)g$. We denote this map by $\omega_P$ again. We define a map $\alpha_P$, for $P \neq \{1\}$, from $\Lambda_O(U_P)^{\times}$ to itself or from $\Lambda_O(U_P)_S^{\times}$ to itself or from $\widehat{\Lambda_O(U_P)_S}^{\times}$ to itself by 
\[
\alpha_P(x) = \frac{x^p}{\prod_{k=0}^{p-1}\omega_P^k(x)}.
\]
If $P \leq G$ is not cyclic or if $P = \{1\}$, we put $\alpha_P(x) = x^p$. We put $\alpha = (\alpha_P)_{P \leq G}$.
\label{defnalpha}
\end{definition}

\begin{definition} For two subgroup $P \leq P' \leq G$ such that $[P',P'] \leq P$, we have maps
\[
tr^{P'}_P : \Lambda_O(U_{P'}^{ab}) \rightarrow \Lambda_O(U_P/[U_{P'},U_{P'}]) \qquad \text{the trace map},
\]
\[
nr^{P'}_P : \Lambda_O(U_{P'}^{ab})^{\times} \rightarrow \Lambda_O(U_P/[U_{P'},U_{P'}])^{\times} \qquad \text{the norm map},
\]
and
\[
\pi^{P'}_P: \Lambda_O(U_P^{ab}) \rightarrow \Lambda_O(U_P/[U_{P'},U_{P'}]) \qquad \text{the natural surjection}.
\]
Similarly, we have trace, norm and the natural surjection maps in the localised case and the case of $p$-adic completions. We again denote them by $tr^{P'}_P$, $nr^{P'}_P$ and $\pi^{P'}_P$ respectively. We also note that if $P$ and $P'$ are cyclic then $\pi^{P'}_P$ is identity.
\end{definition}

\begin{definition}Let $P \leq P' \leq G$ be such that $P$ is a normal subgroup of $P'$. Then define $T_{P,P'}$ (resp. $T_{P,P',S}$ and $\widehat{T_{P,P'}}$) to be the image of the map 
\[
\Lambda_O(U_P) \rightarrow \Lambda_O(U_P)
\]
\[
\text{(resp. } \Lambda_O(U_P)_S \rightarrow \Lambda_O(U_P)_S \qquad \text{and} \qquad \widehat{\Lambda_O(U_P)_S} \rightarrow \widehat{\Lambda_O(U_P)_S} \text{)},
\]
given by $x \mapsto \sum_{g \in P'/P} \tilde{g}x\tilde{g}^{-1}$ ($\tilde{g}$ is a lift of $g$). We denote $T_{P,N_GP}$ (resp. $T_{P,N_GP,S}$ and $\widehat{T_{P,P'}}$) simply by $T_P$ (resp. $T_{P,S}$ and $\widehat{T_P}$).

\label{defntraceideal}
\end{definition}

\begin{lemma} Let $P \leq P' \leq G$ be two subgroups such that $[P':P]=p$. Let $\omega$ be a nontrivial character of $P'/P$. Consider $\omega$ as a character of $U_{P'}^{ab}$ using the surjection $U_{P'}^{ab} \rightarrow P'/P$. Let $R$ be either $\Lambda_O(Z)$, $\Lambda_O(Z)_T$ or $\widehat{\Lambda_O(Z)_T}$. Then the norm and trace map
\[
tr^{P'}_P : R[P'^{ab}]^{\tau} \rightarrow R[P/[P',P']]^{\tau},
\]
\[
nr^{P'}_P : (R[P'^{ab}]^{\tau})^{\times} \rightarrow (R[P/[P',P']]^{\tau})^{\times}
\]
are given by 
\[
tr^{P'}_P(x) = \sum_{k=0}^{p-1} \omega^k(x)
\]
\[
nr^{P'}_P(x) = \prod_{k=0}^{p-1}\omega^k(x).
\]
Here $\omega^k$ is the map from $O[\mu_p] \otimes_O R[P'^{ab}]^{\tau}$ to itself given by mapping $g \in U_{P'}^{ab}$ to $\omega^k(g)g$. 
\label{normandtrace}
\end{lemma}  
 
\noindent{\bf Proof:} Let $P'/P = \{1,g,\ldots, g^{p-1}\}$. Then $R[P'^{ab}]^{\tau}$ is a free $R[P/[P',P']]^{\tau}$-module with basis $\{1,g,\ldots, g^{p-1}\}$. Let $x= a_0 + a_1g + \cdots + a_{p-1}g^{p-1} \in R[P'^{ab}]^{\tau}$ with each $a_i \in R[P/[P',P']]^{\tau}$. The matrix (with respect to the basis $\{1,g,\ldots,g^{p-1}\}$) of the $R[P/[P',P']]^{\tau}$-linear map on $R[P'^{ab}]^{\tau}$ defined by multiplication by $x$ is
\[
\left( \begin{array}{cccc}
a_0 & a_1 & \cdots & a_{p-1} \\
g^pa_{p-1} & a_0 & \cdots & a_{p-2} \\
\vdots & & & \vdots \\
g^pa_1 & g^p a_2 & \cdots & a_0 
\end{array}
\right).
\]
For each $0 \leq k \leq p-1$, $(1,\omega^k(g)g, \cdots, \omega^k(g^{p-1})g^{p-1})$ is an eigenvector of the above matrix with eigenvalue $\omega^k(x)$. Hence the lemma. \qed

\subsection{Statements of the main algebraic theorems}
\subsubsection{The maps $\theta^G$ and $\theta^G_S$}

 Recall that for every subgroup $P$ of $G$, we denote by $U_P$ the inverse image of $P$ in $\mathcal{G}$. Then there is a map $\theta^G_{O,P}$ given by the composition
\[
\theta^G_{O,P}: K_1'(\Lambda_O(\mathcal{G})) \xrightarrow{norm} K_1'(\Lambda_O(U_P)) \rightarrow K_1'(\Lambda_O(U_P^{ab})) = \Lambda_O(U_P^{ab})^{\times},
\]
where the second map is induced by natural surjection. We remark that if $P \in C(G)$, then $U_P$ is abelian and $\theta^G_{O,P}$ is just the norm map.
Let $\theta^G_O$ be the map
\[
\theta^G_O = (\theta^G_{O,P})_{P \leq G} : K_1'(\Lambda_O(\mathcal{G})) \rightarrow \prod_{P \leq G} \Lambda_O(U_P^{ab})^{\times}.
\]

Similarly, we have maps $\theta^G_{O,S}$ and $\widehat{\theta}^G_O$
\[
\theta^G_{O,S} : K_1'(\Lambda_O(\mathcal{G})_S) \rightarrow \prod_{P \leq G} \Lambda_O(U_P^{ab})_S^{\times}.
\]
\[
\widehat{\theta}^G_O : K_1'(\widehat{\Lambda_O(\mathcal{G})_S}) \rightarrow \prod_{P \leq G} \widehat{\Lambda_O(U_P^{ab})_S}^{\times}. 
\]
We denote $\theta^G_O$, $\theta^G_{O,S}$ and $\widehat{\theta}^G_O$ by $\theta^G$, $\theta^G_S$ and $\widehat{\theta}^G$ respectively if $O$ is clear from the context. Our main theorem shows that $\theta^G$ is injective, describes its image and shows that image of $\theta^G_S$ intersected with $\prod_{P \leq G} \Lambda_O(U_P)^{\times}$ is exactly the image of $\theta^G$.

\subsubsection{Definition of the subgroups $\Phi^G$, $\Phi^G_S$ and $\widehat{\Phi^G}$}

The group $\mathcal{G}$ acts on the set $\prod_{P \leq G} U_P^{ab}$ by conjugation. Since $Z$ is central in $\mathcal{G}$, we get an induced action of $G$ on the set $\prod_{P \leq G} U_P^{ab}$. We use it in the following definition.

\begin{definition} Let $\Phi^G_O$ (resp. $\Phi^G_{O,S}$ and $\widehat{\Phi^G_{O,S}}$) be the subgroup of $\prod_{P\leq G}\Lambda_O(U_P^{ab})^{\times}$ (resp. $\prod_{P\leq G} \Lambda_O(U_P^{ab})_S^{\times}$ and $\prod_{P \leq G}\widehat{\Lambda_O(U_P^{ab})_S}^{\times}$) consisting of tuples $(x_P)$ satisfying \\
M1. For any $P \leq P' \leq G$ such that $[P',P'] \leq P$, we have 
\[
nr^{P'}_P(x_{P'}) = \pi^{P'}_P(x_P).
\]
M2. $(x_P)$ is fixed under conjugation action by every $g \in G$. \\
M3. For every $P \leq P' \leq G$ such that $[P':P]=p$, we have 
\[
ver^{P'}_P(x_{P'}) \equiv x_P (\text{mod } T_{P,P'}) (\text{resp. } T_{P,P',S} \text{ and } \widehat{T_{P,P'}}).
\]
M4. For every $P \in C(G)$ and $P \neq \{1\}$ we have the following congruence
\[
\alpha_P(x_P) \equiv   \prod_{P'}\varphi(\alpha_{P'}(x_{P'})) (\text{mod } pT_P) \ (\text{resp. } pT_{P,S} \text{ and } p\widehat{T_P}),
\]
where the product runs through all $P' \in C(G)$ such that $P'^p = P$ and $P' \neq P$. For $P=\{1\}$ we have the congruence
\[
x_{\{1\}}^p \equiv \varphi(x_{\{1\}}) \prod_{P'} \varphi(\alpha_{P'}(x_{P'})) (\text{mod }pT_{\{1\}}) \ (\text{resp. } pT_{\{1\},S} \text{ and } p\widehat{T_{\{1\}}}),
\]
where the product runs through all $P' \in C(G)$ such that $P' \neq \{1\}$ and $P'^p = \{1\}$. We denote $\Phi^G_O$ (resp. $\Phi^G_{O,S}$ and $\widehat{\Phi^G_{O,S}}$) by $\Phi^G$ (resp. $\Phi^G_S$ and $\widehat{\Phi^G_S}$) if $O$ is clear from the context. 

\label{defnofphi}
\end{definition} 

\subsubsection{The theorems}
\begin{theorem}The map $\theta^G$ induces an isomorphism 
\[
K_1'(\Lambda_O(\mathcal{G})) \xrightarrow{\sim} \Phi^G.
\]
\label{theorem1}
\end{theorem}

\begin{theorem} The image of $\theta^G_S$ is contained in $\Phi^G_S$. Hence 
\[
\Phi^G_S \cap \prod_{P \leq G} \Lambda_O(U_P^{ab})^{\times} = Im(\theta^G).
\]
\label{theorem2}
\end{theorem}

\begin{remark} Compare definition \ref{defnofphi} and theorems \ref{theorem1} and \ref{theorem2} with definition \ref{definitionofphi} and proposition \ref{propprimetopk1}. See also remark \ref{remarkaboutphi} after proof of proposition \ref{propprimetopk1}.
\end{remark}

Note that if $O = O_1 \oplus \cdots \oplus O_n$, then under the isomorphism 
\[
\prod_{P\leq G}\Lambda_O(U_P^{ab})^{\times} \cong \oplus_{i=1}^n \prod_{P\leq G}\Lambda_{O_i}(U_P^{ab})^{\times}
\]
(resp. 
\[
\prod_{P\leq G} \Lambda_O(U_P^{ab})_S^{\times} \cong \oplus_{i=1}^n \prod_{P\leq G} \Lambda_{O_i}(U_P^{ab})_S^{\times}
\]
and 
\[
\prod_{P \leq G}\widehat{\Lambda_O(U_P^{ab})_S}^{\times} \cong \oplus_{i=1}^n \prod_{P \leq G}\widehat{\Lambda_{O_i}(U_P^{ab})_S}^{\times}),
\]
$\Phi^G_O$ (resp. $\Phi^G_{O,S}$ and $\widehat{\Phi^G_{O,S}}$) maps to $\oplus_{i=1}^n \Phi^G_{O_i}$ (resp. $\oplus_{i=1}^n \Phi^G_{O_i,S}$ and $\oplus_{i=1}^n \widehat{\Phi^G_{O_i,S}}$). Hence for the proofs of theorems \ref{theorem1} and \ref{theorem2} (i.e in rest of section \ref{sectioncomputationofk1}) we may and do assume that $O$ is a ring of integers in a finite unramified extension of $\mathbb{Q}_p$.

\subsection{An additive theorem}

\subsubsection{The Statement} For every subgroup $P$ of $G$, the Weyl group $W_GP=N_GP/P$ acts $R$-linearly on $R[P^{ab}]^{\tau}$ by conjugation on $P$.  For any $P \leq G$, define a map 
\[
t^G_P: R[Cong(G)]^{\tau} \rightarrow R[P^{ab}]^{\tau}
\]
as follows: let $C(G,P)$ denote any set of left coset representatives of $P$ in $G$. Then 
\[
t^G_P(\overline{g}) = \sum_{x \in C(G,P)} \{ (\overline{x}^{-1})(\overline{g})(\overline{x}) |  x^{-1}gx \in P\}.
\]
This is a well-defined $R$-linear map, independent of the choice of $C(G,P)$. For any $P \in C(G)$, define
\[
\eta_P : R[P]^{\tau} \rightarrow R[P]^{\tau},
\]
by $R$-linearly extending the map 
\[
\eta_P(h)=  \left\{ 
\begin{array}{l l}
h & \quad \mbox{if $h$ is a generator of $P$} \\
0 & \quad \mbox{if not}. \\
\end{array} \right.
\]
In other words $\eta_P(x) = x - \frac{1}{p} \sum_{k=0}^{p-1} \omega_P^k(x)$.

\begin{definition} Define $\beta^G_P : R[Conj(G)]^{\tau} \rightarrow R[P^{ab}]^{\tau}$ by 
\[
\beta^G_P = \left\{
\begin{array}{l l}
\eta_P \circ t^G_P & \text{if $P \in C(G)$} \\
t^G_P & \text{if $P \leq G$ is not cyclic} 
\end{array}
\right.
\]
and $\beta^G_R$ by 
\[
\beta^G_R = (\beta^G_P)_{P\leq G} : R[Conj(G)]^{\tau} \rightarrow \prod_{P\leq G}R[P^{ab}]^{\tau}.
\]
Sometimes we denote $\beta^G_{\Lambda_O(Z)}$ by just $\beta^G$.
\end{definition}

\begin{definition} Let $P \in C(G)$ be a cyclic subgroup of $G$. We define $T_{P,R}$ to be the image of the map 
\[
tr: R[P]^{\tau} \rightarrow R[P]^{\tau} \qquad x \mapsto \sum_{g \in W_GP} (\overline{g})( \overline{x})( \overline{g}^{-1}).
\]
It is an ideal in the ring $(R[P]^{\tau})^{W_GP}$. Hence $T_{P, \Lambda_O(Z)}=T_P$ and $T_{P, \widehat{\Lambda_O(Z)_{(p)}}}=\widehat{T_P}$.
\end{definition}

If $P \leq P' \leq G$ are such that $[P':P]=p$ then $P$ is normal in $P'$ and the commutator subgroup $[P',P']$ is contained in $P$. In addition, if $P$ is a nontrivial cyclic group then $[P',P']$ is a proper subgroup of $P$. In this case the map $\eta_P$ descends to a map on $R[P/[P',P']]^{\tau}$, which we again denote by $\eta_P$, such that the following diagram commutes
\[
\xymatrix{ R[P]^{\tau} \ar[r]^{\eta_P} \ar[d]_{\pi^{P'}_{P}} & R[P]^{\tau} \ar[d]^{\pi^{P'}_P} \\
R[P/[P',P']] \ar[r]_{\eta_P} & R[P/[P',P']]} 
\]
We use this map in the following definition.

\begin{definition} Let $\psi^G_R \subset \prod_{P \leq G} R[P^{ab}]^{\tau}$ be the subgroup consisting of all tuples $(a_P)$ such that \\
A1. Let $P \leq P' \leq G$ such that $[P', P'] \leq P$ and if $P$ is a non-trivial cyclic group then $[P',P'] \neq P$. If $P$ is not cyclic, then $tr^{P'}_P(a_{P'}) = \pi^{P'}_P(a_P)$. If $P$ cyclic and $P'$ is not cyclic then $\eta_P(tr^{P'}_P(a_{P'})) = \pi^{P'}_P(a_P)$. If $P'$ is cyclic, then $tr^{P'}_P(a_{P'})=0$. \\
A2. $(a_P)_{P\in C(G)}$ is invariant under conjugation action by every $g \in G$. \\
A3. For all $P \in C(G)$, $a_P \in T_{P,R}$. \\
Sometimes we denote $\psi^G_{\Lambda_O(Z)}$ by just $\psi^G$.
\end{definition}

\begin{theorem} The homomorphism $\beta^G_R$ induces an isomorphism between $R[Conj(G)]^{\tau}$ and $\psi^G_R$.
\label{additivetheorem}
\end{theorem}

\begin{remark} The injectivity of $\beta^G_R$ follows directly from Artin's induction theorem (theorem 17 in Serre \cite{Serre:representationtheory}) which says that a linear representation of a finite group is a $\mathbb{Q}$-linear combination of representations induced from cyclic subgroups (note that $R[Cong(G)]^{\tau}$ is torsion free). Therefore in fact Artin's induction theorem gives injectivity of $\beta^G_R$ composed with the projection $\prod_{P \leq G} R[P^{ab}]^{\tau} \rightarrow \prod_{P \in C(G)} R[P]^{\tau}$. We prove this injectivity directly. The sufficiency of using cyclic subgroups of $G$ is reflected in our proof in the definition of $\delta$ below. Hence on the additive side it is enough to work with cyclic subgroups of $G$ but on the multiplicative side (i.e. on $K_1$) cyclic subgroups of $G$ are not enough because $K_1$ has non-trivial torsion subgroup. 
\end{remark}

\subsubsection{The Proof} 

\begin{lemma} The image of $\beta^G_R$ is contained in $\psi^G_R$. 
\label{additivecontain}
\end{lemma}
\noindent{\bf Proof:} It is enough to show that $\beta^G_R([\overline{g}]) \in \psi^G_R$, for any $g \in G$, i.e. it satisfies A1, A2 and A3. \\

\noindent A1. Let $P \leq P' \leq G$. First note that $[P',P'] \leq P$ implies that $P$ is normal in $P'$. We first show that the following diagram commutes.
\[
\xymatrix{ R[Conj(G)]^{\tau} \ar[r]^{t^G_{P'}} \ar[d]_{t^G_P} & R[(P')^{ab}]^{\tau} \ar[d]^{tr^{P'}_P} \\
R[P^{ab}]^{\tau} \ar[r]_{\pi^{P'}_P} & R[P/[P',P']]^{\tau}}
\]
Let $h \in P'$, then $tr^{P'}_P(\overline{h}) = 0$ unless $h \in P$ in which case it is 
\[
tr^{P'}_P(\overline{h}) = [P':P]\overline{h} \in R[P/[P',P']]^{\tau}.
\]
On the other hand $\pi^{P'}_P(t^{P'}_P([h])) = 0$ unless $h \in P$ in which case it is 
\[
\pi^{P'}_P(t^{P'}_P([h])) = [P':P]\overline{h} \in R[P/[P',P']]^{\tau}.
\]
Hence for any $g \in G$, we have
\begin{align*}
tr^{P'}_P(t^G_{P'}([g])) & = [P':P]\sum_{x \in C(G, P')} \{\overline{[x^{-1}gx]} : x^{-1}gx \in P\} \\
& = \pi^{P'}_P(t^G_P([g])).
\end{align*} 
This shows A1 when $P$ noncyclic. Now if $P$ is cyclic and $P'$ is not, then the above proof shows that the following diagram commutes
\[
\xymatrix{ R[Conj(G)]^{\tau} \ar[r]^{\beta^G_{P'}} \ar[d]_{t^G_P} & R[(P')^{ab}]^{\tau} \ar[d]^{tr^{P'}_P} \\
R[P]^{\tau} \ar[r]_{\pi^{P'}_{P}} & R[P/[P',P']]}
\]
i.e. $tr^{P'}_P(\beta^G_{P'}([g])) = \pi^{P'}_P(t^G_P([g]))$. Applying $\eta_P$ to both side gives
\begin{align*} 
\eta_P(tr^{P'}_P(\beta^G_{P'}([g]))) & = \eta_P(\pi^{P'}_P(t^G_P([g]))) \\
\eta_P(tr^{P'}_P(\beta^G_{P'}([g]))) & = \pi^{P'}_P(\beta^G_P([g]))
\end{align*}
Lastly, if $P < P'$ are both cyclic then it is clear that $tr^{P'}_P([\bar{h}]) = 0$ unless $ h \in P$, in which case $tr^{P'}_P([\bar{h}])  = [P':P][\bar{h}]$. Since the coefficient of any element $g' \in P'$ in $\beta^G_{P'}([\overline{g}])$ is 0 if $g'$ does not generate $P'$, it is clear that $tr^{P'}_P(\beta^G_{P'}([\overline{g}])) = 0$.  \\
 
\noindent A2. For any $ g, g_1 \in G$, we must show that $\overline{g_1} \beta^G_P([\overline{g}]) \overline{g_1^{-1}} = \beta^G_{g_1Pg_1^{-1}}([\overline{g}])$, for any $P \in C(G)$. 

\begin{align*} 
\overline{g_1}t^G_P([\overline{g}])\overline{g_1^{-1}} & = \overline{g_1}(\sum_{x \in C(G,P)} \{ [\overline{x^{-1}gx}] : x^{-1}gx \in P\}) \overline{g_1^{-1}} \\
& = \sum_{x \in C(G,P)} \{ [\overline{ (g_1x^{-1}g_1^{-1})(g_1gg_1^{-1})(g_1xg_1^{-1})}] : x^{-1}gx \in P\} \\
& = \sum_{x_1 \in C(G, g_1Pg_1^{-1})} \{[\overline{x_1^{-1}(g_1gg_1^{-1})x_1}] : x_1^{-1}g_1gg_1^{-1}x_1 \in g_1Pg_1^{-1}\} \\
& = t^G_{g_1Pg_1^{-1}}([\overline{g_1gg_1^{-1}}]) \\
& = t^G_{g_1Pg_1^{-1}}([\overline{g}]).
\end{align*}
Note that on several occasions above we have used $\bar{g_1} \bar{g} \overline{g_1^{-1}} = \overline{g_1gg_1^{-1}}$. A2 now follows easily. \\

\noindent A3. We must show that $t^G_P([\bar{g}]) \in T_{P,R}$ for any $P \in C(G)$, since $\eta_P(T_{P,R}) \subset T_{P,R}$ for any nontrivial cyclic subgroup $P$. 
\[
t^G_P([\overline{g}]) = t^{N_GP}_P(t^G_{N_GP}([\overline{g}])).
\]
Hence, it is enough to show that $t^{N_GP}_P([\overline{g}]) \in T_{P,R}$ for any $g \in N_GP$. But $t^{N_GP}_P([\overline{g}])$ is non-zero if and only if $ g \in P$, and when $ g \in P$, we have
\begin{align*}
t^{N_GP}_P([\overline{g}]) & = \sum_{x \in C(N_GP, P)} \overline{x^{-1}gx} \\
& = \sum_{x \in W_GP} \overline{x^{-1}gx} \in T_{P,R}.
\end{align*}
This finishes proof of the lemma. \qed

\begin{definition} We define a left inverse $\delta$  of $\beta^G_R$ by
\[
\delta :  \prod_{P\leq G} R[P^{ab}]^{\tau} \rightarrow R[Conj(G)]^{\tau}[\frac{1}{p}],
\]
by putting $\delta = \sum_{P \leq G} \delta_P$ and defining $\delta_P$ to be 0 if $P$ is not cyclic and for $P \in C(G)$ by
\[
\delta_P : R[P^{ab}]^{\tau} \rightarrow R[Conj(G)]^{\tau}[\frac{1}{p}],
\]
\[
x \in R[P^{ab}]^{\tau} \mapsto \frac{1}{[G:P]}[x] \in R[Conj(G)]^{\tau}[\frac{1}{p}].
\]
\end{definition}

\begin{lemma} $\delta \circ \beta^G_R$ is identity on $R[Conj(G)]^{\tau}$. In particular, $\beta^G_R$ is injective.
\end{lemma}

\noindent{\bf Proof:} For any $g \in G$, we show that $\delta(\beta^G_R([\overline{g}])) = [\overline{g}]$. Let $P$ be the cyclic subgroup of $G$ generated by $g$. Let $C$ be the set of all conjugates of $P$ in $G$. Then
\begin{align*} 
\delta(\beta^G_R([\overline{g}])) & = \sum_{P' \in C} \delta_{P'}(\beta^G_R([\overline{g}])) \\
&= \sum_{P'\in C} \frac{1}{[G:P]} [\beta^G_R([\overline{g}])]\\
&= \frac{1}{[G:P]}\sum_{P' \in C} [N_GP':P'] [\overline{g}] \\
&= \frac{1}{[G:N_GP]} \sum_{P' \in C} [\overline{g}] = [\overline{g}]
\end{align*}

\qed

\begin{lemma} The restriction of $\delta$ to the subgroup $\psi^G_R$ is injective and its image lies in $R[Conj(G)]^{\tau}$.
\end{lemma} 

\noindent{\bf Proof:} Let $(a_P) \in \psi^G_R$ be such that $\delta((a_P)) = 0$. We claim that $\delta_P(a_P) = 0$ for each $P \in C(G)$. This follows from two simple observations: firstly, by A1 $\delta_P(a_P)$ and $\delta_{P'}(a_{P'})$ cannot cancel each other unless $P$ and $P'$ are conjugates; but when $P$ and $P'$ are conjugates, $\delta_P(a_P) = \delta_{P'}(a_{P'})$ by A2. Hence $\delta_P(a_P) =0$ for every $P \in C(G)$. 

Let $P \in C(G)$ and $a_P = \sum_{g \in P} r_g \overline{g}$. Then $\delta_P(a_P) = \frac{1}{[G:P]} \sum_{g \in P} r_g [\overline{g}]$. Let $H_0(N_GP,P)$ be the orbit set for the conjugation action of $N_GP$ on $P$. Then
\[
\delta_P(a_P) = \sum_{x \in H_0(N_GP,P)} \sum_{g \in x} r_g[\overline{g}].
\]
If $g, g' \in x \in H_0(N_GP,P)$, then $r_g = r_{g'}$ by A2. Call it $r_x$. Hence, 
\[
\delta_P(a_P) = \sum_{x \in H_0(N_GP, P)} r_x \sum_{g \in x} [\overline{g}] = 0.
\]
Hence $r_x = 0$ for all $x \in H_0(N_GP, P)$. Hence we get that $a_P =0$ for every $P \in C(G)$. 

Next we show that $a_P =0$ for every $P \leq G$. We do this by using induction on order of $P$. Assume that $P'$ is a subgroup of $G$ of smallest order such that $a_{P'} \neq 0$. Then $P'$ cannot be cyclic. Let $a_{P'}  = \sum_{h \in (P')^{ab}} a_hh$. Let $h_0 \in (P')^{ab}$ be such that $a_{h_0} \neq 0$. Let $\tilde{h}_0$ be any lift of $h_0$ to $P'$. Let $P$ be a maximal subgroup of $P'$ containing $\tilde{h}_0$. Then $[P':P]=p$. We have to consider two cases: \\
Case1: when $P$ is cyclic. We claim that $\tilde{h}_0$ is a generator of $P$. Since $P$ is cyclic and $P'/P \cong C_p$, the cyclic group of order $p$, the group $P'$ is isomorphic to $P \rtimes C_p$. If $\tilde{h}_0$ is not a generator of $P$ then $\tilde{h}_0 \in P^p$. Take $P'' = P^p \rtimes C_p$. If $P^p = \{1\}$, then $\tilde{h}_0 = 1$ and $P' = P \times C_p$. Hence $tr^{P'}_{\{1\}}(a_{P'}) = p^2a_{h_0} \neq 0$ and by A1 applied to $\{1\} \leq P' \leq G$ we get , a contradiction to the fact that $a_{\{1\}} = 0$. If $P^p \neq \{1\}$ then we apply A1 to $P'' \leq P' \leq G$.
\[
tr^{P'}_{P''}(a_{P'}) = p\sum_{h \in P''/[P',P']} a_hh.
\]
Since the coefficient $a_{h_0} \neq 0$ we get that $tr^{P'}_{P''}(a_{P'}) = \pi^{P'}_{P''}(a_{P''}) \neq 0$. Hence $a_{P''} \neq 0$, a contradiction because $P'$ is the smallest order group such that $a_{P'} \neq 0$. Hence $\tilde{h}_0$ generates $P$. Since $P'$ is not cyclic $\tilde{h}_0 \neq 1$. Hence $[P',P'] < P$. Hence 
\[
\eta_P(tr^{P'}_P(a_{P'})) = p \Big( \sum_{h \text{ generates } P/[P',P']} a_hh \Big) \neq 0.
\]
But by A1 this must be same as $\pi^{P'}_{P}(a_P) = 0$, a contradiction. \\
Case 2. When $P$ is not cyclic, by A1 $tr^{P'}_P(a_{P'}) = 0$. But the coefficient of $h_0$ in $tr^{P'}_P(a_{P'})$ is $pa_{h_0} \neq 0$ which is a contradiction. \\
Hence $a_{P} = 0$ for all $P \leq G$ and we get injection of $\delta$.

Next we show that the image of $\delta$ restricted to $\psi^G_R$ lies in $R[Cong(G)]^{\tau}$. A3 says that $a_P \in T_{P,R}$ for every $P \in C(G)$. Let $a_P = tr(b_P)$ for some $b_P \in R[P]^{\tau}$. Then
\begin{align*}
\delta_P(a_P) & = \delta_P(\sum_{x \in W_GP} xb_Px^{-1}) \\
& = [N_GP:P] \delta_P(b_P).
\end{align*}
On the other hand
\begin{align*}
\sum_{x \in C(G, N_GP)} \delta_{xPx^{-1}}(a_{xPx^{-1}}) & = [G:N_GP] \delta_P(a_P), \qquad \text{by A2} \\
& = [G:P] \delta_P(b_P) \in R[Conj(G)]^{\tau}.
\end{align*}

\qed
 
\noindent{\bf Proof of theorem \ref{additivetheorem}:} $\delta|_{\psi^G_R}$ is injective and $\delta \circ \beta^G_R$ is identity on $R[Conj(G)]^{\tau}$. We claim that $\beta^G_R \circ \delta$ is identity on $\psi^G_R$. Let $(a_P) \in \psi^G_R$, then $\delta(\beta^G_R(\delta((a_P)))) = \delta((a_P))$. Since the image of $\beta^G_R$ is contained in $\psi^G_R$ and $\delta$ is injective on $\psi^G_R$, we get $\beta^G_R(\delta((a_P))) = (a_P)$. 

\qed

Consider the map 
\[
id_{\mathbb{Q}_p} \otimes \beta^G_R : \mathbb{Q}_p \otimes_{\mathbb{Z}_p} R[Conj(G)]^{\tau} \rightarrow \prod_{P \leq G} \mathbb{Q}_p \otimes R[P^{ab}]^{\tau}.
\]

\begin{proposition} The map $id_{\mathbb{Q}_p} \otimes \beta^G_R$ is injective. Its image consists of all tuples $(a_P)$ satisfying \\
A1. Let $P \leq P' \leq G$ such that $[P', P'] \leq P$ and if $P$ is a non-trivial cyclic group then $[P',P'] \neq P$. If $P$ is not cyclic, then $tr^{P'}_P(a_{P'}) = \pi^{P'}_P(a_P)$. If $P$ is not cyclic, then $tr^{P'}_P(a_{P'}) = \pi^{P'}_P(a_P)$. If $P$ cyclic and $P'$ is not cyclic then $\eta_P(tr^{P'}_P(a_{P'})) = \pi^{P'}_P(a_P)$. If $P'$ is cyclic, then $tr^{P'}_P(a_{P'})=0$. \\
A2. $(a_P)_{P\in C(G)}$ is invariant under conjugation action by every $g \in G$. \\
Hence, if $(id_{\mathbb{Q}_p} \otimes \beta^G_R) (a) = (a_P) \in  \prod_{P \leq G} \mathbb{Q}_p \otimes R[P^{ab}]^{\tau}$ is such that $a_P \in T_{P,R}$ for all $P \in C(G)$, then $a \in R[Conj(G)]^{\tau}$ and $a_P \in R[P^{ab}]^{\tau}$ for all $P \leq G$.
\label{additiveintegrality}
\end{proposition}

\noindent{\bf Proof:} The statement about injectivity and image follows directly from theorem \ref{additivetheorem}. Note that if $a_P$ satisfies A2 then it lies in $\mathbb{Q}_p \otimes T_{P,R}$. 

Now let $(id_{\mathbb{Q}_p} \otimes \beta^G_R) (a) = (a_P)$ be such that $a_P \in T_{P,R}$ for every $P \in C(G)$. The map $\delta$ shows that $a$ is determined by $a_P$'s for cyclic $P$. Because $a_P \in T_{P,R}$ for all $P \in C(G)$, the inverse image $a$ lies in $R[Conj(G)]^{\tau}$ and $a_P \in R[P^{ab}]^{\tau}$ for all $P \leq G$. 
\qed

\subsection{Logarithm for Iwasawa algebras} Logarithms on $K_1$ groups of $p$-adic orders were constructed by R. Oliver and M. Taylor. In this subsection we follow R. Oliver \cite{Oliver:1988} to define logarithm homomorphism on $K_1$-groups of Iwasawa algebras, a straightforward generalisation of the construction of R. Oliver and M. Taylor. In this subsection let $R$ denote the ring $\widehat{\Lambda_O(Z)_{(p)}}$. Let $J_R$ be the Jacobson radical of $R[G]^{\tau}$. Since $\mathcal{G}$ is pro-$p$, the ring $R[G]^{\tau}$ is a local ring and hence $J_R$ is its maximal ideal. We have the series
\[
Log(1+x) = \sum_{i=1}^{\infty} (-1)^{i+1}\frac{x^i}{i},
\]
and 
\[
Exp(x) = \sum_{i=0}^{\infty} \frac{x^i}{i!}.
\]

\begin{lemma} The ideal $J_R/pR[G]^{\tau}$ is a nilpotent ideal of $R[G]^{\tau}/pR[G]^{\tau}$.
\end{lemma}

\noindent{\bf Proof:} Let $k = O/(p)$. We have the following exact sequence
\[
0 \rightarrow J_R/pR[G]^{\tau} \rightarrow Q(k[[\mathcal{G}]]) \rightarrow Q(k[[\Gamma]]) \rightarrow 0.
\]
Let $N$ be the kernel of the map $k[[\mathcal{G}]] \rightarrow k[[\Gamma]]$, and let $I_H$ be the kernel of the map $k[H] \rightarrow k$. Then $N = k[[\mathcal{G}]]I_H$. Since $H$ is a finite $p$-group, we have
\[
N^n = k[[\mathcal{G}]] I_H^n = 0,
\]
for some positive integer $n$. By lemma \ref{lemma1} we can write any element $x \in Q(k[[\mathcal{G}]])$ as $x = \frac{a}{t}$, with $a \in k[[\mathcal{G}]]$ and $t \in k[[Z]]$. Also, $x \in J_R/pR[G]^{\tau}$ if and only if $a \in N$. As $t$ is central, we deduce that $J_R/pR[G]^{\tau}$ is nilpotent. 

\qed

\begin{lemma} Let $I \subset J_R$ be any ideal of $R[G]^{\tau}$. Then \\ 
1) For any $x\in I$, the series $Log(1+x)$ converges to an element in $R[G]^{\tau}[\frac{1}{p}]$. Moreover, for any $u, v \in 1+I$
\begin{equation}
\label{log1}
Log(uv) \equiv Log(u) + Log(v) (\text{mod } [R[G]^{\tau}[\frac{1}{p}], I[\frac{1}{p}]]).
\end{equation}
2) If $I \subset \xi R[G]^{\tau}$, for some central element $\xi$ such that $\xi^p \in p\xi R[G]^{\tau}$, then for all $u,v \in 1+I$, $Log(u)$ and $Log(v)$ converge to an element in $I$ and
\begin{equation}
\label{log2}
Log(uv) \equiv Log(u) + Log(v) (\text{mod } [R[G]^{\tau}, I]).
\end{equation}
In addition, if $I^p \subset pIJ_R$, then the series $Exp(x)$ converges to an element in $1+I$ for all $x \in I$; the maps defined by $Exp$ and $Log$ are inverse bijections between $I$ and $1+I$. Moreover, $Exp([R[G]^{\tau}, I]) \subset E(R[G]^{\tau},I)$ and for any $x,y \in I$, we have
\begin{equation}
\label{exp}
Exp(x+y) \equiv Exp(x)\cdot Exp(y) (\text{mod } E(R[G]^{\tau},I)).
\end{equation}
 
\end{lemma}

\noindent{\bf Proof:} This is an analogue of lemma 2.7 in R. Oliver \cite{Oliver:1988}. We reproduce the proof here indicating the modifications needed in our situation. \\

\noindent{\bf Step 1.} The ring $R[G]^{\tau}$ is $J_R$-adically complete (this follows from the previous lemma), for any $x \in I \subset J_R$, the terms $x^n$ converges to 0 in $I$ and $x^n/n$ converges to 0 in $I[\frac{1}{p}]$. Hence the series $Log(1+x)$ converges in $I[\frac{1}{p}] \subset R[G]^{\tau}[\frac{1}{p}]$.

If $I \subset \xi R[G]^{\tau}$, for some central $\xi$ such that $\xi^p \in p\xi R[G]^{\tau}$, then $I^p \subset pI$. Hence $I^n \subset nI$ for every positive integer $n$. So $x^n/n \in I$. Hence the series $Log(1+x)$ converges to an element in $I \subset R[G]^{\tau}$.

Furthermore, if $I^p \subset pIJ_R$, note that for any positive integer $n$ such that $p^k \leq  n < p^{k+1}$
\[
I^n \subset p^{([n/p]+[n/p^2]+\cdots+[n/p^k])}I J_R^k = n! IJ^k_R.
\]
Recall that $[y]$ denotes the greatest integer less than or equal to $y$ and that $n! p^{-([n/p]+\cdots+[n/p^k])}$ is a $p$-adic unit. Hence $Exp(x)$ converges to an element in $1+I$. The fact that $Log$ and $Exp$ are inverse bijections between $1+I$ and $I$ is formal. \\

\noindent{\bf Step 2.} For any $I \subset J_R$, set 
\[
U(I) = \sum_{m \geq 0, n \geq 1} \frac{1}{m+n} [I^m,I^n] \subset [R[G]^{\tau}[\frac{1}{p}], I[\frac{1}{p}]],
\]
a $R$-submodule of $R[G]^{\tau}[\frac{1}{p}]$. If $I \subset \xi R[G]^{\tau}$, where $\xi$ is a central element such that $\xi^p \in p\xi R[G]^{\tau}$, then $\xi^n \in n\xi R[G]^{\tau}$, and 
\[
U(I) = \langle [r, \frac{\xi^{m+n}}{m+n}s] : m \geq 0, n \geq 1, \xi^mr \in I^m, \xi^ns\in I^n, \xi r, \xi s \in I \rangle \subset [R[G]^{\tau}, I].
\]
So the congruences (\ref{log1}) and (\ref{log2}) will both follow once we have shown that for every $I \subset J_R$ and every $x,y \in I$
\[
Log((1+x)(1+y)) \equiv Log(1+x) + Log(1+y) (\text{mod } U(I)).
\]

For each $n \geq 1$, we let $W_n$ be the set of formal ordered monomials of length $n$ in two variables $a,b$. For $w \in W_n$, set 

$C(w)$ = orbit of $w$ in $W_n$ under cyclic permutations. 

$k(w)$ = number of occurrences of $ab$ in $w$.

$r(w)$ = coefficients of $w$ in $Log(1+a+b+ab) = \sum_{i=0}^{k(w)} (-1)^{n-i-1} \frac{1}{n-i}\binom{k(w)}{i}$.

If $w^{\prime} \in C(w)$, then it is clear that $w(x,y) \equiv w^{\prime}(x,y) (mod \ [I^i, I^j])$ for some $i,j \geq 1$ such that $i+j=n$. So
\begin{align*}
Log(1+x+y+xy) & = \sum_{n=1}^{\infty} \sum_{w\in W_n} r(w)w(x,y) \\
& \equiv \sum_{n=1}^{\infty} \sum_{w \in W_n/C} \Big( \sum_{w^{\prime} \in C(w)} r(w^{\prime})\Big)w(x,y) (mod \ U(I)).
\end{align*}
Let $ k =max \{k(w^{\prime}) : w^{\prime} \in C(w)\}$. Let $ |C(w)| = n/t$. Then $C(w)$ contains $k/t$ elements with exactly $(k-1)$ $ab$'s and $(n-k)/t$ elements with $k$ $ab$'s. Hence
\begin{align*}
\sum_{w^{\prime} \in C(w)} r(w^{\prime}) &= \frac{1}{t} \sum_{i=0}^{k} (-1)^{n-i-1} \frac{1}{n-i}\Big((n-k)\binom{k}{i} + k \binom{k-1}{i}\Big) \\
                                                                        &=\frac{1}{t} \sum_{i=0}^{k} (-1)^{n-i-1} \frac{1}{n-i}\Big((n-k)\binom{k}{i} + (k-i)\binom{k}{i}\Big) \\
                                                                        &= \frac{1}{t}\sum_{i=0}^{k} (-1)^{n-i-1} \binom{k}{i},
\end{align*}
which is 0 unless $k=0$, in which case it is equal to $(-1)^{n-1}\frac{1}{n}$. Thus
\[
Log(1+x+y+xy) \equiv \sum_{n=1}^{\infty} (-1)^{n-1} \Big(\frac{x^n}{n}+\frac{y^n}{n}\Big) = Log(1+x)+Log(1+y) \ (mod\ U(I)).
\]

\noindent{\bf Step 3.} We now prove the congruence (\ref{exp}). $Exp$ and $Log$ induce bijection between $I$ and $1+I$. Hence 
\[
Log(Exp(x)Exp(y)) \equiv x+y \ (mod \ U(I)),
\]
which gives 
\begin{align*}
Exp(x)Exp(y)Exp(x+y)^{-1} & \in Exp(x+y+U(I))Exp(-x-y) \\
& \subset Exp(U(I)) \subset Exp([R[G]^{\tau},I]).
\end{align*}

Hence we only need to prove that $Exp([R[G]^{\tau},I])$ is contained in $E(R[G]^{\tau},I)$. Choose and fix a set of $R$ generators $\{[s_1,v_1], \ldots, [s_m,v_m]\}$ of $[R[G]^{\tau}, I]$, with $s_i \in R[G]^{\tau}$ and $v_i \in I$ and choose an expression $x = \sum_{i=1}^{m}a_i[s_i,v_i]$ for an element $x \in [R[G]^{\tau},I]$. Define 
\[
\psi(x)= \prod_{i=1}^{m}(Exp(a_is_iv_i)Exp(a_iv_is_i)^{-1}).
\]
Note that $\psi(x)$ depends on the choice of the above expression for $x$. 

For any $r \in R[G]^{\tau}$ and any $x \in I$, an identity of Vaserstein gives
\begin{align*}
Exp(rx)Exp(xr)^{-1} & = \Big( 1+ r(\sum_{n=1}^{\infty} \frac{x(rx)^{n-1}}{n!})\Big)\Big(1+(\sum_{n=1}^{\infty}\frac{x(rx)^{n-1}}{n!})r\Big)^{-1} \\
& \in E(R[G]^{\tau}, I).
\end{align*}
Hence $Im(\psi) \subset E(R[G]^{\tau}, I)$. For any $ k \geq 1$ and any $x,y \in p^kI$, 
\[
Exp(x)Exp(y) \equiv Exp(x+y) (\text{mod } U(p^kI) \subset p^{2k}U(I) \subset p^{2k}[R[G]^{\tau},I]).
\]
Also, for any $k,l \geq 1$ and any $ x \in p^kI$, $y \in p^lI$, 
\[
Exp(x)Exp(y) \equiv Exp(y)Exp(x) (\text{mod } [p^kI, p^lI] \subset p^{k+l}[R[G]^{\tau},I]).
\]
So for any $ l \geq k \geq 1$, and any $x \in p^k[R[G]^{\tau}, I]$ and $y\in p^l[R[G]^{\tau}, I]$, choose an expression $x = \sum_{i=1}^{m}a_i[s_i,v_i]$ and $y = \sum_{i=1}^{m}b_i[s_i,v_i]$ and $x + y = \sum_{i=1}^{m}(a_i+b_i)[s_i,v_i]$, then
\[
\psi(x) \equiv Exp(x) (\text{mod } p^{2k}[R[G]^{\tau}, I])
\]
\begin{equation}
\label{psiequation}
\psi(x+y) \equiv \psi(x)\psi(y) \equiv \psi(x)Exp(y) (\text{mod } p^{k+l}[R[G]^{\tau}, I]).
\end{equation}
For arbitrary $u \in Exp(p[R[G]^{\tau}, I])$, define a sequence $x_0, x_1, x_2, \ldots$ in $[R[G]^{\tau}, I]$ by setting
\[
x_0= Log(u) \in p[R[G]^{\tau}, I]; \qquad x_{i+1} = x_i + Log(\psi(x_i)^{-1}u),
\]
We need to choose a representation $x_0 = \sum_{j=1}^{m}a_j[s_j,v_j]$ and at each stage we need to choose a representation $Log(\psi(x_i)^{-1}u) = \sum_{j=1}^m c_j [s_j,v_j]$. This inductively gives a natural representation for all $x_i$'s. By (\ref{psiequation}), applied inductively for all $i \geq 0$,
\[
\psi(x_i) \equiv u, \quad x_{i+1} \equiv x_i \quad (\text{mod } p^{2+i}[R[G]^{\tau},I]).
\]
Then ${x_i}$ converges and we get a natural representation for $\lim_{i \rightarrow \infty} x_i$ in terms as a linear combination of the fixed set of generators. If we use this to define $\psi(\lim_{i \rightarrow \infty} x_i)$, we get $u = \psi(\lim_{i \rightarrow \infty}   x_i)$. This shows that 
\[
Exp(p[R[G]^{\tau}, I]) \subset Im(\psi) \subset E(R[G]^{\tau}, I).
\]

Now define subgroups $D_k$, for all $ k \geq 0$, by setting
\[
D_k = \langle rx-xr : x \in I, r \in R[G]^{\tau}, rx, xr \in IJ_R^k \rangle \subset [R[G]^{\tau}, I]\cap IJ_R^k.
\]
By the hypothesis on $I$, for all $k \geq 0$, 
\begin{align*}
U(IJ^k_R) &= \sum_{m,n \geq 1} \frac{1}{m+n} [(IJ_R^k)^m, (IJ_R^k)^n] \\
&= \langle [r, \frac{\xi^n}{n}s] : n \geq 2, \xi r, \xi s \in IJ_R^k , \xi^nrs, \xi^nsr \in (IJ_R^k)^n \subset nIJ_R^{k+1} \rangle \\
& \subset D_{k+1}.
\end{align*} 
This shows that $Exp(D_k) \subset Exp([R[G]^{\tau}, I])$ are both normal subgroups of $(R[G]^{\tau})^{\times}$. Also, for any $x, y \in IJ_R^k$,
\begin{equation}
\label{uik}
Exp(x)Exp(y) \equiv Exp(x+y) (\text{mod } Exp(U(IJ_R^k)) \subset Exp(D_{k+1})).
\end{equation}
For any $k \geq 0$ and any $x \in D_k$, if we write $x = \sum (r_ix_i - x_ir_i)$, where $r_i \in R[G]^{\tau}, x_i \in I$ and $r_ix_i, x_ir_i \in IJ_R^k$, then by above 
\begin{align*}
Exp(x) & \equiv \prod (Exp(r_ix_i)Exp(x_ir_i)^{-1}) (\text{mod } Exp(D_{k+1})) \\
& \equiv 1 (\text{mod } E(R[G]^{\tau}, I)).
\end{align*}
In other words, $Exp(D_k) \subset E(R[G]^{\tau}, I)Exp(D_{k+1})$ for all $k \geq 0$. But for a large enough $k$, $D_k \subset p[R[G]^{\tau}, I]$. Hence we get 
\[
Exp(R[G]^{\tau}, I]) \subset Exp(D_0) \subset E(R[G]^{\tau}, I)Exp(p[R[G]^{\tau}, I]) \subset E(R[G]^{\tau}, I).
\]
 \qed

\begin{proposition} For any ideal $I \subset J_R$ of $R[G]^{\tau}$, the $p$-adic logarithm $Log(1+x)$, for any $x\in I$, induces a unique homomorphism 
\[
log_I : K_1(R[G]^{\tau}, I) \rightarrow (I/[R[G]^{\tau},I])\otimes_{\mathbb{Z}_p}\mathbb{Q}_p.
\]
If, furthermore, $I \subset \xi R[G]^{\tau}$ for some central $\xi$ such that $\xi^p \in p\xi R[G]^{\tau}$, then the logarithm induces a homomorphism 
\[
log_I : K_1(R[G]^{\tau}, I) \rightarrow I/[R[G]^{\tau},I],
\]
and $log_I$ is an isomorphism if $I^p \subset pIJ_R$.
\label{logdefn}
\end{proposition}  

\noindent{\bf Proof:} This is an analogue of theorem 2.8 of R. Oliver \cite{Oliver:1988}. By the previous lemma 
\begin{equation}
\label{log}
L : 1+I \xrightarrow{Log} I[\frac{1}{p}] \xrightarrow{proj} I[\frac{1}{p}]/[R[G]^{\tau}, I[\frac{1}{p}]],
\end{equation}
is a homomorphism.  
For each $ n \geq 1$, let
\begin{equation}
\label{trlog}
Tr_n: M_n(I[\frac{1}{p}])/[M_n(R[G]^{\tau}[\frac{1}{p}]), M_n(I[\frac{1}{p}])] \rightarrow I[\frac{1}{p}]/[R[G]^{\tau}[\frac{1}{p}], I[\frac{1}{p}]]
\end{equation}
be the homomorphism induced by the trace map. Then (\ref{log}), applied to the ideal $M_n(I) \subset M_n(R[G]^{\tau})$, induces a homomorphism 
\begin{align*}
L_n: 1+M_n(I) = GL_n(R[G]^{\tau},I) & \xrightarrow{Log} M_n(I[\frac{1}{p}])/[M_n(R[G]^{\tau}[\frac{1}{p}]), M_n(I[\frac{1}{p}])] \\ & \xrightarrow{Tr_n} I[\frac{1}{p}]/[R[G]^{\tau}, I[\frac{1}{p}]].
\end{align*}
For any $n$, and any $u \in 1+M_n(I)$ and $r \in GL_n(R[G]^{\tau})$. 
\[
L_n([r,u]) = L_n(rur^{-1}) - L_n(u) = Tr_n(rLog(u)r^{-1}) - Tr_n(Log(u)) = 0.
\]
So $L_{\infty} = \cup L_n$ factors through a homomorphism 
\begin{align*}
log_I : K_1(R[G]^{\tau}, I) = GL(R[G]^{\tau}, I)/[ & GL(R[G]^{\tau}), GL(R[G]^{\tau}, I)] \\&\rightarrow I[\frac{1}{p}]/[R[G]^{\tau}[\frac{1}{p}], I[\frac{1}{p}]].
\end{align*}

If $I \subset \xi R[G]^{\tau}$, for some central element $\xi$, such that $\xi^p \in p\xi R[G]^{\tau}$, then the same argument by the second part of the above lemma gives a homomorphism 
\[
log_I : K_1(R[G]^{\tau}, I) \rightarrow I/[R[G]^{\tau}, I].
\]
If in addition, $I^p \subset pIJ_R$, then $Log$ is bijective and $Log^{-1}([R[G]^{\tau}, I]) \subset E(R[G]^{\tau}, I])$, by the last part of lemma above. Hence $log_I$ is an isomorphism. 

\qed

Now consider the case of $\Lambda_O(\mathcal{G})$. Let $F$ be the field of fractions of $O$. For a profinite group $P$ let $F[[P]]$ denote $\ilim{U} F[P/U]$, where $U$ runs through open normal subgroups of $P$. Let $G_n = \mathcal{G}/Z^{p^n}$. Then for each $n \geq 0$ we have, by R. Oliver (\cite{Oliver:1988} chapter 2), a unique homomorphism 
\[
log: K_1(O[G_n]) \rightarrow F[Cong(G_n)].
\]
Since logarithm kills torsion we have a homomorphism
\[
log: K_1'(O[G_n]) \rightarrow F[Cong(G_n)].
\]
Taking the inverse limit over $n$'s gives 
\[
log : K_1'(\Lambda_O(\mathcal{G})) \rightarrow \ilim{n} F[Cong(G_n)] = F[[\mathcal{G}]].
\]

\begin{remark} In the case when $R = \widehat{\Lambda_{O}(Z)_{(p)}}$ is it is not clear if we can extend the homomorphism $log_{J_R}$ to $K_1(R[G]^{\tau})$ is any canonical fashion. However, for our purposes $log_{J_R}$ suffices.
\end{remark}

\subsection{Integral logarithm} Again we follow R. Oliver \cite{Oliver:1988} to construct the integreal logarithm homomorphism from $K_1$ of Iwasawa algebras. This is  a straightforward generalisation of the integral logarithm on $K_1$ of $p$-adic group rings of finite groups contructed by R. Oliver and M. Taylor. 

\noindent As before $R$ denotes either $\Lambda_O(Z)$ or $\widehat{\Lambda_O(Z)_{(p)}}$. 

\begin{definition} Let $\varphi$ be the map on $R$ induced by the Frobenius map on $O$ and the $p$-power map on $Z$. We extend this to a map, still denoted by $\varphi$, to $R[Conj(G)]^{\tau}$ by mapping $[\overline{g}]$ to $[\overline{g}^p]$ (remember that $\overline{g}^p$ is not the same as $\overline{g^p}$). 
\label{integrallogforiwasawaalgebra}
\end{definition}

\subsubsection{Integral logarithm homomorphism for $\Lambda(\mathcal{G})$} 

For every $n \geq 0$, let $G_n = \mathcal{G}/Z^{p^n}$. The integral logarithm map defined by R. Oliver and M. Taylor is
\[
L: K_1'(O[G_n]) \rightarrow O[Conj(G_n)],
\]
defined as $L = log - \frac{\varphi}{p}log$, where $\varphi: O[Conj(G_n)] \rightarrow O[Conj(G_n)]$ is the map induced by Frobenius on $O$ and the $p$-power map on $G_n$.  The kernel of $L$ is $K_1'(O[G_n])_{tor}$ which is equal to $\mu(O) \times G_n^{ab}$ by a theorem of G. Higman \cite{Higman:1940} and C.T.C. Wall \cite{Wall:1974}. Here $\mu(O)$ denotes the torsion subgroup of $O^{\times}$. In fact, we have an exact sequence (R. Oliver \cite{Oliver:1988}, theorem 6.6)
\begin{equation}
\label{fp}
1 \rightarrow \mu(O) \times G_n^{ab} \rightarrow K_1'(O[G_n]) \xrightarrow{L} O[Cong(G_n)] \xrightarrow{\omega}  G_n^{ab} \rightarrow 1.
\end{equation}
Here the map $\omega$ is defined by
\[
\sum_{g \in Conj(G_n)} a_gg = \prod_{g} (g)^{tr_{O/\mathbb{Z}_p}(a_g)}.
\]

\begin{definition}  We define the integral logarithm map $L$ on $K_1'(\Lambda_O(\mathcal{G}))$ by taking the inverse limit of the sequence (\ref{fp}) over all $n$'s and get an exact sequence
\[
1 \rightarrow \mu(O) \times \mathcal{G}^{ab} \rightarrow K_1'(\Lambda_O(\mathcal{G})) \xrightarrow{L} \Lambda_O(Z)[Cong(G)]^{\tau} \xrightarrow{\omega}  \mathcal{G}^{ab} \rightarrow 1.
\]
\label{integrallogforlambda}
 \end{definition}

\subsubsection{Integral logarithm homomorphism for $\widehat{\Lambda(\mathcal{G})_S}$} Let $J$ denote the kernel of the natural surjection
\[
\widehat{\Lambda_O(\mathcal{G})_S} \rightarrow \widehat{\Lambda_O(\Gamma)_{(p)}}.
\]
Since $\widehat{\Lambda_O(\mathcal{G})_S}$ is local, the following maps are surjective
\[
\widehat{\Lambda_O(\mathcal{G})_S}^{\times} \rightarrow K_1'(\widehat{\Lambda_O(\mathcal{G})_S}),
\]
and 
\[
1+J \rightarrow K_1(\widehat{\Lambda_O(\mathcal{G})_S}, J).
\]
We have an exact sequence 
\[
1 \rightarrow 1+J \rightarrow \widehat{\Lambda_O(\mathcal{G})_S}^{\times} \rightarrow \widehat{\Lambda_O(\Gamma)_{(p)}}^{\times} \rightarrow 1,
\]
which splits using the distinguished embedding $\Gamma \hookrightarrow \mathcal{G}$. Hence any element $x \in \widehat{\Lambda_O(\mathcal{G})_S}^{\times}$ can be written uniquely as $x = uy$, for $u \in 1+J$ and $y \in \widehat{\Lambda_O(\Gamma)_{(p)}}^{\times}$. As a result we get the following
 
 \begin{lemma} Every $x \in K_1'(\widehat{\Lambda_O(\mathcal{G})_S})$ can be written as a product $x = uy$, where $ y \in K_1'(\widehat{\Lambda_O(\Gamma)_{(p)}})$ and $u$ lies in the image of $K_1(\widehat{\Lambda_O(\mathcal{G})_S}, J)$ in $K_1'(\widehat{\Lambda_O(\mathcal{G})_S})$.
\end{lemma}

\qed

\begin{lemma} For any $y \in K_1(\widehat{\Lambda_O(\Gamma)_{(p)}})$, we have
\[
\frac{y^p}{\varphi(y)} \equiv 1 (\text{mod } p\widehat{\Lambda_O(\Gamma)_{(p)}})
\]
Hence the logarithm of $\frac{y^p}{\varphi(y)}$ is defined.
\end{lemma}
\noindent{\bf Proof:} The ring $\widehat{\Lambda_O(\Gamma)_{(p)}}/p\widehat{\Lambda_O(\Gamma)_{(p)}}$ is isomorphic to the domain $O/(p)[[\Gamma]]$ (Recall that $(p)$ is the maximal ideal of $O$ since we have assumed that $O$ is unramified extension of $\mathbb{Z}_p$). Let $O/(p) = \mathbb{F}_q$. We also note that $\mathbb{F}_q[[\Gamma]] \cong \mathbb{F}_q[[X]]$. On this ring the $\varphi$ map is the one induced by $X \mapsto (1+X)^p-1 = X^p$ and by $p$-power Frobenius on $\mathbb{F}_q$. Let $ \overline{y} \in \mathbb{F}_q[[X]]$. Write $\overline y = \sum_{i=0}^{\infty} a_iX^i$. Then
\[
\overline{y}^p = (\sum_{i} a_iX^i)^p = \sum_i a_i^p X^{ip} = \varphi(\overline{y}),
\]
Hence the lemma.
\qed
 
\begin{definition} We define the integral logarithm $L$ on $K_1'(\widehat{\Lambda_O(\mathcal{G})_S})$ as follows: write any $x \in K_1'(\widehat{\Lambda_O(\mathcal{G})_S})$ as $x = uy$, with $y \in K_1'(\widehat{\Lambda_O(\Gamma)_{(p)}})$ and $u$ lies in the image of $K_1(\widehat{\Lambda_O(\mathcal{G})_S}, J)$ in $K_1'(\widehat{\Lambda_O(\mathcal{G})_S})$. Define
\[
L(x) = L(uy) = L(u)+L(y) = log(u)- \frac{\varphi}{p}log(u) + \frac{1}{p}log(\frac{y^p}{\varphi(y)}).
\]
\end{definition}

\begin{proposition} $L$ induces a homomorphism 
\[
L: K_1'(\widehat{\Lambda_O(\mathcal{G})_S}) \rightarrow \widehat{\Lambda_O(Z)_{(p)}}[Conj(G)]^{\tau}
\]
and $L$ is independent of the choice of splitting of $\mathcal{G} \rightarrow \Gamma$.
\end{proposition}
\noindent{\bf Proof:} First we show that image of $L$ lies in $\widehat{\Lambda_O(Z)_{(p)}}[Conj(G)]^{\tau}$. Let $x \in K_1'(\widehat{\Lambda_O(\mathcal{G})_S})$. Write it as $x=uy$, with  $y \in K_1'(\widehat{\Lambda_O(\Gamma)_{(p)}})$ and $u$ in the image of $K_1(\widehat{\Lambda_O(\mathcal{G})_S}, J)$ in $K_1'(\widehat{\Lambda_O(\mathcal{G})_S})$. We will show that $L(u)$ and $L(y)$ lie in $\widehat{\Lambda_O(Z)_{(p)}}[Conj(G)]^{\tau}$. The fact that $L(y) = \frac{1}{p}log(\frac{y^p}{\varphi(y)})$ lies in $\widehat{\Lambda_O(Z)_{(p)}}$ follows from proposition \ref{logdefn}. Now let $u=1-v$ and consider $L(1-v)$
\begin{align*}
L(1-v) &= -(v+\frac{v^2}{2} +\frac{v^3}{3}+ \cdots)+(\frac{\varphi(v)}{p}+\frac{\varphi(v)}{2p}+\cdots) \\
& \equiv \sum_{k=1}^{\infty} \frac{1}{pk}(v^{pk}-\varphi(v^k)) \quad (\text{mod } \widehat{\Lambda_O(Z)_{(p)}}[Cong(G)]^{\tau})
\end{align*}
It suffices that $pk| (v^{pk} - \varphi(v^k))$ for all $k$; which will follow from 
\[
p^n|(v^{p^n}-\varphi(v^{p^{n-1}})),
\]
for all $n \geq 1$. Write $v = \sum_{g \in G} r_g\bar{g}$, where $r_g \in \widehat{\Lambda_O(Z)_{(p)}}$. Set $ q =p^n$ and consider a typical term in $v^q$:
\[
r_{g_1}\cdots r_{g_q} \bar{g}_1\cdots \bar{g}_q.
\]
Let $\mathbb{Z}/q\mathbb{Z}$ acts cyclically permuting the $g_i$'s, so that we get a total of $p^{n-t}$ conjugate terms, where $p^t$ is the number of cyclic permutations leaving each term invariant. Then $\bar{g}_1\cdots \bar{g}_q$ is a $p^t$-th power, and the sum of the conjugate terms has the form
\[
p^{n-t} \hat{r}^{p^t} \hat{g}^{p^t} \in \widehat{\Lambda_O(Z)_{(p)}}[Conj(G)]^{\tau},
\]
where $\hat{r} = \prod_{j=1}^{p^{n-t}} r_{g_j}$ and $\hat{g} = \prod_{j=1}^{p^{n-t}} g_j$. Here we may do the multiplication after rearrangement because the twisting map $\tau$ is symmetric. If $t=0$, then this is a multiple of $p^n$. If $t >0$, then there is a corresponding term $p^{n-t} \hat{r}^{p^{t-1}} \hat{g}^{p^{t-1}}$ in the expansion of $v^{p^{n-1}}$. It remains only to show that
\[
p^{n-t}\hat{r}^{p^t}\hat{g}^{p^t} \equiv p^{n-t} \varphi(\hat{r}^{p^{t-1}}\hat{g}^{p^{t-1}}) = p^{n-t} \varphi(\hat{r}^{p^{t-1}}) \hat{g}^{p^t} (\text{mod } p^n).
\]
But $p^t |(\hat{r}^{p^t} - \varphi(\hat{r}^{p^{t-1}}))$, since $p|(\hat{r}^p - \varphi(\hat{r}))$. 

Next we show that $L$ is independent of the choice of splitting $\mathcal{G} \rightarrow \Gamma$. Since $\mathcal{G}$ is one dimensional, there are only finitely many splittings $\Gamma \rightarrow \mathcal{G}$. Let $A$ be the intersection of image of $\Gamma$ under all the splittings. Then $A$ is an open subgroup of $\mathcal{G}$. It is clear that for any $x \in Im(K_1(\widehat{\Lambda_O(A)_S}) \rightarrow K_1'(\widehat{\Lambda_O(\mathcal{G})_S}))$, $L(x)$ is independent of the choice of splitting $\Gamma \rightarrow \mathcal{G}$.  

Let $x \in Im(K_1(\widehat{\Lambda_O(\mathcal{G})_S},J_R) \rightarrow K_1'(\widehat{\Lambda_O(\mathcal{G})_S}))$. Recall that $J_R$ is the Jacobson radical of $\widehat{\Lambda_O(\mathcal{G})_S}$. Let $x=u_1y_1=u_2y_2$ be two different decompositions of $x$ given by two different splittings of $\mathcal{G} \rightarrow \Gamma$. Then $log(y_i)$ is defined as it is defined for $x$ and $u_i$ (for $i=1,2$) and $log(u_1) + log(y_1) = log(u_2)+log(y_2)$. Hence $L(x)$ is independent of the choice of splitting $\Gamma \rightarrow \mathcal{G}$.

Hence $L$ is well defined on $Im(K_1(\widehat{\Lambda_O(\mathcal{G})_S},J_R)) Im(K_1(\widehat{\Lambda_O(A)_S})) \subset K_1'(\widehat{\Lambda_O(\mathcal{G})_S})$. Now note that $K_1'(\widehat{\Lambda_O(\mathcal{G})_S})/Im(K_1(\widehat{\Lambda_O(\mathcal{G})_S},J_R)) Im(K_1(\widehat{\Lambda_O(A)_S}))$ is torsion group since $\widehat{\Lambda_O(\mathcal{G})_S}^{\times}/(1+J_R)\widehat{\Lambda_O(A)_S}^{\times}$ is a torsion group. Hence $L$ is independent of the choice of splitting of $\mathcal{G} \rightarrow \Gamma$.
\qed

\begin{remark} It would be interesting to find the kernel and cokernel of $L$ in this case. Since we do not need it we will not investigate it here (but see Ritter-Weiss \cite{RitterWeiss:2010b}).
\end{remark}

\subsection{Relation between the maps $\theta$ and $\beta$} 

For any subgroup $P$ of $G$, the map $t^G_P$ naturally extends to a map 
\[
t^G_P : F[[Z]][Cong{G}]^{\tau} \rightarrow F[[Z]][Cong(P)]^{\tau}.
\]

\begin{lemma} For any $P \leq G$, we have the following diagram 
\[
\xymatrix{ K_1'(\Lambda_O(\mathcal{G})) \ar[rr]^{log} \ar[d]_{\theta^G_P} & & F[[Z]][Conj(G)]^{\tau} \ar[d]^{t^G_P} \\
K_1'(\Lambda_O(U_P)) \ar[rr]_{log} & & F[[Z]][Conj(P)]^{\tau}}
\]
In the case of $p$-adic completions the following diagram commutes
\[
\xymatrix{ K_1(\widehat{\Lambda_O(\mathcal{G})_S}, J) \ar[rr]^{log} \ar[d]_{\widehat{\theta^G_P}} & & \widehat{\Lambda_O(Z)_{(p)}}[Conj(G)]^{\tau}[\frac{1}{p}] \ar[d]^{t^G_P} \\
K_1(\widehat{\Lambda_O(U_P)_S},J) \ar[rr]_{log} & & \widehat{\Lambda_O(Z)_{(p)}}[Conj(P)]^{\tau}[\frac{1}{p}]}
\]
\label{relationthetat}
\end{lemma}
\noindent{\bf Proof:} The first assertion is theorem 6.8 in R. Oliver \cite{Oliver:1988} which we refer for details. The second case is similar : for any $u \in K_1(\widehat{\Lambda_O(\mathcal{G})_S}, J)$, the result follows from the expressions
\[
log(u) = \lim_{n \rightarrow \infty} \frac{1}{p^n}(u^{p^n}-1); 
\]
and  
\[
 \widehat{\theta^G_P}(u) = \lim_{n \rightarrow \infty}(1+t^G_P(u^{p^n}-1))^{1/p^n}.
\]
\qed

\begin{lemma} We have the following commutative diagram for all $P \in C(G)$ and $P \neq \{1\}$
\[
\xymatrix{ \Lambda_O(U_P)^{\times} \ar[r]^{log} \ar[d]_{\alpha_P} &  F[[U_P]] \ar[d]^{p\eta_P} \\ 
 \Lambda_O(U_P)^{\times} \ar[r]_{log} &  F[[U_P]]}.
\]
And in the case of $p$-adic completion the following diagram commutes
\[
\xymatrix{  K_1(\widehat{\Lambda_O(U_P)_S}, J) \ar[rr]^{log} \ar[d]_{\alpha_P} & &  \widehat{\Lambda_O(U_P)_S}[\frac{1}{p}] \ar[d]^{p\eta_P}\\
K_1(\widehat{\Lambda_O(U_P)_S},J) \ar[rr]_{log} & &  \widehat{\Lambda_O(U_P)_S}[\frac{1}{p}] }
\]
\label{relationalphaeta}
\end{lemma}

\noindent{\bf Proof:} Note that the following diagram commutes
\[
\xymatrix{ \Lambda_{O[\mu_p]}(U_P)^{\times} \ar[r]^{log} \ar[d]_{\omega_P} & F(\mu_p)[[U_P]] \ar[d]^{\omega_P} \\
\Lambda_{O[\mu_p]}(U_P)^{\times} \ar[r]_{log} & F(\mu_p)[[U_P]] }
\]
and
\[
\eta_P = \frac{1}{p}(p - \sum_{k=0}^{p-1} \omega_P^k).
\]
Hence we get the first claim of the lemma. The second one is similar. 
\qed

\begin{definition} We define the map $v^G_P : \prod_{C \leq G} F[[Z]][C^{ab}] \rightarrow  F[[Z]][P^{ab}]$ as follows: \\
If $P$ is not cyclic 
\[
v^G_P((x_C)) =\Big(\sum_{P'} \frac{[P':P'^p]}{[P:(P')^p]} \varphi(x_{P'})\Big) 
\]
where sum ranges over all $P' \in C(G)$, such that $P'^p \leq P$. \\
If $P$ is cyclic
\[
v_P^G((x_C)) = \sum_{P'}[P':P'^p] \varphi(x_{P'}) = p\sum_{P'} \varphi(x_{P'}),
\]
where the sum ranges over all $P' \in C(G)$ such that $P' \neq P'^p = P$. The empty sum is taken to be 0. \\
Put $v^G = (v^G_P)_P$. We denote the analogous map in the case of $p$-adic completions again by $v^G$.
\end{definition}

\begin{lemma} For $P \neq \{1\}$, the following diagram commutes 
\[
\xymatrix{ F[[Z]][Conj(G)]^{\tau} \ar[rr]^{\varphi} \ar[d]_{\beta^G} & & F[[Z]][Conj(G)]^{\tau} \ar[d]^{\beta^G_P} \\
\prod_{C\leq G} F[[Z]][C^{ab}] \ar[rr]_{v^G_P} & &  F[[Z]][P^{ab}]}
\]
In the case $P = \{1\}$, the following diagram commutes
\[
\xymatrix{ F[[Z]][Conj(G)]^{\tau} \ar[rr]^{\varphi} \ar[d]_{\beta^G} & & F[[Z]][Conj(G)]^{\tau} \ar[d]^{\beta^G_P} \\
\prod_{C\leq G} F[[Z]][C^{ab}] \ar[rr]_{\varphi + v^G_{1}} & &  F[[Z]]}
\]
Here $\varphi$ in the lower row is just the map $\varphi: F[[Z]] \rightarrow F[[Z]]$. Analogous result holds for the $p$-adic completions. 
\label{relationbetaphi}
\end{lemma} 

\noindent{\bf Proof:} We assume that $P \neq \{1\}$ and $P$ is cyclic. For $g \in G$, we must show that $\beta^G_P(\varphi([\overline{g}])) = v^G_P(\beta^G([\overline{g}]))$. We start with the left hand side 
\[
\beta^G_P(\varphi([\overline{g}])) = \eta_P(t^G_P([\overline{g}^p])),
\]
which is non-zero if and only if some conjugate of $g^p$ generates $P$. Hence we may assume that $g^p$ itself is a generator of $P$. Then
\begin{align*}
v^G_P(\beta^G([\overline{g}])) & = (\sum_{P' \text{ s.t. } P'^p=P} [P':P'^p] \varphi(\beta^G_{P'}([\overline{g}]))) \\
& = (\sum_{P' \text{ s.t. } P'^p=P} [P':P'^p] \varphi(t^G_{P'}([\overline{g}]))) \\
& = (\sum_{P' \text{ s.t. } P'^p=P} [P':P'^p] \varphi(\sum_{x \in C(G,P')} \{ [\overline{x^{-1}gx}] | x^{-1}gx \in P'\})) \\
&= (\sum_{P' \text{ s.t. } P'^p=P} [P':P'^p] \sum_{x \in C(G,P')} \{[\overline{x^{-1}g^px}] | x^{-1}gx \in P'\}) \\
&= \sum_{P' \text{ s.t. } P'^p=P}  (\sum_{x \in C(G,P)} \{[\overline{x^{-1}g^px}] | x^{-1}gx \in P'\}) \\
& = \sum_{x \in C(G,P)} \sum_{P' \text{ s.t. } P'^p=P} \{[\overline{x^{-1}g^px}] | x^{-1}gx \in P'\} \\
& = \sum_{x \in C(G,P)} \{[\overline{x^{-1}g^px}] | x^{-1}g^px \in P\} \\
& = \beta^G_P(\varphi([\overline{g}])).
\end{align*}

Next we assume that $P$ is not cyclic. Then $\beta^G_P(\varphi([\overline{g}]))$ is non-zero if and only if some conjugate of $g^p$ lies in $P$. Hence we assume that $g^p \in P$ and $g \neq 1$. Let $C$ be the set of all conjugates of the subgroup $\langle g^p \rangle$ in $P$.
\begin{align*} 
v^G_P(\beta^G([\overline{g}])) &= \Big(\sum_{P' \text{ s.t. } P'^p <P} \frac{[P':P'^p]}{[P:P'^p]} \varphi (\beta^{G}_{P'}([\overline{g}]))\Big) \\
& = \Big(\sum_{P' \text{ s.t. } P'^p \in C} \frac{[P':P'^p]}{[P:P'^p]} \varphi(t^G_{P'}([\overline{g}]))\Big) \\
&= \sum_{P' \text{ s.t. } P'^p \in C} \Big( \frac{1}{[P:P'^p]}\sum_{x \in C(G, P'^p)} \{ [\overline{x^{-1}g^px}] : x^{-1}gx \in P'\}\Big) \\
&= \sum_{x \in C(G,P)} \{ [\overline{x^{-1}g^px}] : x^{-1}g^px \in P\} \\
&= \beta^G_P(\varphi([\overline{g}]))
\end{align*} 

It remains to show the equality when $g=1$. In this case
\begin{align*}
v^G_P(\beta^G(1)) & = \frac{1}{|P|}|G| \\
& = [G:P] \\
& = \beta^G_P(\varphi(1)).
\end{align*}

Now consider the remaining case $P = \{1\}$. We assume that $g$ is of order $p$ which is the only possibly nontrivial case. Let $C$ be the set of all conjugates in $G$ of the subgroup $\langle g \rangle$.
\begin{align*} 
v^G_{\{1\}}(\beta^G([\overline{g}])) & =  \sum_{P' \in C} \Big([P':P'^p] \varphi(\beta^G_{P'}([\overline{g}])) \Big) \\
& =   \sum_{P' \in C} \Big( [P':P'^p] \sum_{x \in C(G, P')} \{ [\overline{x^{-1}g^px}] : x^{-1}gx \in P'\} \Big) \\
& = \sum_{P' \in C} \Big( [P':P'^p] |N_GP'| \Big) \\
& = p|G| \\
& = \beta^G_{\{1\}}(\varphi([\overline{g}])) = \beta^G_{\{1\}}(\varphi([\overline{g}])) - \varphi(\beta^G_{\{1\}}([\overline{g}])).
\end{align*}
The last equality holds because clearly $\beta^G_{\{1\}}([\overline{g}]) = 0$. 
\qed

\begin{definition} We define the map $u^G_P: \prod_{C \leq G} \Lambda_O(U_C^{ab})^{\times} \rightarrow \Lambda_O(U_P^{ab})^{\times}$ as follows: \\
If $P$ is not a cyclic subgroup of $G$, then define
\[
u^G_P((x_C)) = \prod_{P'} \varphi(x_{P'})^{|P'|},
\]
where the product ranges over all $P' \in C(G)$ such that $P'^p \leq P$. \\
If $P$ is cyclic then define 
\[
u^G_P((x_C)) = \prod_{P'} \varphi(x_{P'}),
\]
where the product is taken over all $P' \in C(G)$ such that $P' \neq P'^p = P$. The empty product is taken to be 1. \\
Put $u^G = (u^G_P)_P$. We denote the analogous map in the case of $p$-adic completions again by $u^G$.
\end{definition}

\begin{lemma} For a noncyclic subgroup $P$ of $G$ the following diagram commutes
\[
\xymatrix{ \prod_{C \leq G} \Lambda_O(U_C^{ab})^{\times} \ar[r]^{log} \ar[d]_{u^G_P} & \prod_{C \leq G} F[[Z]][C^{ab}] \ar[d]^{|P|v^G_P} \\
 \Lambda_O(U_P^{ab})^{\times} \ar[r]_{log} &  F[[Z]][P^{ab}]}
\]
For $P \in C(G)$ the following diagram commutes
\[
\xymatrix{ \prod_{C \leq G} \Lambda_O(U_C^{ab})^{\times} \ar[r]^{log} \ar[d]_{u^G_P} & \prod_{C \leq G} F[[Z]][C^{ab}] \ar[d]^{\frac{1}{p}v^G_P} \\
\Lambda_O(U_P^{ab})^{\times} \ar[r]_{log} & F[[Z]][P^{ab}]}
\]
Similarly for the $p$-adic completions we have: Let $R = \widehat{\Lambda_O(Z)}_{(p)}$. Let $J_P$ denote the kernel of the natural map $R[P^{ab}] \rightarrow \widehat{\Lambda_O(\Gamma)}_{(p)}$. Then
for a noncyclic subgroup $P$ of $G$ the following diagram commutes
\[
\xymatrix{ \prod_{C \leq G} 1+J_C \ar[r]^{log} \ar[d]_{u^G_P} & \prod_{C \leq G} \mathbb{Q}_p \otimes J_C \ar[d]^{|P|v^G_P} \\
 1+J_P \ar[r]_{log} &  \mathbb{Q}_p \otimes J_P}
\]
For $P \in C(G)$ the following diagram commutes
\[
\xymatrix{ \prod_{C \leq G} 1+J_C  \ar[r]^{log} \ar[d]_{u^G_P} & \prod_{C \leq G} \mathbb{Q}_p \otimes J_C \ar[d]^{\frac{1}{p}v^G_P} \\
 1+J_P \ar[r]_{log} & \mathbb{Q}_p \otimes J_P}
\]
\label{relationuv}
\end{lemma}

\noindent{\bf Proof:} This is because $log$ is natural with respect to the homomorphisms of Iwasawa algebras induced by group homomorphisms. 
\qed

\begin{lemma} Let $R$ denote either the ring $\Lambda_O(Z)$ or the ring $\widehat{\Lambda_O(Z)_{(p)}}$. Let $x \in K_1'(R[G]^{\tau})$. Then for every noncyclic subgroup $P$ of $G$, we have
\[
\alpha_P(\theta^G_P(x))^{p|P|} \equiv u^G_P(\alpha(\theta^G(x))) (\text{mod } p)
\]
In particular, $log$ of $\frac{\alpha_P(\theta^G_P(x))^{p|P|}}{u^G_P(\alpha(\theta^G(x)))}$ is defined.
\label{weakcong}
\end{lemma}

\noindent{\bf Proof:} If $C$ is a nontrivial cyclic subgroup of $G$, then $\alpha_C(\theta^G_C(x)) \equiv 1 (\text{mod }p)$ and 
\[
u^G_P(\alpha(\theta^G(x))) \equiv \varphi(\alpha_{\{1\}}(\theta^G_{\{1\}}(x))) (\text{mod }p).
\]
The result will follow if we show that 
\[
\theta^G_P(x)^{|P|} \equiv \theta^G_{\{1\}}(x) (\text{mod p})
\]
But this is proposition 2.3 in \cite{SchneiderVenjakob:2010}.
\qed

\begin{remark} Let $x \in K_1'(R[G]^{\tau})$. Then for any nontrivial cyclic subgroup $P$ of $G$ 
\[
\frac{\alpha_P(\theta^G_P(x))}{u^G_P(\alpha(\theta^G(x)))} \equiv 1 (\text{mod }p).
\]
and 
\[
\frac{\alpha_{\{1\}}(\theta^G_{\{1\}}(x))}{\varphi(\theta^G_{\{1\}}(x))u^G_{\{1\}}(\alpha(\theta^G(x)))} \equiv 1 (\text{mod }p).
\]
\end{remark}

\begin{proposition} Let $x \in K_1'(\Lambda_O(\mathcal{G}))$. The additive and the multiplicative sides are related by the following formula: if $P$ is a noncyclic subgroup of $G$ then  
\[
\beta^G_P(L(x)) = \frac{1}{p^2|P|}log(\frac{\alpha_P(\theta^G_P(x))^{p|P|}}{u^G_P(\alpha(\theta^G(x)))}).
\]
If $P \in C(G)$ but $P \neq \{1\}$ then
\[
\beta^G_P(L(x)) = \frac{1}{p} log \Big( \frac{\alpha_P(\theta^G_P(x))}{u^G_P(\alpha(\theta^G(x)))} \Big).
\]
If $P = \{1\}$, then
\[
\beta^G_{\{1\}}(L(x)) = \frac{1}{p} log \Big( \frac{\alpha_{1}(\theta^G_{\{1\}}(x))}{\varphi(\theta^G_{\{1\}}(x))u^G_{\{1\}}(\alpha(\theta^G(x)))}\Big).
\]
Analogous relation holds in the case of $p$-adic completions.
\label{relation}
\end{proposition}
\noindent{\bf Proof:} Let $x \in K_1'(\Lambda_O(\mathcal{G}))$. Let $P$ be a noncyclic subgroup of $G$. Recall the definition of $\alpha$, definition \ref{defnalpha}. Consider the left hand side of the above equation
\begin{align*} 
\beta^G_P(L(x)) &= \beta^G_P(log(x) - \frac{\varphi}{p} log(x)) \\
& =  \frac{1}{p}log(\alpha_P(\theta^G_P(x))) - \beta^G_P(\frac{\varphi}{p}log(x)) \hspace{.5cm} \text{(lemmas \ref{relationthetat} and \ref{relationalphaeta})} \\
& = \frac{1}{p}log(\alpha_P(\theta^G_P(x))) - \frac{1}{p}v^G_P(\beta^G(log(x))) \hspace{.5cm} \text{(lemma \ref{relationbetaphi})} \\
& = \frac{1}{p}log(\alpha_P(\theta^G_P(x))) - \frac{1}{p^2}v^G_P(log(\alpha(\theta^G(x)))) \hspace{.5cm} \text{(lemmas \ref{relationthetat} and \ref{relationalphaeta})} \\
& = \frac{1}{p} log(\alpha_P(\theta^G_P(x)) - \frac{1}{p^2|P|}log(u^G_P(\alpha(\theta^G(x)))) \hspace{.5cm} \text{(lemma \ref{relationuv})}\\
& = \frac{1}{p^2|P|} log\Big( \frac{\alpha_P(\theta^G_P(x))^{p|P|}}{u^G_P(\alpha(\theta^G(x)))} \Big).
\end{align*}
The proof for $P \in C(G)$ but $P \neq \{1\}$ is similar. \\
Now assume that $P = \{1\}$. In this case $\beta^G_{\{1\}} = t^G_{\{1\}}$ and
\begin{align*}
\beta^G_{\{1\}}(L(x)) &= t^G_{\{1\}}(log(x)) - \frac{1}{p}t^G_{\{1\}}(\varphi(log(x))) \\
&= log(\theta^G_{\{1\}}(x)) - \frac{1}{p}(\varphi+v^G_{\{1\}})(\beta^G(log(x))) \hspace{.5cm} \text{(lemmas \ref{relationthetat} and \ref{relationbetaphi})} \\
&=log(\theta^G_{\{1\}}(x)) - \frac{1}{p}\varphi(log(\theta^G_{\{1\}}(x))) -\frac{1}{p} v^G_{\{1\}}(log(\alpha(\theta^G(x)))) \\
& = \frac{1}{p} log\Big( \frac{\theta^G_{\{1\}}(x)^p}{\varphi(\theta^G_{\{1\}}(x))u^G_{\{1\}}(\alpha(\theta^G(x)))}\Big) \hspace{.5cm} \text{(lemma \ref{relationuv})} \\
& = \frac{1}{p} log\Big(\frac{\alpha_{\{1\}}(\theta^G_{\{1\}}(x))}{\varphi(\theta^G_{\{1\}}(x))u^G_{\{1\}}(\alpha(\theta^G(x)))}\Big).
\end{align*}

Now we prove the proposition for $x \in K_1'(\widehat{\Lambda_O(\mathcal{G})_S})$. Write $x=uy$, where $y \in K_1'(\widehat{\Lambda_O(\Gamma)_{(p)}})$ and $u$ lies in the image of $K_1(\widehat{\Lambda_O(\mathcal{G})_S},J)$ in $K_1'(\widehat{\Lambda_O(\mathcal{G})_S})$. The proof for $u$ is exactly as above. In fact the proof works for any $u$ in the image of $K_1(\widehat{\Lambda_O(\mathcal{G})_S},J_R)$ (recall that $J_R$ is the Jacobson radical of $\widehat{\Lambda_O(\mathcal{G})_S}$). Hence we only need to show the result for $y$. We first show it for $y \in \widehat{\Lambda_O(Z)_{(p)}}$. In this case $\theta^G_P(y) = y^{[G:P]}$ for every $P \leq G$. Also, for any nontrivial cyclic subgroup $P$ of $G$
\[
\alpha_P(\theta^G_P(y)) = 1.
\]
Since $L(y) \in \widehat{\Lambda_O(Z)_{(p)}}$, for any noncyclic subgroup $P$ of $G$ or if $P = \{1\}$, we have
\[
\beta^G_P(L(y)) = [G:P]L(y).
\]
For any nontrivial cyclic subgroup $P$ of $G$
\[
\beta^G_P(L(y)) = 0.
\]
With these observations the result is easily verified. Now for any $y \in \widehat{\Lambda_O(\Gamma)_{(p)}}$, there is a non-negative integer $r$ such that $\varphi^r(y) \in \widehat{\Lambda_O(Z)_{(p)}}$. We know that the result holds for $\frac{y^{p^r}}{\varphi^r(y)}$ (as it is congruent to 1 mod $p$ and hence is in the image of $K_1(\widehat{\Lambda_O(\mathcal{G})_S},J_R)$) and it holds for $\varphi^r(y)$. Hence it holds for $y^{p^r}$. Since $\beta^G_P$ takes values in torsion free abelian group, the proposition holds for $y$. 
\qed

\subsection{Proof of the main algebraic theorems}  

\begin{lemma} The image of $\theta^G$ is contained in $\Phi^G$. Similarly, image of $\widehat{\theta^G_S}$ is contained in $\widehat{\Phi^G_S}$. 
\label{multiplicativecontain}
\end{lemma}

\noindent{\bf Proof:} Here we write a proof for the first claim of the lemma only. The proof for the second claim is exactly the same. We must show that for any $ x \in K_1(\Lambda_O(\mathcal{G})),$ the element $\theta^G(x)$ satisfies M1, M2, M3 and M4. \\

M1 is clear if the following diagram commutes
\[
\xymatrix{ K_1'(\Lambda_O(\mathcal{G})) \ar[rrd]^{norm^{G}_{P}} \ar[d]_{norm^G_{P'}} & & \\
K_1'(\Lambda_O(U_{P'})) \ar[rr]_{norm^{P'}_P} \ar[d] & & K_1'(\Lambda_O(U_P)) \ar[d]\\
\Lambda_O(U_{P'}^{ab})^{\times} \ar[rd]_{nr^{P'}_P} & & \Lambda_O(U_P^{ab})^{\times} \ar[ld]^{\pi^{P'}_P} \\
& \Lambda_O(U_P/[U_{P'},U_{P'}])^{\times} &}
\]
To see that this diagram commutes we note that the $\Lambda_O(U_{P})$-basis of $\Lambda_O(U_{P'})$ can be taken to be the same as the $\Lambda_O(U_P/[U_{P'},U_{P'}])$-basis of $\Lambda_O(U_{P'}^{ab})$.

 M2. Pick a lift $\tilde{x}$ of $x$ in $\Lambda_O(\mathcal{G})^{\times}$. Recall how the map $\theta^G_P$ is defined for any $P \leq G$. The ring $\Lambda_O(\mathcal{G})$ is a free $\Lambda_O(U_P)$-module of rank $[G:P]$. We take the $\Lambda_O(U_P)$-linear map on $\Lambda_O(\mathcal{G})$ defined by multiplication by $x$ on the right. Let $A_P(x)$ be the matrix (with entires in $\Lambda_O(U_P)$) of this map with respect to some basis $C(G,P)$. For any $g \in G$, the set $g C(G,P) g^{-1}$ is a basis for the $\Lambda_O(U_{gPg^{-1}})$-module $\Lambda_O(\mathcal{G})$. The matrix for the $\Lambda_O(U_{gPg^{-1}})$-linear map on $\Lambda_O(\mathcal{G})$ induced by multiplication on the right by $x$ with respect to this basis is $gA_P(x)g^{-1}$. Hence
\[
g\theta^G_P(x)g^{-1} = \theta^G_{gPg^{-1}}(x).
\]

M3. Let $P \leq P' \leq G$ be such that $[P':P]=p$. Let $n^{P'}_P$ be the norm map from $K'_1(\Lambda_O(U_{P'}))$ to $K'_1(\Lambda_O(U_P))$. Consider the following diagram
\[
\xymatrix{ K'_1(\Lambda_O(\mathcal{G})) \ar[d]_{norm} \ar[rd]^{norm} & \\
K'_1(\Lambda_O(U_{P'})) \ar[r]^{n^{P'}_P} \ar[d]_{\pi_{P'}} & K'_1(\Lambda_O(U_P)) \ar[d]^{\pi_P} \\
\Lambda_O(U_{P'}^{ab})^{\times} & \Lambda_O(U_P^{ab})^{\times}}
\]
It is enough to show that for any $x \in K'_1(\Lambda_O(U_{P'}))$
\[
ver^{P'}_P(\pi_{P'}(x)) \equiv \pi_P(n^{P'}_P(x)) (\text{mod } T_{P,P'}).
\]
Let $g \in U_{P'}-U_P$. Then $\{1,g,g^2,\cdots, g^{p-1}\}$ forms a basis the $\Lambda_O(U_P)$-module $\Lambda_O(U_{P'})$. Write $x \in \Lambda_O(U_{P'})$ as $x = \sum_{k=0}^{p-1} x_k g^k$ with $x_k \in \Lambda_O(U_P)$ for all $k$. Let $\sigma$ denote the automorphism of $\Lambda_O(U_P)$ given by $y \mapsto gyg^{-1}$. Then norm of $x$ in $K'_1(\Lambda_O(U_P))$ is class of the matrix
\[
\left( \begin{array}{cccc}
x_0 & x_1 & \cdots & x_{p-1} \\
g^p\sigma(x_{p-1}) & \sigma(x_0) & \cdots & \sigma(x_{p-2}) \\
\vdots &  & & \vdots \\
g^p\sigma^{p-1}(x_1) & g^p\sigma^{p-1}(x_2) & \cdots & \sigma^{p-1}(x_0)
\end{array}
\right).
\]
$\pi_P(n^{P'}_P(x))$ is determinant of this matrix. Hence if $S_p$ denotes the group of symmetries of the set $\{0,1,2,\ldots, p-1\}$ and for any $\delta \in S_p$, let $e_{\delta}$ be the number of $k$'s such that $\delta(k) < k$, then
\[
\pi_P(n^{P'}_P(x)) = \sum_{\delta \in S_p}g^{e_{\delta}p} \Big( \prod_{k=0}^{p-1} \sigma^k(\pi_P(x_{\delta(k)-k})) \Big).
\]
If $\sigma(\prod_{k=0}^{p-1}\sigma^k(\pi_P(x_{\delta(k)-k}))) = \prod_{k=0}^{p-1}\sigma^k(\pi_P(x_{\delta(k)-k}))$, then $\delta(k)-k$ is constant for all $0 \leq k \leq p-1$. Hence 
\[
\pi_P(n^{P'}_P(x)) \equiv \sum_{i=0}^{p-1} g^{ip} \prod_{k=0}^{p-1} \sigma^k(\pi_P(x_i)) (\text{mod } T_{P,P'})
\]
Now for each $0 \leq i \leq p-1$ write $x_i = \sum_{h \in P} a_hh$, with each $a_h \in \Lambda_O(Z)$. Then $ver^{P'}_P(\pi_{P'}(x_i)) = \pi_P(\sum_{h\in P} \varphi(a_h) \prod_{k=0}^{p-1}\sigma^k(h))$. On the other hand 
\begin{align*}
\prod_{k=0}^{p-1} \sigma^k(\pi_P(x_i)) & = \prod_{k=0}^{p-1} \pi_P(\sigma^k(x_i)) \\
& = \prod_{k=0}^{p-1} \pi_P(\sigma^k(\sum_{h \in P}a_hh)) \\
& = \pi_P( \prod_{k=0}^{p-1} (\sum_{h \in P} a_h\sigma^k(h))) \\
& \equiv \pi_P(\sum_{h \in P} a_h^p \prod_{k=0}^{p-1} \sigma^k(h)) (\text{mod } T_{P,P'}) \\
& \equiv \pi_P(\sum_{h \in P} \varphi(a_h) \prod_{k=0}^{p-1} \sigma^k(h)) (\text{mod } T_{P,P'}) \\
& = ver^{P'}_P(\pi_{P'}(x_i)).
\end{align*}
Hence we have
\begin{align*}
ver^{P'}_P(\pi_{P'}(x)) & = \sum_{i=0}^{p-1} ver^{P'}_P(\pi_{P'}(x_i))g^{ip} \\
& \equiv \sum_{i=0}^{p-1} \prod_{k=0}^{p-1}\sigma^k(\pi_P(x_i)) g^{ip} (\text{mod } T_{P,P'})\\
& = \pi_P(n^{P'}_P(x)).
\end{align*}

M4. Let $P$ be a nontrivial cyclic subgroup of $G$. By proposition \ref{relation} we have
\[
\beta^G_P(L(x)) = \frac{1}{p}log\Big( \frac{\alpha_P(\theta^G_P(x))}{u^G_P(\alpha(\theta^G(x)))} \Big).
\]
By lemma \ref{additivecontain} we get that $p\beta^G_P(L(x)) \in pT_P$. But $log$ induces an isomorphism between $1+pT_P$ and $pT_P$. Hence M4 is satisfied for nontrivial cyclic subgroups of $G$. The verification of M4 for $P = \{1\}$ is similar.

\qed

\begin{lemma} There is a map $\mu(O) \times \mathcal{G}^{ab} \rightarrow \prod_{P\leq G} \Lambda_O(U_P^{ab})^{\times}$ which is injective and fits in the following commutative diagram
\[
\xymatrix{ \mu(O)\times \mathcal{G}^{ab} \ar@{^{(}->}[r] \ar@{=}[d] & K'_1(\Lambda_O(\mathcal{G})) \ar[d]_{\theta^G} \\
\mu(O)\times \mathcal{G}^{ab} \ar@{^{(}->}[r] & \prod_{P\leq G} \Lambda_O(U_P^{ab})^{\times} }
\]
\end{lemma}

\noindent{\bf Proof:} We define the claimed map as the composition
\[
\mu(O) \times \mathcal{G}^{ab} \xrightarrow{id \times ver} \mu(O) \times U_P \hookrightarrow \Lambda_O(U_P^{ab})^{\times}.
\]
Injectivity is obvious. Since the norm homomorphism on $K_1$ groups is same as the transfer when restricted to $\mathcal{G}^{ab}$, we also get the commutativity of the diagram in the lemma.

\qed

\begin{lemma} Recall that we denote $\psi^G_{\Lambda_O(Z)}$ from the additive theorem by $\psi^G$. Then there is a surjective map 
\[
\omega : \psi^G \rightarrow  \mathcal{G}^{ab},
\]
which makes the following diagram commute
\[
\xymatrix{ \Lambda_O(Z)[Conj(G)]^{\tau} \ar[rr]^{\omega} \ar[d]_{\beta^G} & &  \mathcal{G}^{ab} \ar[d] \\
\psi^G \ar[rr]_{\omega} & & \mathcal{G}^{ab} }
\]
where the map $\omega$ in the top row is the one from definition \ref{integrallogforlambda}.
\end{lemma}

\noindent{\bf Proof:} This is a trivial corollary of theorem \ref{additivetheorem}. 

\qed

\begin{lemma} The map $\mathcal{L} :=(\mathcal{L}_P) : \Phi^G \rightarrow \psi^G$, defined by 
\[
\mathcal{L}_P((x_C)) = \left\{ 
\begin{array}{l l}
\frac{1}{p^2|P|}log\Big(\frac{\alpha_P(x_P)^{p|P|}}{u^G_P(\alpha((x_C)))}\Big) & \text{if $P$ is not cyclic} \\
\frac{1}{p}log\Big(\frac{\alpha_P(x_P)}{u^G_P(\alpha((x_C)))}\Big) & \text{if $P$ is cyclic and $P \neq \{1\}$} \\
\frac{1}{p} log \Big( \frac{\alpha_{\{1\}}(x_{\{1\}})}{\varphi(x_{\{1\}})u^G_{\{1\}}(\alpha((x_C)))} \Big) & \text{if $P= \{1\}$}
\end{array} \right.
\]
sits in the following exact sequence
\[
1 \rightarrow \mu(O) \times \mathcal{G}^{ab} \rightarrow \Phi^G \xrightarrow{\mathcal{L}} \psi^G \xrightarrow{\omega}  \mathcal{G}^{ab} \rightarrow 1.
\]
The map $\mu(O) \times \mathcal{G}^{ab} \rightarrow \Phi^G \subset \prod_{P \leq G} \Lambda_O(U_P^{ab})^{\times}$ is the composition 
\[
\mu(O) \times \mathcal{G}^{ab} \rightarrow \prod_{P \leq G} \mu(O) \times U_P^{ab} \hookrightarrow \prod_{P \leq G} \Lambda_O(U_P^{ab})^{\times},
\]
where the first map is identity on $\mu(O)$ and the transfer homomorphism from $\mathcal{G}^{ab}$ to $U_P^{ab}$ for each $P \leq G$.
\end{lemma}

\noindent{\bf Proof:} Image of $\mathcal{L}$ is contained in $\prod_{P\leq G} \mathbb{Q}_p \otimes \Lambda(U_P^{ab})$. We must verify that it is actually contained in $\psi^G$, i.e. it lies in $\prod_{P \leq G} \Lambda(U_P^{ab})$ and satisfies A1, A2, and A3. We first prove the following

\begin{lemma} Let $P \leq P' \leq G$ such that $[P':P]=p$. Let $C \in C(G)$ be such that $C^p$ is contained in $P'$ but not in $P$. Then
\[
nr^{P'}_P(\varphi(\alpha_C(x_C))) = 1
\]
in $\Lambda(U_P/[U_{P'},U_{P'}])$.
\end{lemma}

\noindent{\bf Proof:} Recall that
\[
\alpha_C(x_C) = \frac{x_C^p}{\prod_{k=0}^{p-1}\omega_C^k(x_C)}.
\]
Hence
\begin{align*}
\varphi(\alpha_C(x_C)) & = \frac{\varphi(x_C)^p}{\prod_{k=0}^{p-1}\varphi(\omega_C^k(x_C))} \\
& = \frac{\varphi(x_C)^p}{\prod_{k=0}^{p-1}\omega^k_{C^p}(\varphi(x_C))}
\end{align*}
But $nr^{P'}_P(\varphi(x_C)) = \prod_{k=0}^{p-1} \omega^k_{C^p}(\varphi(x_C))$ by lemma (\ref{normandtrace}). Hence we get the result. 
\qed

Using M1 we now show that the image of $\mathcal{L}$ satisfies A1. Let $P \leq P' \leq G$ such that $[P', P'] \leq P$ and if $P$ is a non-trivial cyclic group then $[P',P'] \neq P$. Then we consider following 3 cases: \\
Case 1: $P$ is not cyclic. Then
\begin{align*}
tr^{P'}_P(\mathcal{L}_{P'}((x_C))) & = tr^{P'}_P\Big(\frac{1}{p^2|P'|}log\Big(\frac{\alpha_{P'}(x_{P'})^{p|P'|}}{u^G_{P'}(\alpha((x_C)))}\Big)\Big) \\
& = tr^{P'}_P\Big(\frac{1}{p^2|P'|}log\Big(\frac{(x_{P'})^{p^2|P'|}}{\prod_{C'} \varphi(\alpha_{C'}(x_{C'}))^{|C'|}}\Big) \Big) \\
& = \frac{1}{p^2|P'|} log \Big( \frac{nr^{P'}_P(x_{P'})^{p^2|P'|}}{nr^{P'}_P(\prod_{C'} \varphi(\alpha_{C'}(x_{C'}))^{|C'|})} \Big) \\
& = \frac{1}{p^2|P'|} log \Big( \frac{\pi^{P'}_P(x_P)^{p^2|P'|}}{\prod_{C} \varphi(\alpha_C(x_C))^{p|C|} } \Big) \\
& = \pi^{P'}_P( \frac{1}{p^2|P|} log \Big( \frac{\alpha_P(x_P)^{p|P|}}{\prod_{C} \varphi(\alpha_C(x_C))^{|C|}} \Big) \\
& = \pi^{P'}_P(\mathcal{L}_P((x_C))). 
\end{align*}
In the products $C'$ runs through all cyclic subgroups of $G$ such that $C'^p \leq P'$ and $C$ runs through all cyclic subgroups of $G$ such that $C^p \leq P$. \\
Case 2: When $P$ is cyclic but $P'$ is not. Then
\begin{align*}
\eta_P(tr^{P'}_P(\mathcal{L}_{P'}((x_C)))) &= \pi^{P'}_P \Big( \eta_P \Big[ \frac{1}{p^2|P|} log \Big( \frac{x_P^{p^2|P|}}{\prod_C \varphi(\alpha_C(x_C))^{|C|}} \Big) \Big] \Big)
\end{align*}
Here $C$ runs through all cyclic subgroups of $G$ such that $C^p \leq P$. Now note that 
\[
\alpha_P(\varphi(\alpha_C(x_C))) = \left\{ \begin{array}{l l} 
1 & \text{if $C^p \neq P$} \\
\varphi(\alpha_C(x_C))^p & \text{if $C^p=P$}
\end{array} \right.
\]
Hence
\begin{align*}
\eta_P(tr^{P'}_P(\mathcal{L}_{P'}((x_C)))) & = \pi^{P'}_P\Big[ \frac{1}{p} log \Big( \frac{\alpha_P(x_P)}{\prod_C \varphi(\alpha_C(x_C))} \Big) \Big] \\
& = \pi^{P'}_P(\mathcal{L}_P((x_C))).
\end{align*}
Here in the product $C$ runs through all cyclic subgroup of $G$ such that $C^p =P$. \\
Case 3: When $P'$ is cyclic. This case is easy by above lemma. 

A2 and A3 are clear from M2 and M4 respectively. We now show that image of $\mathcal{L}$ is actually contained in $\prod_{P \leq G} \Lambda(U_P^{ab})$. As $\mathcal{L}_P((x_C)) \in T_P$ for all $P \in C(G)$ this follows from proposition \ref{additiveintegrality}. Hence image of $\mathcal{L}$ is contained in $\psi^G$. 

Let us show exactness at $\Phi^G$. Let $(x_P) \in Ker(\mathcal{L})$. Since $log$ is an isomorphism on $pT_P$, we have
\[
\frac{\alpha_P(x_P)}{u^G_P(\alpha((x_C)))}=1,
\]
for all $P \in C(G)$. Consider first the cyclic subgroup $P$ for which there is no $P' \in C(G)$ with $P'^p=P$. Then the denominator of the above expression is 1. Hence we have $\alpha_P(x_P)=1$. By induction we get that $\alpha_P(x_P) = 1$ for all $P \in C(G)$ and consequently $\alpha_P(x_P) =1$ for all $P \leq G$. The only torsion in $\Lambda_O(U_P^{ab})^{\times}$ is $\mu(O) \times U_P^{ab}$ by a theorem of Higman \cite{Higman:1940}. But for $P$ non-cyclic, $\alpha_P(x_P) = x_P^p$. Hence $x_P \in \mu(O) \times U_P^{ab}$ for all non-cyclic subgroups $P$. Also, $x_{\{1\}}^p/\varphi(x_{\{1\}}) =1$, hence by the exact sequence in definition (\ref{integrallogforlambda}) $x_{\{1\}} \in \mu(O) \times U_{\{1\}}$.

Let $P \leq P' \in C(G)$ with $[P':P] =p$. Then $U_P$ is an index $p$ subgroup of $U_{P'}$. Hence $nr^{P'}_P(x) = \prod_{k=0}^{p-1}\omega_{P'}^k(x)$ by lemma (\ref{normandtrace}). This is precisely the denominator of $\alpha_{P'}(x)$ for any nontrivial cyclic subgroup $P'$ of $G$. Using M1 and using induction on the order of cyclic subgroups (this time starting with the trivial subgroup) we conclude that $x_P \in \mu(O) \times U_P$ for all cyclic subgroups $P$. Now we use induction on index of subgroup $P$ in $G$. Let $P \leq G$ be of index $p$. Then M3 implies that $ver^G_P(x_G) \equiv x_P (\text{mod } T_{P,G})$. Since $x_G$ and $x_P$ are in $\mu(O) \times \mathcal{G}^{ab}$ and $\mu(O) \times U_P^{ab}$ respectively, the congruence implies equality. One proceeds by induction on $[G:P]$ to conclude that $(x_P)$ must be in the image of $\mu(O) \times \mathcal{G}^{ab}$. This shows that $ker(\mathcal{L})$ is contained in the image of $\mu(O) \times \mathcal{G}^{ab}$. The converse is trivial.

Next we show exactness at $\psi^G$. First note that the map
\[
\omega: \psi^G \rightarrow \mathcal{G}^{ab}
\]
is given by $\omega((a_P)_{P\leq G}) = \omega(a_G) \in \mathcal{G}^{ab}$, where the second $\omega$ is 
\[
\omega: \Lambda_O(\mathcal{G}^{ab}) \rightarrow \mathbb{G}^{ab}
\]
in definition (\ref{integrallogforlambda}).

\noindent{\bf Claim:} For any $(x_C)_C \in \Phi^G$, we have $\mathcal{L}_G((x_C)_C) = \frac{1}{p}log\Big(\frac{x_G^p}{\varphi(x_G)}\Big)$.

\noindent{\bf Proof of the claim:} Let $\mathcal{L}_G((x_C)_C) = (a_P)_{P\leq G} \in \psi^G$. Then 
\[
a_G = \sum_{P \in C(G)}\frac{1}{[G:P]}a_P,
\]
in $\Lambda_O(\mathcal{G}^{ab})$. Hence
\begin{align*}
\mathcal{L}_G((x_C)_C) & = \sum_{P \in C(G)} \frac{1}{[G:P]} \mathcal{L}_P((x_C)_C)\\
\frac{1}{p^2|G|}log\Big(\frac{x_G^{p^2|G|}}{u_G(\alpha((x_C)_C))}\Big) &= \frac{1}{p|G|}log\Big(\frac{x_{\{1\}}^p}{\varphi(x_{\{1\}}\prod_{|P'|=p}\varphi(\alpha_{P'}(x_{P'}))} \Big) \\ &+ \sum_{P \in C(G)} \frac{1}{p[G:P]} log\Big(\frac{\alpha_P(x_P)}{\prod_{P'\in C(G):[P':P]=p}\varphi(\alpha_{P'}(x_{P'}))}\Big) \\
&= \frac{1}{p^2|G|} log \Big( \frac{x_{\{1\}}^{p^2}}{\varphi(x_{\{1\}})^p\prod_{P'}\varphi(\alpha_{P'}(x_{P'}))^p}\Big) \\
&+ \frac{1}{p^2|G|}log\Big( \prod_{P} \frac{\alpha_P(x_P)^{p|P|}}{\prod_{P'}\varphi(\alpha_{P'}(x_{P'}))^{p|P|}}\Big) \\
& = \frac{1}{p^2|G|} log \Big(\frac{\prod_{P\in C(G)} \alpha_P(x_P)^{p|P|}}{\varphi(x_{\{1\}})^p \prod_{P \in C(G): P\neq\{1\}}\varphi(\alpha_P(x_P))^{|P|}}\Big) \\
&= \frac{1}{p^2|G|} log\Big(\frac{\prod_{P \in C(G)} \alpha_P(x_P)^{p|P|}}{u_G(\alpha((x_C)_C))}\Big).
\end{align*}
Hence
\[
log(x_G^{p|G|}) = log(\prod_{P \in C(G)} \alpha_P(x_P)^{|P|}).
\]
Applying $\varphi$ to both side gives 
\[
log(\varphi(x_G)^{p|G|}) = log(\prod_{P \in C(G)} \varphi(\alpha_P(x_P))^{|P|}) = log(u_G(\alpha(x_C)_C)).
\]
Hence
\[
\mathcal{L}_G((x_C)_C) = \frac{1}{p^2|G|}log\Big(\frac{x_G^{p^2|G|}}{u_G((\alpha(x_C)_C))}\Big) = \frac{1}{p} log\Big(\frac{x_G^p}{\varphi(x_G)}\Big).
\]
This finishes proof of the claim. Hence 
\[
\omega(\mathcal{L}((x_C)_C)) = \omega(\mathcal{L}_G((x_C)_C)) = \omega\Big(\frac{1}{p}log\Big(\frac{x_G^p}{\varphi(x_G)}\Big)\Big) = 1
\]
by exact sequence in definition \ref{integrallogforlambda}. Hence the image of $\mathcal{L}$ is contained in the kernel of $\omega$. On the other hand, let $a \in \psi^G$ be in the kernel of $\omega$. Then by the previous lemma, the additive theorem and the exact sequence in definition \ref{integrallogforlambda} we get $x \in K_1'(\Lambda_O(\mathcal{G}))$ such that $\beta^G(L(x)) = a$. We finish the proof by using the commutativity of the following diagram
\[
\xymatrix{ K_1'(\Lambda_O(\mathcal{G})) \ar[rr]^{L} \ar[d]_{\theta^G} & & \Lambda_O(Z)[Conj(G)]^{\tau} \ar[d]^{\beta^G} \\
\Phi^G \ar[rr]_{\mathcal{L}} & & \psi^G}
\]
Hence we have $\mathcal{L} (\theta^G(x)) = a$ and we get the exactness at $\psi^G$.	

\qed

\noindent{\bf Proof of theorem \ref{theorem1}:} From the previous lemmas we have the following commutative diagram
\[
\xymatrix{ \mu(O) \times \mathcal{G}^{ab} \ar[r] \ar@{=}[d] & K'_1(\Lambda_O(\mathcal{G})) \ar[r] \ar[d]_{\theta^G} & \Lambda_O(Z)[Conj(G)]^{\tau} \ar[r] \ar[d]^{\beta^G}_{\sim} &  \mathcal{G}^{ab} \ar@{=}[d]  \\
\mu(O) \times \mathcal{G}^{ab} \ar[r] & \Phi^G \ar[r] & \psi^{G} \ar[r] &  \mathcal{G}^{ab} }
\]
The theorem now follows by the five lemma. 

\qed

\noindent{\bf Proof of theorem \ref{theorem2}:} It follows from lemma \ref{multiplicativecontain} and the fact  
\[
\widehat{\Phi^G_S} \cap \prod_{P \leq G} \Lambda_O(U_P)_S^{\times} = \Phi^G_S,
\]
that the image of $\theta^G_S$ is contained in $\Phi^G_S$. Moreover, from 
\[
\Phi^G_S \cap \prod_{P \leq G} \Lambda_O(U_P^{ab})^{\times} = \Phi^G,
\] 
and theorem \ref{theorem1}, we obtain that 
\[
\Phi^G_S \cap \prod_{P \leq G} \Lambda_O(U_P^{ab})^{\times} = Im(\theta^G).
\]

\qed

\begin{corollary} Let $\mathcal{G}$ be a compact $p$-adic Lie group with a quotient isomorphic to $\mathbb{Z}_p$ and let $O$ be the ring of integers in a finite unramified extension of $\mathbb{Q}_p$. Then $K_1'(\Lambda_O(\mathcal{G}))$ injects in $K_1'(\Lambda_O(\mathcal{G})_S)$.
\end{corollary}

\section{Arithmetic side of the proof}
\label{sectioncongruences} To prove the main theorem (theorem \ref{mainconjecture}) we need to prove (by virtue of the reductions in section \ref{sectionreductions}; more precisely theorems \ref{reductiontodim1}, \ref{reductiontohyperelementary}, \ref{lhyperelementarymainconjecture}, \ref{theoremreductiontoprop}) the main conjecture for admissible $p$-adic Lie extension $F_{\infty}/F$ satisfying $\mu=0$ hypothesis such that $Gal(F_{\infty}/F)  = \Delta \times \mathcal{G}$, where $\mathcal{G}$ is a pro-$p$ p-adic Lie group of dimension one and $\Delta$ is a finite cyclic group of order prime to $p$.

\begin{notation} Recall that we have fixed an open central pro-cyclic subgroup $Z$ of $\mathcal{G}$ and put $G= \mathcal{G}/Z$. Also recall $C(G)$ is the set of cyclic subgroups of $G$ and for every $P \leq G$ we put $U_P$ to be the inverse image of $P$ in $\mathcal{G}$. We denote the field $F_{\infty}^{\Delta \times U_P}$ by $F_P$ and denote the field $F_{\infty}^{[U_P,U_P]}$ by $K_P$. Then $K_P/F_P$ is an abelian extension with $Gal(K_P/F_P) = \Delta \times U_P^{ab}$. We denote the field $F_{\infty}^{\Delta \times N_{\mathcal{G}}U_P}$ by $F_{W_GP}$. Note that $F_{W_GP} \subset F_P$ and $Gal(F_P/F_{W_GP}) \cong W_GP$. We denote the Deligne-Ribet, Cassou-Nogues, Barsky $p$-adic zeta function for the abelian extension $K_P/F_P$ by $\zeta_P$. It is an element in $\Lambda(\Delta \times U_P^{ab})_S^{\times}$. Let $\zeta_0$ be the $p$-adic zeta function of the extension $F_{\infty}/F_{\infty}^{\Delta \times Z^p}$. Recall that we have fixed a finite set $\Sigma$ of finite primes of $F$ containing all primes which ramify in $F_{\infty}$. Let $\Sigma_P$ denote the set of primes of $F_P$ lying above $\Sigma$. Let $r_P := [F_P:\mathbb{Q}]$ and $d_P = [F_P:F]$. If a group $P$ acts on a set $X$, then we denote the stabiliser of $x \in X$ by $P_x$. 
\end{notation}

\subsection{The strategy of Burns and Kato} 

\begin{lemma} Let $\rho$ be an irreducible Artin representation of $\mathcal{G}$. Then there is a one dimensional representation $\chi$ of $\mathcal{G}$ inflated from $\Gamma$ such that $\rho \otimes \chi$ is trivial on $Z$.
\end{lemma}

\noindent{\bf Proof:} We use induction on the order of $\mathcal{G}/Z$. By proposition 24 in Serre \cite{Serre:representationtheory} either  \\
a) $\rho$ restricted to $Z$ is isotypic (i.e. direct sum of isomorphic irreducible representations) OR\\
b) $\rho$ is induced from an irreducible representation of a proper subgroup $A$ of $\mathcal{G}$ containing $Z$. 

In case a) let $\rho|_{\Gamma} = \oplus_{i}\chi_i$ Define $\chi = \chi_i^{-1}$ for any $i$ (note that $\chi_i|_Z = \chi_j|_Z$ for any $i,j$). Then $\rho \otimes \chi$ is trivial on $Z$.

In case b) Say $\rho = Ind^{\mathcal{G}}_A(\eta)$. Let $r$ be such that image of $A$ in $\Gamma$ is $\Gamma^{p^r}$. By induction hypothesis we can find a $\chi$ inflated from $\Gamma^{p^r}$ such that $\eta \otimes \chi$ is trivial on $Z$. We may extend $\chi$ to $\tilde{\chi}$ on $\Gamma$. Then
\[
Ind^{\mathcal{G}}_A(\eta \otimes \chi) = Ind^{\mathcal{G}}_A(\eta) \otimes \tilde{\chi} = \rho \otimes \tilde{\chi}.
\]
Since $\eta \otimes \chi |_{Z}$ is trivial and $Z$ is central, $Ind^{\mathcal{G}}_A(\eta \otimes \chi)|_Z = (\rho \otimes \tilde{\chi})|_Z$ is trivial.
\qed

\begin{proposition} With the notations as above, the main conjecture is true for $F_{\infty}/F$ if and only if $(\zeta_P)_{P} \in \Phi^{G}_S$. Note that here $\Phi^G_S = \Phi^G_{\mathbb{Z}_p[\Delta], S}$.
\label{propburnsstrategy}
\end{proposition}

\noindent{\bf Proof:} Let $f \in K_1'(\Lambda(\Delta \times \mathcal{G})_S)$ be any element such that 
\[
\partial(f) = -[C(F_{\infty}/F)].
\]
Let $\theta^G_S(f)=(f_P)_P$ in $\prod_{P \leq G} \Lambda(\Delta \times U_P^{ab})_S^{\times}$. Then $(f_P)_P \in \Phi_S^{G}$  by theorem \ref{theorem2}. Define $u_P = \zeta_Pf_P^{-1}$. As $\partial(f_P) = \partial(\zeta_P) = -[C(K_P/F_P)]$, we have $u_P \in K_1'(\Lambda(\Delta \times U_P^{ab}))$. Moreover, if $(\zeta_P)_P \in \Phi_S^{G}$, then $(u_P)_P \in \Phi^{G}$. Then by theorem \ref{theorem1} there is a unique $u \in K_1'(\Lambda(\Delta \times \mathcal{G}))$ such that $\theta(u) = (u_P)$. Define $\zeta= \zeta(F_{\infty}/F)  =uf$. We claim that $\zeta$ is the $p$-adic zeta function satisfying the main conjecture for $F_{\infty}/F$. It is clear that $\partial(\zeta) = -[C(F_{\infty}/F)]$. We now show the interpolation property. Let $\rho$ be an irreducible Artin representation of $\mathcal{G}$. Let $\sigma$ be a one dimensional representation of $\mathcal{G}$ given by the previous lemma i.e. such that $\rho \otimes \sigma$ is trivial on $Z$. Then $\rho \otimes \sigma = Ind^{\mathcal{G}}_{U_P}(\eta)$ for some $P \leq G$ (by theorem 16 Serre \cite{Serre:representationtheory}). We denote the restriction of $\sigma$ to $U_P$ by the same letter $\sigma$. Hence $\rho = Ind^{\mathcal{G}}_{U_P}(\eta) \otimes \sigma^{-1} = Ind^{\mathcal{G}}_{U_P}(\eta \otimes (\sigma^{-1}))$. Then for any character $\chi$ of $\Delta$ and any positive integer $r$ divisible by $[F_{\infty}(\mu_p):F_{\infty}]$, we have
\begin{align*} 
\zeta(\chi\rho\kappa_F^r) &  = \zeta_P(\chi \eta \sigma^{-1} \kappa_{F_P}^r) \\
& = L_{\Sigma_P}(\chi\eta\sigma^{-1}, 1-r) \\
& = L_{\Sigma}(\chi \rho, 1-r) 
\end{align*}
Hence $\zeta$ satisfies the required interpolation property. \qed

Hence we need to show the following 

\begin{theorem} The tuple $(\zeta_P)_{P\leq G}$ in the set $\prod_{P \leq G} \Lambda(\Delta \times U_P^{ab})_S^{\times}$ satisfies \\
M1. For all $P \leq P' \leq G$ such that $[P',P'] \leq P$, we have 
\[
nr^{P'}_P(\zeta_{P'}) = \pi^{P'}_P(\zeta_P).
\]
M2. The tuple $(\zeta_P)_P$ is fixed by every $g \in \mathcal{G}$ under the conjugation action. \\
M3. For $P \leq P' \leq G$ such that $[P':P]=p$, we have
\[
ver^{P'}_P(\zeta_{P'}) \equiv \zeta_P (\text{mod } T_{P,P',S}).
\]
M4. For every nontrivial cyclic subgroup $P$ of $G$ the following congruence
\[
\frac{\alpha_P(\zeta_P)}{\prod \varphi(\alpha_{P'}(\zeta_{P'}))} \equiv 1 (\text{mod } pT_{P,S}),
\]
where the product is over all $P' \in C(G)$ such that $P'^p = P$ and $P' \neq P$, holds. \\
For $P = \{1\}$ the following congruence holds
\[
\frac{\zeta_{\{1\}}^p}{\varphi(\zeta_{\{1\}}) \prod_{P'}\varphi(\alpha_{P'}(\zeta_{P'}))} \equiv 1 (\text{mod }pT_{\{1\},S}),
\]
where $P'$ runs through nontrivial cyclic subgroups of $G$ of order $p$. 
\label{theoremrequiredcongruences}
\end{theorem}

\begin{proposition} The tuple $(\zeta_P)_P$ in the theorem satisfies M1. and M2.
\label{m1m2}
\end{proposition}
\noindent{\bf Proof} Let $P \leq P'$ be two subgroups of $G$ such that $[P',P'] \leq P$. Then we must show that 
\[
nr^{P'}_P(\zeta_{P'}) = \pi^{P'}_P(\zeta_P)
\]
in $\Lambda(\Delta \times U_P/[U_{P'},U_{P'}])_S$. Let $\rho$ be an irreducible Artin representation of $U_P/[U_{P'},U_{P'}]$ and let $r$ be any positive integer divisible by $[F_{\infty}(\mu_p):F_{\infty}]$. Then for any character $\chi$ of $\Delta$, we have
\begin{align*}
nr^{P'}_P(\zeta_{P'})(\chi\rho \kappa_{F_P}^r) & = \zeta_{P'}(\chi Ind^{U_{P'}}_{U_P} (\rho) \kappa_{F_{P'}}^r) \\
& = L_{\Sigma_{P'}}(\chi Ind^{U_{P'}}_{U_P}(\rho), 1-r) \\
& = L_{\Sigma_P}(\chi \rho, 1-r)
\end{align*}
On the other hand
\begin{align*}
\pi^{P'}_P(\zeta_P)(\chi \rho \kappa_{F_P}^r) & = \zeta_{P}(\chi \rho\kappa_{F_P}^r) \\
&= L_{\Sigma_P}(\chi\rho, 1-r).
\end{align*}
Since both $nr^{P'}_P(\zeta_{P'})$ and $\pi^{P'}_P(\zeta_P)$ interpolate the same values on a dense subset of representations of $\Delta \times U_P/[U_{P'},U_{P'}]$, they must be equal. This shows that the tuple $(\zeta_P)_P$ satisfies M1. 

Next we show that the tuple $(\zeta_P)_P$ satisfies M2 i.e. for all $g \in \mathcal{G}$
\[
g(\zeta_P)g^{-1} = \zeta_{gPg^{-1}}
\]
in $g\Lambda(\Delta \times U_P^{ab})g^{-1} = \Lambda(\Delta \times gU_P^{ab}g^{-1}) = \Lambda(\Delta \times U_{gPg^{-1}}^{ab})$. We let $\rho$ be any one dimensional Artin representation of $gU_Pg^{-1}$ and $r$ be any positive integer divisible by $[F_{\infty}(\mu_p):F_{\infty}]$. Then for any character $\chi$ of $\Delta$, we have
\begin{align*} 
g(\zeta_P)g^{-1}(\chi\rho\kappa_{F_{gPg^{-1}}}^r) & = \zeta_P(\chi g\rho g^{-1} \kappa_{F_P}^r) \\
& = L_{\Sigma_P}(\chi g\rho g^{-1}, 1-r) \\
& = L_{\Sigma}(\chi Ind^{\mathcal{G}}_{U_P}(g\rho g^{-1}), 1-r)
\end{align*}
On the other hand
\begin{align*}
\zeta_{gPg^{-1}}(\chi \rho \kappa_{F_{gPg^{-1}}}^r) & = L_{\Sigma_P}(\chi \rho, 1-r) \\
& = L_{\Sigma}(\chi Ind^{\mathcal{G}}_{U_{gPg^{-1}}}(\rho), 1-r).
\end{align*}
But $Ind^{\mathcal{G}}_{U_P}(g\rho g^{-1}) = Ind^{\mathcal{G}}_{U_{gPg^{-1}}}(\rho)$. Hence $g(\zeta_P)g^{-1}$ and $\zeta_{gPg^{-1}}$ interpolate the same values on a dense subset of representations of $\Delta \times gU_P^{ab}g^{-1}$ and so must be equal. This proves that the tuple $(\zeta_P)_P$ satisfies M2.
\qed

The rest of the paper is devoted to proving that $(\zeta_P)_P$ satisfies M3 and M4. 

\subsection{Basic congruences} \label{subsectionbasiccongruences} The congruence M4 is multiplicative and does not yield directly to the method of Deligne-Ribet. In this section we state certain additive congruences which yield to the Deligne-Ribet method as we show in the following sections. These congruences are then used in the last section to prove M4. 

Let $\mathfrak{p}$ be the maximal ideal of $\mathbb{Z}_p[\mu_p]$.

\begin{proposition} For every $P \leq P' \leq G$ such that $[P':P] = p$, we have
\begin{equation} 
ver^{P'}_P(\zeta_{P'}) \equiv \zeta_P (\text{mod }T_{P,P',S}).
\label{congruence1}
\end{equation}
\label{propcongruence1}
\end{proposition}

\begin{proposition} For every $P \in C(G)$ such that there is no $P' \in C(G)$ with $P \subset P'$, we have
\begin{equation}
\zeta_P \equiv \omega_P^k(\zeta_P) (\text{mod }\mathfrak{p}T_{P,S}).
\label{congruence2}
\end{equation}
\label{propcongruence2}
\end{proposition}

\begin{proposition} If $P \in C(G)$ such that there is a $P' \in C(G)$ with $P'^p =P$, we have
\begin{equation}
\zeta_P \equiv \sum_{P'} \varphi(\zeta_{P'}) (\text{mod }T_{P,S}),
\label{congruence3}
\end{equation}
where $P'$ runs through all cyclic subgroups of $G$ such that $P' \neq P$ and $P'^p = P$.
\label{propcongruence3}
\end{proposition}

\begin{proposition} If $P \in C(G)$ and $P' \in C(G)$ such that $P' \neq P$ and $P'^p = P$, then 
\begin{equation}
\zeta_P \equiv \varphi(\zeta_{P'}) (\text{mod } T_{P, N_{G}P',S}).
\label{congruence4}
\end{equation}
\begin{equation} 
\zeta_0 \equiv \varphi(\zeta_{\{1\}}) (\text{mod } p|G|)
\label{congruence5}
\end{equation}
\label{propcongruence4}
\end{proposition}

The congruence (\ref{congruence1}) is of course M3. Other congruences will be put together in subsection \ref{subsectionm3frombasiccongruences} to prove M4. We prove the above propositions in section (\ref{sectionofproofofcongruences}).

\subsection{$L$-values} \label{subsectionlvalues} Let $j \geq 0$. Let $x \in \Delta \times U_P^{ab}/Z^{p^j}$. Then we define $\delta^{(x)} : \Delta \times U_P^{ab} \rightarrow \mathbb{C}$ to be the characteristic function of the coset $x$ of $Z^{p^j}$ in $\Delta \times U_P^{ab}$. Define the partial zeta function by 
\[
\zeta(\delta^{(x)},s) = \sum_{\mathfrak{a}} \frac{\delta^{(x)}(g_{\mathfrak{a}})}{N(\mathfrak{a})^s}, \hspace{1cm} \text{for } Re(s) > 1,
\]
where the sum is over all ideals $\mathfrak{a}$ of $O_{F_P}$ which are prime to $\Sigma_P$, the Artin symbol of $\mathfrak{a}$ in $\Delta \times U_P^{ab}$ is denoted by $g_{\mathfrak{a}}$ and the absolute norm of the ideal $\mathfrak{a}$ is denoted by $N(\mathfrak{a})$. A well known theorem of Klingen \cite{Klingen:1962} and Seigel \cite{Seigel:1970} says that the function $\zeta(\delta^{(x)},s)$ has analytic continuation to the whole complex place except for a simple pole at $s=1$, and that $\zeta(\delta^{(x)}, 1-k)$ is rational for any even positive integer $k$. 

If $\epsilon$ is a locally constant function on $\Delta \times U_P^{ab}$ with values in a $\mathbb{Q}$ vector space $V$, say for a large enough $j$
\[
\epsilon \equiv \sum_{ x \in \Delta \times U_P^{ab}/Z^{p^j}} \epsilon(x) \delta^{(x)}.
\]
Then the special value $L_{\Sigma_P}(\epsilon, 1-k)$ can be canonically defined as 
\begin{equation}
L_{\Sigma_P}(\epsilon, 1-k) = \sum_{x \in \Delta \times U_P^{ab}/Z^{p^j}} \epsilon(x) \zeta(\delta^{(x)}, 1-k) \in V.
\label{lvalue}
\end{equation}
If $\epsilon$ is an Artin character of degree 1, then $L_{\Sigma_P}(\epsilon,1-k)$ is of course the value at $1-k$ of the complex $L$-function associated to $\epsilon$ with Euler factors at primes in $\Sigma_P$ removed. If $\epsilon$ is a locally constant $\mathbb{Q}_p$-values function on $\Delta \times U_P^{ab}$, then for any positive integer $k$ is divisible by $[F_{\infty}(\mu_p):F_{\infty}]$ and any $u \in U_P^{ab}$, we define 
\begin{equation}
\Delta_P^u(\epsilon, 1-k) = L_{\Sigma_P}(\epsilon,1-k) - \kappa(u)^kL_{\Sigma_P}(\epsilon_u,1-k),
\label{equationforDelta}
\end{equation}
where $\epsilon_u$ is a locally constant function defined by $\epsilon_u(g) = \epsilon(ug)$, for all $g \in \Delta \times U_P^{ab}$.

\subsection{Approximation to $p$-adic zeta functions} We get a sequence of elements in certain group rings which essentially approximate the abelian $p$-adic zeta functions $\zeta_P$. These group rings are obtained as follows. Recall that $\kappa$ is the $p$-adic cyclotomic character of $F$. Let $f$ be a positive integer such that $\kappa^{p-1}(Z) = 1+p^f\mathbb{Z}_p$.

\begin{definition} Let $P\leq P' \leq G$ be such that $P$ is normal in $P'$. Define $T_{P,P',j}$ be the image of the map
\[
tr_{P,P',j}: \mathbb{Z}_p[\Delta \times U_P^{ab}/Z^{p^j}]/(p^{f+j}) \rightarrow \mathbb{Z}_p[\Delta \times U_P^{ab}/Z^{p^j}]/(p^{f+j}),
\]
given by
\[
x \mapsto \sum_{g \in P'/P} gxg^{-1}.
\]
We denote $T_{P,N_GP,j}$ simply by $T_{P,j}$.
\end{definition}

\begin{lemma} We have an isomorphism
\[
\Lambda(\Delta \times U_P^{ab}) \xrightarrow{\sim} \ilim{j \geq 0} \mathbb{Z}_p[\Delta \times U_P^{ab}/Z^{p^j}]/(p^{f+j}).
\]
If $P\leq P'\leq G$ are such that $P$ is normal in $P'$, then under this isomorphism $T_{P,P'}$ maps isomorphically to $\ilim{j}T_{P,P',j}$. 
\label{inverselimit}
\end{lemma}

\noindent{\bf Proof:} We prove the surjectivity first. Given any 
\[
(x_j)_j \in \ilim{j \geq 0} \mathbb{Z}_p[\Delta \times U_P^{ab}/Z^{p^j}]/(p^{f+j}),
\]
 we construct a canonical $\tilde{x}_j \in \mathbb{Z}_p[\Delta \times U_P^{ab}/Z^{p^j}]$ as follows: for every $t \geq j$, let $\overline{x}_t$ be the image of $x_t \in \mathbb{Z}_p[\Delta \times U_P^{ab}/Z^{p^t}]/(p^{f+t})$ in $\mathbb{Z}_p[\Delta \times U_P^{ab}/Z^{p^j}]/(p^{f+t})$. Then $(\overline{x}_t)_{t \geq j}$ for an inverse system. We define $\tilde{x}_j$ to be the limit of $\overline{x}_t$ in $\mathbb{Z}_p[\Delta \times U_P^{ab}/Z^{p^j}]$. The tuple $(\tilde{x}_j)_{j \geq 0}$ forms an inverse system. We define $x$ to be their limit in $\Lambda(\Delta \times U_P^{ab})$. This is an inverse image of $(x_j)_{j \geq 0}$ in $\Lambda(\Delta \times U_P)$. This construction also proves the injectivity of the map. 
 
 To prove the second assertion we use the following exact sequence
 \[
 0 \rightarrow Ker(tr_{P,P',j}) \rightarrow \mathbb{Z}_p[\Delta \times U_P^{ab}/Z^{p^j}]/(p^{f+j}) \rightarrow T_{P,P',j} \rightarrow 0.
 \]
 Passing to the inverse limit over $j$ gives 
 \[
 0 \rightarrow \ilim{j} Ker(tr_{P,P',j}) \rightarrow \Lambda(\Delta \times U_P^{ab}) \rightarrow \ilim{j}T_{P,P',j} \rightarrow 0.
 \]
Exactness on the right is because all the abelian groups involved are finite. Hence $T_{P,P'} \cong \ilim{j} T_{P,P',j}$. 
\qed

\begin{proposition} (Ritter-Weiss) For any $j \geq 0$ and any positive integer $k$ divisible by $[F_{\infty}(\mu_p):F_{\infty}]$, the natural surjection of $\Lambda(\Delta \times U_P^{ab})$ onto $\mathbb{Z}_p[\Delta \times U_P^{ab}/Z^{p^j}]/(p^{f+j})$, maps $(1-u)\zeta_P \in \Lambda(\Delta \times U_P^{ab})$ to 
\[
\sum_{x \in U_P^{ab}/Z^{p^j}} \Delta_P^u(\delta^{(x)},1-k)\kappa(x)^{-k}x  \quad (\text{mod } p^{f+j}).
\]
In particular, we are claiming that the inverse limit is independent of the choice of $k$. Also note that since $x$ is a coset of $Z^{f+j}$ in $\Delta \times U_P^{ab}$, the value $\kappa(x)^k$ is well defined only module $p^{f+j}$.
\label{lemmaapprox}
\end{proposition}

\noindent{\bf Proof:} Since $\zeta_P$ is a pseudomeasure, $(1-u)\zeta_P$ lies in $\Lambda(\Delta \times U_P^{ab})$. We prove the proposition in 3 steps: first we show that the sums form an inverse system. Second we show that the inverse limit is independent of the choice of $k$. And thirdly we show that it interpolates the same values as $(1-u)\zeta_P$. 

{\bf Step1:} Let $j \geq 0$ be an integer. Let 
\[
\pi: \mathbb{Z}_p[\Delta \times U_P^{ab}/Z^{p^{j+1}}]/(p^{f+j+1}) \rightarrow \mathbb{Z}_p[\Delta \times U_P^{ab}/Z^{p^j}]/(p^{f+j}),
\]
denote the natural projection. Then 
\begin{align*}
& \pi\Big( \sum_{x \in \Delta \times U_P^{ab}/Z^{p^{j+1}}} \Delta_P^u(\delta^{(x)}, 1-k) \kappa(x)^{-k}x \Big) \\
& = \sum_{y \in \Delta \times U_P^{ab}/Z^{p^j}} \Big( \sum_{x \in yZ^{p^j}/Z^{p^{j+1}}} \Delta_P^u(\delta^{(x)}, 1-k) \kappa(x)^{-k} \pi(x)  \Big) (\text{mod }p^{f+j})\\
& = \sum_{y \in \Delta \times U_P^{ab}/Z^{p^j}} \Big( \kappa(y)^{-k}y \sum_{x \in Z^{p^j}/Z^{p^{j+1}}} \Delta_P^u(\delta^{(x)}, 1-k) \Big) (\text{mod }p^{f+j})\\
& = \sum_{y \in \Delta \times U_P^{ab}/Z^{p^j}} \Delta^u_P(\delta^{(y)}, 1-k) \kappa(y)^{-k}y (\text{mod }p^{f+j}).
\end{align*}
Here the second equality is because for any $x \in \Delta \times U_P^{ab}/Z^{p^{j+1}}$ we have $\kappa(x)^k \equiv \kappa(y)^k (\text{mod } p^{f+j})$ if $\pi(x) = y$. This shows that the sums form an inverse system.

{\bf Step2:} The inverse limit would be independent of the choice of $k$ if we show that for any two positive integers $k$ and $k'$ divisible by $[F_{\infty}(\mu_p):F_{\infty}]$, we have
\[
\sum_{x \in \Delta \times U_P^{ab}/Z^{p^j}} \Delta_P^u(\delta^{(x)}, 1-k) \kappa(x)^{-k}x \equiv \sum_{x \in \Delta \times U_P^{ab}/Z^{p^j}} \Delta_P^u(\delta^{(x)}, 1-k') \kappa(x)^{-k'}x (\text{mod } p^{f+j}).
\]
Or equivalently that,
\begin{equation}
\Delta_P^u(\delta^{(x)}, 1-k) \kappa(x)^{-k} \equiv \Delta_P^u(\delta^{(x)}, 1-k') \kappa(x)^{-k'} (\text{mod }p^{f+j}),
\label{deltacongruence}
\end{equation}
for all $x \in \Delta \times U_P^{ab}/Z^{p^j}$. Choose a locally constant function $\eta: \Delta \times U_P^{ab} \rightarrow \mathbb{Z}_p^{\times}$ such that $\eta^{[F_{\infty}(\mu_p):F_{\infty}]} \equiv \kappa^{[F_{\infty}(\mu_p):F_{\infty}]} (\text{mod }p^{f+j})$. Define the functions $\epsilon_k$ and $\epsilon_{k'}$ from $\Delta \times U_P^{ab}$ to $\mathbb{Q}_p$ by
\[
\epsilon_k = \frac{1}{p^{f+j}}\eta(x)^{-k}\delta^{(x)},
\]
and 
\[
\epsilon_{k'} = \frac{1}{p^{f+j}}\eta(x)^{-k'}\delta^{(x)}.
\]
Then the function $(\epsilon_k \kappa^{k-1} - \epsilon_{k'}\kappa^{k'-1})$ takes values in $\mathbb{Z}_p$. Hence the congruence (\ref{deltacongruence}) is satisfied by Deligne-Ribet (\cite{DeligneRibet:1980}, theorem 0.4).

{\bf Step3:} Let 
\[
\zeta_u = \ilim{j \geq 0} \Big( \sum_{x \in \Delta \times U_P^{ab}/Z^{p^j}} \Delta^u_P(\delta^{(x)}, 1-k)\kappa(x)^{-k}x (\text{mod } p^{f+j}) \Big) \in \Lambda(\Delta \times U_P^{ab}).
\]
Let $\epsilon$ be a locally constant function on $\Delta  \times U_P^{ab}$ factoring through $\Delta \times U_P^{ab}/Z^{p^j}$ for some $j \geq 0$. Note that for every $i \geq j$
\begin{align*}
& \sum_{x \in \Delta \times U_P^{ab}/Z^{p^i}} \Delta_P^u(\delta^{(x)}, 1-k) \epsilon(x) \\
& = \sum_{x \in \Delta \times U_P^{ab}/Z^{p^i}}L(\delta^{(x)}, 1-k) \epsilon(x) - \sum_{x \in \Delta \times U_P^{ab}/Z^{p^i}} \kappa(u)^k L(\delta^{(x)}_u, 1-k)\epsilon(x) \\
 & = \sum_{x \in \Delta \times U_P^{ab}/Z^{p^i}}L(\delta^{(x)}, 1-k) \epsilon(x) - \sum_{x \in \Delta \times U_P^{ab}/Z^{p^i}} \kappa(u)^k L(\delta^{(u^{-1}x)}, 1-k)\epsilon(x) \\
  & = \sum_{x \in \Delta \times U_P^{ab}/Z^{p^i}}L(\delta^{(x)}, 1-k) \epsilon(x) - \sum_{x \in \Delta \times U_P^{ab}/Z^{p^i}} \kappa(u)^k L(\delta^{(x)}, 1-k)\epsilon(ux) \\
& = L(\epsilon, 1-k) - \kappa(u)^k L(\epsilon_u, 1-k) \\
&= \Delta_P^u(\epsilon, 1-k).
\end{align*}
Then by definition of $\zeta_u$, for any $i \geq j$, we have
\begin{align*}
\zeta_u(\kappa^{k}\epsilon) & \equiv \sum_{ x \in \Delta \times U_P^{ab}/Z^{p^j}} \Delta_P^u(\delta^{(x)}, 1-k) \epsilon(x) (\text{mod } p^{f+i}) \\
& \equiv  \Delta_P^u(\epsilon, 1-k) (\text{mod } p^{f+i}).
\end{align*}

On the other hand, by definition of the $p$-adic zeta function or the construction of Deligne-Ribet (see discussion after theorem 0.5 in Deligne-Ribet \cite{DeligneRibet:1980}) we have
\[
(1-u)\zeta_P(\kappa^k\epsilon) = \Delta_P^u(\epsilon, 1-k).
\]
Hence $(1-u)\zeta_P = \zeta_u$ because they interpolate the same values on all on cyclotomic twists of locally constant functions. This finishes the proof of the proposition. 
\qed

\subsection{A sufficient condition to prove the basic congruences}

\label{sectionsuff}
\begin{lemma} Let $y$ be a coset of $Z^{p^j}$ in $\Delta \times U_P^{ab}$. Then for any $u\in Z$ and for any $g \in \mathcal{G}$, we have
\[
\Delta_P^u(\delta^{(y)}, 1-k) = \Delta_{gPg^{-1}}^u(\delta^{(gyg^{-1})}, 1-k).
\]
\label{invDelta}
\end{lemma}

\noindent{\bf Proof:} It is sufficient to show that $\zeta(\delta^{(y)}, 1-k) = \zeta(\delta^{(gyg^{-1})},1-k)$ because of the following:
\[
\Delta_P^u(\delta^{(y)},1-k) = \zeta(\delta^{(y)}, 1-k) - \kappa^k(u) \zeta(\delta^{(y)}_u, 1-k),
\]
\[
\Delta_{gPg^{-1}}^u(\delta^{(gyg^{-1})},1-k) = \zeta(\delta^{(gyg^{-1})}, 1-k) - \kappa^k(u) \zeta(\delta^{(gyg^{-1})}_u, 1-k).
\]
But
\[
\delta^{(y)}_u = \delta^{(u^{-1}y)} \qquad \text{and} \qquad \delta^{(gyg^{-1})}_u = \delta^{(u^{-1}gyg^{-1})} = \delta^{(gu^{-1}yg^{-1})}.
\]
Now to show that $\zeta(\delta^{(y)}, 1-k) = \zeta(\delta^{(gyg^{-1})},1-k)$, note that for $Re(s) >1$
\begin{align*}
\zeta(\delta^{(gyg^{-1})},s) & = \sum_{\mathfrak{a}} \frac{\delta^{(gyg^{-1})}(g_{\mathfrak{a}})}{N(\mathfrak{a})^s} \\
& = \sum_{\mathfrak{a}} \frac{\delta^{(y)}(g_{\mathfrak{a}^{g}})}{N(\mathfrak{a}^{g})^s} \\
& = \zeta(\delta^{(y)}, s).
\end{align*}
Since $\zeta(\delta^{(ygy^{-1})},s)$ and $\zeta(\delta^{(y)},s)$ are meromorphic functions agreeing on the right half plane, we deduce $\zeta(\delta^{(gyg^{-1})},1-k) = \zeta(\delta^{(y)},1-k)$, as required. 
\qed

\begin{proposition} To prove the congruence (\ref{congruence1}) in proposition (\ref{propcongruence1}) it is sufficient to prove the following: for any $j \geq 1$ and any coset $y$ of $Z^{p^j}$ in $\Delta \times U_P^{ab}$ fixed by $P'$ and any $u \in Z$ 
\begin{equation}
\Delta_P^{u^p}(\delta^{(y)},1-k) \equiv \Delta_{P'}^u(\delta^{(y)} \circ ver^{P'}_P,1-pk) (\text{mod } p\mathbb{Z}_p)
\label{equationdelta1}
\end{equation}
for all positive integers $k$ divisible by $[F_{\infty}(\mu_p):F_{\infty}]$.
\label{propsuffcongruence1}
\end{proposition}

\noindent{\bf Proof:} By lemma \ref{lemmaapprox} the image of $(1-u^p)\zeta_P$ in $\mathbb{Z}_p[\Delta \times U_P^{ab}/Z^{p^j}]/(p^{f+j-1})$ is 
\begin{equation}
\sum_{y \in \Delta \times U_P^{ab}/Z^{p^j}} \Delta_P^{u^p}(\delta^{(y)},1-k) \kappa(y)^{-k}y (\text{mod } p^{f+j-1}).
\label{equationzeta1}
\end{equation}
And the image of $(1-u)\zeta_{P'}$ in $\mathbb{Z}_p[\Delta \times U_{P'}^{ab}/Z^{p^{j-1}}]/(p^{f+j-1})$ is
\[
\sum_{x \in \Delta \times U_{P'}^{ab}/Z^{p^{j-1}}} \Delta_{P'}^u(\delta^{(x)},1-pk) \kappa(x)^{-pk}x (\text{mod } p^{f+j-1}).
\]
Let $V_{P,P'}$ be the kernel of the homomorphism $ver^{P'}_P:U_{P'}^{ab} \rightarrow U_P^{ab}$. Then $V_{P,P'} \cap Z = \{1\}$ which implies that the map 
\[
\Delta \times U_{P'}^{ab}/V_{P,P'}Z^{p^{j-1}} \rightarrow \Delta \times U_P^{ab}/Z^{p^j}
\]
induced by $ver^{P'}_P$ is injective. Moreover $\kappa(V_{P,P'})^k = \{1\}$. Hence the image of $ver^{P'}_P((1-u)\zeta_{P'}) = (1-u^p)ver^{P'}_P(\zeta_{P'})$ in $\mathbb{Z}_p[\Delta \times U_P^{ab}/Z^{p^j}]/(p^{f+j-1})$ is 
\[
\sum_{x \in \Delta \times U_{P'}^{ab}/V_{P',P}Z^{p^{j-1}}} \Delta_{P'}^u(\delta^{(x)},1-pk)\kappa(x)^{-pk}ver^{P'}_P(x) (\text{mod }p^{f+j-1}),
\]
which can be written as
\begin{equation}
\sum_{y \in \Delta \times U_P^{ab}/Z^{p^j}} \Delta_{P'}^u(\delta^{(y)}\circ ver^{P'}_P,1-pk) \kappa(y)^{-k}y (\text{mod }p^{f+j-1})
\label{equationzeta2}
\end{equation}
because if $y \notin Im(ver^{P'}_P)$, then $\delta^{(y)} \circ ver^{P'}_P \equiv 0$ and if $y = ver^{P'}_P(x)$, then $\kappa(y)^k = \kappa(x)^{pk}$. Subtracting (\ref{equationzeta2}) from (\ref{equationzeta1}) gives
\begin{equation}
\sum_{y \in \Delta \times U_P^{ab}/Z^{p^j}} \Big(\Delta_P^{u^p}(\delta^{(y)}\circ ver^{P'}_P,1-k) - \Delta_{P'}^u(\delta^{(y)},1-pk) \Big)\kappa(y)^{-k} y (\text{mod } p^{f+j-1}).
\label{sum}
\end{equation}
If $y$ is fixed by $P'$ then $\Big(\Delta_P^{u^p}(\delta^{(y)},1-k) - \Delta_{P'}^u(\delta^{(y)},1-pk) \Big)\kappa(y)^{-k} y \equiv py \equiv 0 (\text{mod }T_{P,P',j})$ under equation (\ref{equationdelta1}). On the other hand if $y$ is not fixed by $P'$, then the full $y$ orbit of under action of $P'$ in the above sum is
\begin{align*}
& \sum_{g \in P'/P} (\Delta_P^{u^p}(\delta^{(gyg^{-1})},1-k) - \Delta_{P'}^u(\delta^{(y)}\circ ver^{P'}_P,1-pk))\kappa(gyg^{-1})^{-k}gyg^{-1} \\
= & \Big(\Delta_P^{u^p}(\delta^{(y)},1-k)-\Delta_{P'}^u(\delta^{(y)}\circ ver^{P'}_P,1-pk)\Big) \kappa(y)^{-k} \sum_{g \in P'/P} gyg^{-1} \\
\in & T_{P,P',j}.
\end{align*}
The second equality is by lemma \ref{invDelta}. Hence the sum in (\ref{sum}) lies in $T_{P,P',j} (\text{mod }p^{f+j-1})$. By lemma \ref{inverselimit} $(1-u^p)(\zeta_P-ver^{P'}_P(\zeta_{P'})) \in T_{P,P'}$. By as $u$ is a central element congruence (\ref{congruence1}) holds. \qed 

\begin{remark} The proof of following three proposition is very similar to the above proof.
\end{remark}

\begin{proposition} To prove congruence (\ref{congruence2}) in proposition (\ref{propcongruence2}) it is sufficient to show the following: for any $j \geq 0$ and any coset $y$ of $Z^{p^j}$ in $\Delta \times U_P$ whose image in $P$ is a generator of $P$, and any $u \in Z$
\begin{equation}
\Delta_P^{u^{d_P}}(\delta^{(y)},1-k) \equiv 0 (\text{mod } |(W_GP)_y| \mathbb{Z}_p),
\label{equationdelta2}
\end{equation}
for all positive integers $k$ divisible by $[F_{\infty}(\mu_p):F_{\infty}]$.
\label{propsuffcongruence2}
\end{proposition}

\noindent{\bf Proof:} Let $v=u^{d_P}$. Then by lemma \ref{lemmaapprox} the image of $(1-v)\zeta_P - \omega_P^k((1-v)\zeta_P)$ in $\mathbb{Z}_p[\Delta \times U_P/Z^{p^j}]/(p^{f+j})$ is 
\begin{equation}
\sum_{y \in \Delta \times U_P/Z^{p^j}} \Delta_P^v(\delta^{(y)},1-k) \kappa(y)^{-k}(y-\omega_P^k(y)) (\text{mod } p^{f+j}).
\label{equationzeta3}
\end{equation}
In the image of $y$ in $P$ is not a generator of $P$, then $y-\omega_P^k(y)=0$. For $y$ whose image in $P$ is a generator of $P$, we look at the $W_GP$ orbit of $y$ in expression (\ref{equationzeta3}). It is 
\begin{align*}
&\sum_{g \in W_GP/(W_GP)_y}  \Delta_P^v(\delta^{(gyg^{-1})},1-p) \kappa(gyg^{-1})^{-k}(gyg^{-1} - \omega_P^k(gyg^{-1}) (\text{mod }p^{f+j}) \\
=& \Delta_P^v(\delta^{(y)},1-k)\kappa(y)^{-k} \sum_{g \in W_GP/(W_GP)_y} (gyg^{-1}- \omega_P^k(gyg^{-1})) 
\end{align*}
which lies in $\mathfrak{p} T_{P,j}$ under equation (\ref{equationdelta2}) and then the sum in expression (\ref{equationzeta3}) lies in $\mathfrak{p}T_{P,j}$. Then by lemma \ref{inverselimit} $(1-v)(\zeta_P - \omega_P^k(\zeta_P)) \in \mathfrak{p}T_P$. As $v$ is a central element congruence (\ref{congruence2}) holds. \qed

\begin{proposition} To prove congruence (\ref{congruence3}) in proposition (\ref{propcongruence3}) it is sufficient to prove the following: for any $j \geq 0$ and any coset $y$ of $Z^{p^j}$ in $\Delta \times U_P$, and any $u \in Z$
\begin{equation}
\Delta_P^{u^{d_P}}(\delta^{(y)}, 1-k) \equiv \sum_{P'} \Delta_{P'}^{u^{d_P/p}}(\delta^{(y)} \circ \varphi, 1-pk) (\text{mod } |(W_{G}P)_y| \mathbb{Z}_p),
\label{equationdelta3}
\end{equation}
for all positive integers $k$ divisible by $[F_{\infty}(\mu_p):F_{\infty}]$. The sum in the congruence runs over all $P' \in C(G)$ such that $P'^p = P$ and $P' \neq P$. 
\label{propsuffcongruence3}
\end{proposition}

\noindent{\bf Proof:} Let $v= u^{d_P/p}$. By lemma \ref{lemmaapprox} the image of $(1-v^p)\zeta_P$ in $\mathbb{Z}_p[\Delta \times U_P/Z^{p^j}]/(p^{f+j-1})$ is 
\begin{equation}
\sum_{y \in \Delta \times U_P/Z^{p^j}} \Delta_P^{v^p}(\delta^{(y)},1-k) \kappa(y)^{-k}y (\text{mod } p^{f+j-1}).
\label{equationzeta4}
\end{equation}
And the image of $(1-v)\zeta_{P'}$ in $\mathbb{Z}_p[\Delta \times U_{P'}/Z^{p^{j-1}}]/(p^{f+j-1})$ is
\[
\sum_{x \in \Delta \times U_{P'}/Z^{p^{j-1}}} \Delta_{P'}^v(\delta^{(x)},1-pk) \kappa(x)^{-pk}x (\text{mod } p^{f+j-1}).
\]
Let $V_{P,P'}$ be the kernel of the homomorphism $\varphi:U_{P'} \rightarrow U_P$. Then $V_{P,P'} \cap Z = \{1\}$ which implies that the map 
\[
\Delta \times U_{P'}/V_{P,P'}Z^{p^{j-1}} \rightarrow \Delta \times U_P/Z^{p^j}
\]
induced by $\varphi$ is injective. Moreover, $\kappa(V_{P,P'})^k=\{1\}$. Hence the image of $\sum_{P'}\varphi((1-v)\zeta_{P'}) = (1-v^p)\sum_{P'}\varphi(\zeta_{P'})$ in $\mathbb{Z}_p[\Delta \times U_P/Z^{p^j}]/(p^{f+j-1})$ is 
\[
\sum_{P'}\sum_{x \in \Delta \times U_{P'}/V_{P',P}Z^{p^{j-1}}} \Delta_{P'}^v(\delta^{(x)},1-pk)\kappa(x)^{-pk}\varphi(x) (\text{mod }p^{f+j-1}),
\]
which can be written as
\begin{equation}
\sum_{y \in \Delta \times U_P/Z^{p^j}} \sum_{P'} \Delta_{P'}^v(\delta^{(y)}\circ \varphi,1-pk) \kappa(y)^{-k}y (\text{mod }p^{f+j-1})
\label{equationzeta5}
\end{equation}
because if $y \notin Im(\varphi)$, then $\delta^{(y)} \circ \varphi \equiv 0$ and if $y = \varphi(x)$, then $\kappa(y)^k = \kappa(x)^{pk}$. Subtracting (\ref{equationzeta5}) from (\ref{equationzeta4}) gives 
\begin{equation} 
\sum_{y \in \Delta \times U_P/Z^{p^j}} \Big(\Delta_P^{v^p}(\delta^{(y)},1-k) - \sum_{P'} \Delta_{P'}^v(\delta^{(y)} \circ \varphi, 1-pk) \Big) \kappa(y)^{-k} y (\text{mod } p^{f+j-1})
\label{equationzeta6}
\end{equation}
Now we take the orbit of $y$ in the sum in (\ref{equationzeta6}) under the action of $W_GP$. It is 
\begin{align*} 
\Big(\Delta_P^{v^p}(\delta^{(y)},1-k)-\sum_{P'}\Delta_{P'}^v(\delta^{(y)}\circ \varphi,1-pk)\Big) \kappa(y)^{-k} \sum_{g \in W_GP/(W_GP)_y} gyg^{-1} 
\end{align*}
which lies in $T_{P,j} (\text{mod }p^{f+j-1})$ under equation (\ref{equationdelta3}) and then the sum in expression (\ref{equationzeta6}) lies in $T_{P,j}(\text{mod }p^{f+j-1})$. Then by lemma (\ref{inverselimit}) $(1-v^p)(\zeta_P- \sum_{P'}\varphi(\zeta_{P'})) \in T_P$. As $v$ is a central element congruence (\ref{congruence3}) holds. \qed

\begin{proposition} To prove congruence (\ref{congruence4}) in proposition (\ref{propcongruence4}) it is sufficient to prove the following: for any $j \geq 0$, any coset $y$ of $Z^{p^j}$ in $\Delta \times U_P$ and any $u$ in $Z$ 
\begin{equation}
\Delta_P^{u^{pd_{P'}}}(\delta^{(y)},1-k) \equiv \Delta_{P'}^{u^{d_{P'}}}(\delta^{(y)} \circ \varphi, 1-pk) (\text{mod } |(N_{G}P'/P)_y| \mathbb{Z}_p),
\label{equationdelta4}
\end{equation}
for all positive integers $k$ divisible by $[F_{\infty}(\mu_p):F_{\infty}]$.

To prove the congruence (\ref{congruence5}) in proposition (\ref{propcongruence4}) it is sufficient to show the following: for any $j \geq 1$, any coset $y$ of $Z^{p^j}$ in $\Delta \times Z^p$ and any $u$  in $Z^p$
\begin{equation}
\Delta_0^{u^{p|G|}}(\delta^{(y)}, 1-k) \equiv \Delta_{\{1\}}^{u^{|G|}}(\delta^{(y)} \circ \varphi, 1-pk) (\text{mod } p|G| \mathbb{Z}_p),
\label{equationdelta5}
\end{equation}
for any positive integer $k$ divisible by $[F_{\infty}(\mu_p):F_{\infty}]$.
\label{propsuffcongruence4}
\end{proposition}

\noindent{\bf Proof:} We will only prove the first assertion. Proof of the second one exactly the same. Let $v= u^{d_{P'}}$. By lemma \ref{lemmaapprox} the image of $(1-v^p)\zeta_P$ in $\mathbb{Z}_p[\Delta \times U_P/Z^{p^j}]/(p^{f+j-1})$ is 
\begin{equation}
\sum_{y \in \Delta \times U_P/Z^{p^j}} \Delta_P^{v^p}(\delta^{(y)},1-k) \kappa(y)^{-k}y (\text{mod } p^{f+j-1}).
\label{equationzeta7}
\end{equation}
And the image of $(1-v)\zeta_{P'}$ in $\mathbb{Z}_p[\Delta \times U_{P'}/Z^{p^{j-1}}]/(p^{f+j-1})$ is
\[
\sum_{x \in \Delta \times U_{P'}/Z^{p^{j-1}}} \Delta_{P'}^v(\delta^{(x)},1-pk) \kappa(x)^{-pk}x (\text{mod } p^{f+j-1}).
\]
Let $V_{P,P'}$ be the kernel of the homomorphism $\varphi:U_{P'} \rightarrow U_P$. Then $V_{P,P'} \cap Z = \{1\}$ which implies that the map 
\[
\Delta \times U_{P'}/V_{P,P'}Z^{p^{j-1}} \rightarrow \Delta \times U_P/Z^{p^j}
\]
induced by $\varphi$ is injective. Moreover $\kappa(V_{P,P'})^k=\{1\}$. Hence the image of $\varphi((1-v)\zeta_{P'}) = (1-v^p)\varphi(\zeta_{P'})$ in $\mathbb{Z}_p[\Delta \times U_P/Z^{p^j}]/(p^{f+j-1})$ is 
\[
\sum_{x \in \Delta \times U_{P'}/V_{P',P}Z^{p^{j-1}}} \Delta_{P'}^v(\delta^{(x)},1-pk)\kappa(x)^{-pk}\varphi(x) (\text{mod }p^{f+j-1}),
\]
which can be written as
\begin{equation}
\sum_{y \in \Delta \times U_P/Z^{p^j}} \Delta_{P'}^v(\delta^{(y)}\circ \varphi,1-pk) \kappa(y)^{-k}y (\text{mod }p^{f+j-1})
\label{equationzeta8}
\end{equation}
because if $y \notin Im(\varphi)$, then $\delta^{(y)} \circ \varphi \equiv 0$ and if $y = \varphi(x)$, then $\kappa(y)^k = \kappa(x)^{pk}$. Subtracting (\ref{equationzeta8}) from (\ref{equationzeta7}) we get
\begin{equation}
\sum_{y \in \Delta \times U_P/Z^{p^j}} \Big( \Delta_P^{v^p}(\delta^{(y)},1-k) - \Delta_{P'}^v(\delta^{(y)}\circ \varphi,1-pk) \Big) \kappa(y)^{-k} y (\text{mod }p^{f+j-1}).
\label{equationzeta9}
\end{equation}
Now for a fixed $y \in \Delta \times U_P/Z^{p^j}$ we take the orbit of $y$ in this sum under the action of $N= N_GP/P'$. It is
\[
\Big( \Delta_P^{v^p}(\delta^{(y)},1-k) - \Delta_{P'}^v (\delta^{(y)} \circ \varphi,1-pk) \Big) \kappa(y)^{-k} \sum_{g \in N/N_y} gyg^{-1}
\]
which lies in $T_{P,N_GP',j} (\text{mod } p^{f+j-1})$ under equation (\ref{equationdelta4}) and then the sum in (\ref{equationzeta9}) lies in $T_{P,N_GP',j} (\text{mod } p^{f+j-1})$. Then by lemma (\ref{inverselimit}) $(1-v^p)(\zeta_P - \varphi(\zeta_{P'})) \in T_{P,N_GP'}$. As $v$ is a central element congruence (\ref{congruence4}) holds. \qed

\subsection{Hilbert modular forms} \label{subsectionhmf} In this section we briefly recall the basic notions in the theory of Hilbert modular forms. Let $L$ be an arbitrary totally real number field of degree $r$ over $\mathbb{Q}$. Let $\mathfrak{H}_L$ be the Hilbert upper half plane of $L$. Let $\Sigma$ be a finite set of finite primes of $L$ containing all primes above $p$. Let $\kappa$ be the $p$-adic cyclotomic character of $L$. Let $\mathfrak{f}$ be an integral ideal of $L$ with all its prime factors in $\Sigma$. We put $GL^+_2(L \otimes \mathbb{R})$ for the group of all $2 \times 2$ matrices with totally positive determinant. For any even positive integer $k$, the group $GL^+_2(L \otimes \mathbb{R})$ acts on functions $f: \mathfrak{H}_L \rightarrow \mathbb{C}$ by
\[
f|k \left( \begin{array}{ll} a & b \\ c & d \end{array} \right)(\tau) = \mathcal{N}(ad-bc)^{k/2}\mathcal{N}(c\tau+d)^{-k}f(\frac{a\tau+d}{c\tau+d}),
\]
where $\mathcal{N} : L \otimes \mathbb{C} \rightarrow \mathbb{C}$ is the norm map. Set
\[
\Gamma_{00}(\mathfrak{f}) = \{ \left( \begin{array}{ll} a & b \\ c & d \end{array}\right) \in SL_2(L): a,d \in 1+\mathfrak{f}, b \in \mathfrak{D}^{-1}, c \in \mathfrak{fD} \}, 
\]
where $\mathfrak{D}$ is the different of $L/\mathbb{Q}$. A \emph{Hilbert modular form} $f$ of weight $k$ on $\Gamma_{00}(\mathfrak{f})$ is a holomorphic function $f: \mathfrak{H}_L \rightarrow \mathbb{C}$ (which we assume to be holomorphic at $\infty$ if $L=\mathbb{Q}$) satisfying 
\[
f|_kM = f \qquad \text{for all } M \in \Gamma_{00}(\mathfrak{f}).
\]
The space of all Hilbert modular forms of weight $k$ on $\Gamma_{00}(\mathfrak{f})$ is denoted by $M_k(\Gamma_{00}(\mathfrak{f}), \mathbb{C})$. Since $f$ is invariant under the translation $\tau \mapsto \tau + b$ (for $b \in \mathfrak{D}^{-1}$), we may expand $f$ as a Fourier series to get the standard $q$-expansion
\[
f(\tau)= c(0,f) + \sum_{\mu} c(\mu,f)q^{\mu}_L,
\]
where $\mu$ runs through all totally positive elements in $O_L$ and $q^{\mu}_L = e^{2\pi i tr_{L/\mathbb{Q}}(\mu \tau)}$. 

\subsection{Restrictions along diagonal} Let $L'$ be another totally real number field containing $L$. Let $r'$ be the degree of $L'$ over $L$. The inclusion of $L$ in $L'$ induces  maps $\mathfrak{H}_L \xrightarrow{*} \mathfrak{H}_{L'}$ and $SL_2(L \otimes \mathbb{R}) \xrightarrow{*} SL_2(L' \otimes\mathbb{R})$. For a holomorphic function $f : \mathfrak{H}_{L'} \rightarrow \mathbb{C}$, we define the ``restriction along diagonal" $R_{L'/L}f : \mathfrak{H}_L \rightarrow \mathbb{C}$ by $R_{L'/L}f(\tau) = f(\tau^*)$. We then have
\[
(R_{L'/L}f)|_{r'k}M = R_{L'/L}(f|_kM^*),
\]
for any $M \in SL_2(L \otimes \mathbb{R})$. Let $\mathfrak{f}$ be an integral ideal of $L$, then $R_{L'/L}$ induces a map
\[
R_{L'/L} : M_k(\Gamma_{00}(\mathfrak{f}O_{L'}), \mathbb{C}) \rightarrow M_{r'k}(\Gamma_{00}(\mathfrak{f}), \mathbb{C}).
\]
If the standard $q$-expansion of $f$ is 
\[
c(0,f) + \sum_{\nu \in O_{L'}^+} c(\nu,f) q^{\nu}_{L'},
\]
then the standard $q$-expansion of $R_{L'/L}f$ is
\[
c(0,f) + \sum_{\mu \in O_L^+} \Big(\sum_{\nu : tr_{L'/L}(\nu) = \mu}c(\nu,f)\Big) q^{\mu}_L.
\]
Here $O_L^+$ and $O_{L'}^+$ denotes totally positive elements of $O_L$ and $O_{L'}$ respectively.

\subsection{Cusps} Let $\mathbb{A}_L$ be the ring of finite adeles of $L$. Then by strong approximation
\[
SL_2(\mathbb{A}_L) = \widehat{\Gamma_{00}(\mathfrak{f})}\cdot SL_2(L).
\]
Any $M \in SL_2(\mathbb{A}_L)$ can be written as $M_1M_2$ with $M_1 \in \widehat{\Gamma_{00}(\mathfrak{f})}$ and $M_2 \in SL_2(L)$. We define $f|_kM$ to be $f|_kM_2$. Any $\alpha$ in $\mathbb{A}_L^{\times}$ determines a \emph{cusp}. We let 
\[
f|_{\alpha} = f|_k \left( \begin{array}{ll} \alpha & 0 \\ 0 & \alpha^{-1} \end{array} \right).
\]
\emph{The $q$-expansion of $f$ at the cusp determined by $\alpha$} is defined to the standard $q$-expansion of $f|_{\alpha}$. We write it as
\[
c(0,\alpha,f) + \sum_{\mu} c(\mu, \alpha,f) q_L^{\mu},
\]
where the sum is restricted to all totally positive elements of $L$ which lie in the square of the ideal ``generated" by $\alpha$. 

\begin{lemma} Let $\mathfrak{f}$ be an integral ideal in $L$. Let $f \in M_k(\Gamma_{00}(\mathfrak{f}O_{L'}), \mathbb{C})$. Then the constant term of the $q$-expansion of $R_{L'/L}f$ at the cusp determined by $\alpha \in \mathbb{A}_L^{\times}$ is equal to the constant term of the $q$-expansion of $f$ at the cusp determined by $\alpha^* \in \mathbb{A}_{L'}^{\times}$ i.e. 
\[
c(0,\alpha,R_{L'/L}f) = c(0,\alpha^*,f).
\]
\label{constanttermofrestriction}
\end{lemma}

\noindent{\bf Proof:} The $q$-expansion of $f$ at the cusp determined by $\alpha^*$ is the standard $q$-expansion of $f|_{\alpha^*}$. Similarly, the $q$-expansion of $R_{L'/L}f$ at the cusp determined by $\alpha$ is the standard $q$-expansion of $(R_{L'/L}f)|_{\alpha}$. But $(R_{L'/L}f)|_{\alpha} = R_{L'/L}(f|_{\alpha^*})$. 
\qed

\subsection{A Hecke operator} 
\begin{lemma} Let $\beta \in O_L$ be a totally positive element. Assume that $\mathfrak{f} \subset \beta O_L$. Then there is a Hecke operator $U_{\beta}$ on $M_k(\Gamma_{00}(\mathfrak{f}), \mathbb{C})$ so that for $f \in M_k(\Gamma_{00}(\mathfrak{f}), \mathbb{C})$ has the standard $q$-expansion of $f|_kU_{\beta}$ is
\[
c(0,f) + \sum_{\mu} c(\mu \beta,f) q^{\mu}_L.
\]
\end{lemma}

\noindent{\bf Proof:} The claimed operator $U_{\beta}$ is the one defined by $ \left( \begin{array}{l l} \beta & 0 \\ 0 & 1 \end{array} \right)$. Then
\[
\Gamma_{00}(\mathfrak{f})  \left( \begin{array}{l l} \beta & 0 \\ 0 & 1 \end{array} \right) \Gamma_{00}(\mathfrak{f}) = \cup_{b} \Gamma_{00}(\mathfrak{f}) \left( \begin{array}{l l} 1 & b \\ 0 & \beta \end{array} \right).
\]
where $b$ ranges over all coset representatives of $\beta \mathfrak{D}$  in $\mathfrak{D}$ and the union is a disjoint union. Define
\[
f|_kU_{\beta}(\tau)  = \mathcal{N}(\beta)^{k/2-1} \sum_{b} f|_k\left( \begin{array}{l l} 1 & b \\ 0 & \beta \end{array} \right)(\tau),
\]
where $b$ runs through the set of coset representatives of $\beta \mathfrak{D}$ in $\mathfrak{D}$. Then

\begin{align*}
f|_kU_{\beta}(\tau) & = \mathcal{N}(\beta)^{k/2-1} \sum_{b} f|_k\left( \begin{array}{l l} 1 & b \\ 0 & \beta \end{array} \right)(\tau) \\
& = \mathcal{N}(\beta)^{k/2-1}\mathcal{N}(\beta)^{k/2} \mathcal{N}(\beta)^{-k}\sum_b f(\frac{\tau + b }{\beta}) \\
& = \mathcal{N}(\beta)^{-1} \sum_b \Big(c(0,f) + \sum_{\mu} c(\mu,f) e^{2 \pi i tr_{L/\mathbb{Q}}(\mu(\beta^{-1}\tau + \beta^{-1}b))} \Big) \\
& = c(0,f) + \mathcal{N}(\beta)^{-1} \sum_{\mu} c(\mu,f) e^{2 \pi i tr_{L/\mathbb{Q}}(\mu \tau / \beta)} ( \sum_b e^{2 \pi i tr_{L/\mathbb{Q}}(\mu b/ \beta)} )
\end{align*}
The sum $\sum_{b} e^{2 \pi i tr_{L/\mathbb{Q}}(\mu b/ \beta)}= 0$ unless $\mu \in \beta O_L$. On the other hand, if $\mu \in \beta O_L$, then $\sum_b e^{2 \pi i tr_{L/\mathbb{Q}}(\mu b/ \beta)} = \mathcal{N}(\beta)$. Hence we get 
\[
f|_kU_{\beta} (\tau) = c(0,f) + \sum_{\mu} c(\mu \beta,f) q_L^{\mu}.
\]
\qed

\subsection{Eisenstein series} The following proposition is proven by Deligne-Ribet (\cite{DeligneRibet:1980}, proposition 6.1).

\begin{proposition} Let $L_{\Sigma}$ be the maximal abelian totally real extension of $L$ unramified outside $\Sigma$. Let $\epsilon$ be a locally constant $\mathbb{C}$-valued function on $Gal(L_{\Sigma}/L)$. Then for every even positive integer $k$ \\
(i) There is an integral ideal $\mathfrak{f}$ of $L$ with all its prime factors in $\Sigma$, and a Hilbert modular form $G_{k,\epsilon}$ in $M_k(\Gamma_{00}(\mathfrak{f}), \mathbb{C})$ with standard $q$-expansion 
\[
2^{-r}L(\epsilon, 1-k)+ \sum_{\mu}\Big( \sum_{\mathfrak{a}} \epsilon(g_{\mathfrak{a}}) N(\mathfrak{a})^{k-1}\Big) q^{\mu}_L,
\]
where the first sum ranges over all totally positive $\mu \in O_L$, and the second sum ranges over all integral ideals $\mathfrak{a}$ of $L$ containing $\mu$ and prime to $\Sigma$. Here $g_{\mathfrak{a}}$ is the image of $\mathfrak{a}$ under the Artin symbol map. $N(\mathfrak{a})$ denotes norm of the ideal $\mathfrak{a}$. \\
(ii) Let $q$-expansion of $G_{k, \epsilon}$ at the cusp determined by any $\alpha \in \mathbb{A}_L^{\times}$ has constant term 
\[
N^k((\alpha))2^{-r}L(\epsilon_g, 1-k),
\]
where $(\alpha)$ is the ideal of $L$ generated by $\alpha$ and $N((\alpha))$ is its norm. The element $g$ is the image of $(\alpha)$ under the Artin symbol map (see for instance 2.22 in Deligne-Ribet \cite{DeligneRibet:1980}). The locally constant function $\epsilon_g$ is given by
\[
\epsilon_g(h) = \epsilon(gh) \hspace{1cm} \text{for all $h \in Gal(L_{\Sigma}/L)$}.
\]
\label{propeseries}
\end{proposition}

\subsection{The $q$-expansion principle} \label{subsectionqexpansionprinciple}
Let $f \in M_k(\Gamma_{00}(\mathfrak{f}), \mathbb{Q}))$ i.e. $c(\mu,\alpha,f) \in \mathbb{Q}$ for all $\mu \in O_L^+ \cup \{0\}$ and all $\alpha \in \mathbb{A}_L^{\times}$. Suppose the standard $q$-expansion of $f$ has all non-constant coefficients in $\mathbb{Z}_{(p)}$, Let $\alpha \in \mathbb{A}_{L}^{\times}$ be a finite adele. Then 
\[
c(0,f) - N(\alpha_p)^{-k}c(0,\alpha,f) \in \mathbb{Z}_p.
\]
Here $\alpha_p \in L \otimes_{\mathbb{Q}}\mathbb{Q}_p$ is the $p$th component of $\alpha$ and $N: L \otimes_{\mathbb{Q}} \mathbb{Q}_p \rightarrow \mathbb{Q}_p$ is the norm map. This is the $q$-expansion principle of Deligne-Ribet (see \cite{DeligneRibet:1980} 0.3 and 5.13-5.15). 

\begin{remark} Hence if $u$ is the image in $Gal(L_{\Sigma}/L)$ of an id\`{e}le $\alpha$ under the Artin symbol map, then using the equation $N((\alpha))^kN(\alpha_p)^{-k} = \kappa(u)^k$, we get
\[
c(0,G_{k,\epsilon}) - N(\alpha_p)^{-k} c(0,\alpha,G_{k,\epsilon}) = 2^{-r}\Delta^{u}(\epsilon, 1-k),
\]
for any positive even integer $k$.
\label{differenceofconstantterms}
\end{remark}

\subsection{Proof of the sufficient conditions in section \ref{sectionsuff}}

\begin{proposition} The sufficient condition in proposition (\ref{propsuffcongruence1}) for proving proposition (\ref{propcongruence1}) holds. Hence M3 holds.
\end{proposition}

\noindent{\bf Proof:} We must show that for any $P \leq P' \leq G_p$ such that $[P':P]=p$ and any $j \geq 0$, any coset $y$ of $Z^{p^j}$ in $\Delta \times U_P^{ab}$ fixed by $P'$ and any $u$ in $Z$, we have the congruence 
\[
\Delta_P^{u^p}(\delta^{(y)},1-k) \equiv \Delta_{P'}^u(\delta^{(y)} \circ ver^{P'}_P, 1-pk) (\text{mod } p\mathbb{Z}_p),
\]
for all positive integers $k$ divisible by $[F_{\infty}(\mu_p):F_{\infty}]$. Choose an integral ideal $\mathfrak{f}$ of $F_{P'}$ such that the Hilbert Eisenstein series $G_{k,\delta^{(y)}}$ and $G_{pk, \delta^{(y)} \circ ver^{P'}_P}$, given by proposition (\ref{propeseries}), on $\mathfrak{H}_P$ and $\mathfrak{H}_{F_{P'}}$ respectively are defined over $\Gamma_{00}(\mathfrak{f}O_{F_P})$ and $\Gamma_{00}(\mathfrak{f})$ respectively. Moreover, we may assume that all prime ideals dividing $\mathfrak{f}$ lie in $\Sigma_{F_{P'}}$ and $\mathfrak{f} \subset pO_{F_{P'}}$. Define $E$ by 
\[
E=R_{F_P/F_{P'}}(G_{k,\delta^{(y)}})|_{pk}U_p - G_{pk, \delta^{(y)} \circ ver^{P'}_P}.
\]
Let $\alpha \in \mathbb{A}_{F_{P'}}^{\times}$ whose image in $\Delta \times U_{P'}^{ab}$ under the Artin symbol map coincides with $u$. Then by lemma \ref{constanttermofrestriction} and remark \ref{differenceofconstantterms}
\begin{align*}
&c(0,E) - N(\alpha_p)^{-pk}c(0,\alpha,E) \\ & = 2^{-r_P} \Delta_P^{u^p}(\delta^{(y)},1-k) - 2^{-r_{P'}} \Delta_{P'}^u(\delta^{(y)} \circ ver^{P'}_P,1-pk).
\end{align*}
Note that the image of $\alpha^*$ in $\Delta \times U_P^{ab}$ under the Artin symbol map is $u^p$. Since $2^{-r_P} \equiv 2^{-r_{P'}} (\text{mod }p)$ it is enough to prove, using the $q$-expansion principle that the non-constant terms of the standard $q$ expansion of $E$ all lie in $p \mathbb{Z}_{(p)}$ i.e. for all $\mu \in O_{F_{P'}}^+$
\begin{align*}
& c(\mu,E) = c(p\mu, R_{F_P/F_{P'}}(G_{k.\delta^{(y)}})) - c(\mu,G_{pk,\delta^{(y)}\circ ver^{P'}_P}) \\
&=\sum_{(\mathfrak{b},\eta)} \delta^{(y)}(g_{\mathfrak{b}}) N(\mathfrak{b})^{k-1} - \sum_{\mathfrak{a}} \delta^{(y)}(g_{\mathfrak{a}O_{F_P}})N(\mathfrak{a})^{pk-1} \in p\mathbb{Z}_{(p)}
\end{align*}
Here the pairs $(\mathfrak{b},\eta)$ runs through all integral ideals $\mathfrak{b}$ of $F_P$ which are prime to $\Sigma_{F_P}$ and contains the totally positive element $\eta \in O_{F_P}$ and $tr_{F_P/F_{P'}}(\eta) = p\mu$. The ideal $\mathfrak{a}$ runs through all integral ideals of $F_{P'}$ prime to $\Sigma_{F_{P'}}$ and contains $\mu$. The group $P'/P$ acts trivially on the pair $(\mathfrak{b},\eta)$ if and only if there is an $\mathfrak{a}$ such that $\mathfrak{a}O_{F_P} = \mathfrak{b}$ and $\eta \in O_{F_{P'}}$. In this case
\begin{align*}
& \delta^{(y)}(g_{\mathfrak{b}})N(\mathfrak{b})^{k-1} - \delta^{(y)}(g_{\mathfrak{a}O_{F_P}}) N(\mathfrak{a})^{pk-1} \\
 = & \delta^{(y)}(g_{\mathfrak{b}})(N(\mathfrak{a})^{p(k-1)} - N(\mathfrak{a})^{pk-1}) \\
\in & p \mathbb{Z}_{(p)}.
\end{align*}
On the other hand, if $P'/P$ does not act trivially on the $(\mathfrak{b},\eta)$, then the orbit of $(\mathfrak{b}, \eta)$ under the action of $P'/P$ in the above sum is 
\begin{align*}
& \sum_{g \in P'/P} \Big( \delta^{(y)}(gg_{\mathfrak{b}}g^{-1}) N(\mathfrak{b}^g)^{k-1} \Big) \\
= & |P'/P| \delta^{(y)}(g_{\mathfrak{b}})N(\mathfrak{b})^{k-1} \\
\in & p\mathbb{Z}_{(p)}.
\end{align*}
Here in the second equality we use that $\delta^{(y)}(gg_{\mathfrak{b}}g^{-1}) = \delta^{(y)}(g_{\mathfrak{b}})$ because $y$ is fixed under the action of $P'$. This proves the proposition. \qed

\begin{lemma} Let $P \in C(G)$ be such that there is no $P' \in C(G)$ with $P'^p = P$. Let $P \lneq N \leq N_GP$. Then the image of the transfer homomorphism 
\[
ver: N^{ab} \rightarrow P
\]
is a proper subgroup of $P$.
\label{imageunderver}
\end{lemma}

\noindent{\bf Proof:} Recall the definition of transfer map. Let $g \in N$. Let $\{x_1,\ldots,x_n\}$ be the double coset representatives $\langle g \rangle \backslash N / P$. Let $m$ be the smallest integer such that $g^m=1$. Then a set of left coset representatives of $P$ in $N$ is 
\[
\{1,g,\ldots,g^{m-1},x_1,gx_1,\ldots,g^{m-1}x_1,\ldots,x_n,gx_n,\ldots,g^{m-1}x_n\}.
\]
for all $0 \leq i \leq m-1$ and $1 \leq j \leq n$, we define $h_{ij}(g) \in P$ by
\[
g(g^ix_j) = g^{i'}x_{j'} h_{ij}(g).
\]
for a unique $0 \leq i' \leq m-1$ and $1 \leq j' \leq n$. Then
\[
h_{ij}(g) = \left\{
\begin{array}{l l}
1 & \text{if $i \leq m-2$} \\
x_j^{-1}g^mx_j & \text{if $i=m-1$}
\end{array} \right.
\]
Hence $ver(g) = \prod_{j=1}^{n}x_j^{-1}g^mx_j$. If $g \notin P$ then $g^m$ is not a generator of $P$ because of our assumption on $P$. Hence $ver(g)$ is not generator of $P$. On the other hand if $g \in P$ then 
\[
ver(g) = \prod_{x \in N/P} x^{-1}gx.
\]
Since both $P$ and $N$ are $p$-groups, the action of $N$ on the subgroup of order $p$ of $P$ is trivial. If $p^r$ is the order of $g$ then $N$ acts trivially on $g^{p^{r-1}}$. Hence 
\[
ver(g)^{p^{r-1}} = \prod_{x \in N/P} x^{-1}g^{p^{r-1}}x = \prod_{x \in N/P}g^{p^{r-1}} =1.
\]
Hence $ver(g)$ is not a generator of $P$. 
\qed

\begin{proposition} The sufficient condition in proposition (\ref{propsuffcongruence2}) for proving proposition (\ref{propcongruence2}) holds.
\end{proposition}

\noindent{\bf Proof:} We must show that for any $P \in C(G)$ such that there is no $P' \in C(G)$ with $P \subsetneq P'$ and any $j \geq 0$, any coset $y$ of $Z^{p^j}$ in $\Delta \times U_P$ whose image in $P$  is a generator of $P$ and any $u$ in $Z$ we have
\[
\Delta_P^{u^{d_P}}(\delta^{(y)},1-k) \equiv 0 (\text{mod } |(W_GP)_y| \mathbb{Z}_p),
\]
for any positive integer $k$ divisible by $[F_{\infty}(\mu_p):F_{\infty}]$. Choose an integral ideal $\mathfrak{f}$ of $O_{F_{W_P}}$ such that the Hilbert Eisenstein series $G_{k, \delta^{(y)}}$ over $\mathfrak{H}_{F_P}$, given by proposition (\ref{propeseries}), is defined on $\Gamma_{00}(\mathfrak{f}O_{F_P})$. Define 
\[
E= R_{F_P/F_{W_P}}(G_{k,\delta^{(y)}}).
\]
Then $E$ is a Hilbert modular form of weight $d_Pk$ on $\Gamma_{00}(\mathfrak{f})$. Let $\alpha$ be a finite id\`{e}le of $F_{W_P}$ whose image under the Artin symbol map coincides with $u$. Then by lemma \ref{constanttermofrestriction} and remark \ref{differenceofconstantterms}, we have
\[
c(0,E)-N(\alpha_p)^{-d_Pk}c(0,\alpha,E)= 2^{-r_P} \Delta_P^{u^{d_P}}(\delta^{(y)},1-k).
\]
Hence, using the $q$-expansion principle, it is enough to prove that the non-constant terms of the standard $q$-expansion of $E$ all lie in $|(W_GP)_y|\mathbb{Z}_{(p)}$ i.e. for any $\mu \in O_{F_{W_P}}^+$, 
\[
c(\mu,E) = \sum_{(\mathfrak{b},\nu)} \delta^{(y)}(g_{\mathfrak{b}})N(\mathfrak{b})^{k-1} \in |(W_GP)_y| \mathbb{Z}_{(p)}, 
\]
where $(\mathfrak{b},\nu)$ runs through all integral ideals $\mathfrak{b}$ of $F_P$ which are prime to $\Sigma_{F_P}$ and $\nu \in \mathfrak{b}$ is totally positive with $tr_{F_P/F_{W_P}}(\nu) = d_P\mu$. The group $(W_GP)_y$ acts on the pairs $(\mathfrak{b},\nu)$. Let $V$ be the stabiliser of $(\mathfrak{b},\nu)$. Then there is an integral ideal $\mathfrak{c}$ of $F_V:= F_P^V$ and a totally positive element $\eta$ of $O_{F_V}$ such that $\mathfrak{c}O_{F_P} = \mathfrak{b}$ and $\nu = \eta$. If $V$ is nontrivial group then $\delta^{(y)}(g_{\mathfrak{b}}) = 0$ by lemma (\ref{imageunderver}). On the other hand, if $V$ is trivial, then the $(W_GP)_y$ orbit of $(\mathfrak{b},\nu)$ in the above sum is 
\begin{align*}
& \sum_{g \in (W_GP)_y} \delta^{(y)}(gg_{\mathfrak{b}}g^{-1}) N(\mathfrak{b}^g)^{k-1} \\
= & |(W_GP)_y| \delta^{(y)}(g_{\mathfrak{b}}) N(\mathfrak{b})^{k-1} \\
\in & |(W_GP)_y| \mathbb{Z}_{(p)}.
\end{align*}
Here the second equality uses the fact that $\delta^{(y)}(gg_{\mathfrak{b}}g^{-1}) = \delta^{(y)}(g_{\mathfrak{b}})$ for any $g \in (W_GP)_y$. This proves the proposition. \qed

\begin{proposition} The sufficient condition in proposition (\ref{propsuffcongruence3}) for proving proposition (\ref{propcongruence3}) holds.
\end{proposition}

\noindent{\bf Proof:} We have to show that for any $P \in C(G)$ such that there exists a $P' \in C(G)$ with $P'^p=P$, any $j \geq 0$, any coset $y$ of $Z^{p^j}$ in $\Delta \times U_P$ and any $u$ in $Z$, we have
\[
\Delta_P^{u^{d_P}}(\delta^{(y)},1-k) \equiv \sum_{P'} \Delta_{P'}^{u^{d_P/p}}(\delta^{(y)} \circ \varphi, 1-pk) (\text{mod }|(W_{G_p}P)_y|\mathbb{Z}_p),
\]
for any positive integer $k$ divisible by $[F_{\infty}(\mu_p):F_{\infty}]$. The sum here runs through all $P' \in C(G)$ such that $[P':P]=p$. 

Choose an integral ideal $\mathfrak{f}$ of $F_{W_GP}$ such that the Hilbert Eisenstein series $G_{k,\delta^{(y)}}$ and $G_{pk, \delta^{(y)}\circ\varphi}$, given by proposition (\ref{propeseries}), on $\mathfrak{H}_{F_P}$ and $\mathfrak{H}_{F_{P'}}$ respectively are defined over $\Gamma_{00}(\mathfrak{f}O_{F_P})$ and $\Gamma_{00}(\mathfrak{f}O_{F_{P'}})$ respectively for every $P' \in C(G)$ such that $[P':P]=p$. We may assume that all prime factors of $\mathfrak{f}$ are in $\Sigma_{F_{W_{G_p}P}}$ and $\mathfrak{f} \subset d_PO_{F_{W_{G_p}P}}$. Define
\[
E = R_{F_P/F_{W_{G}P}}(G_{k,\delta^{(y)}})|_{d_Pk}U_{d_P} - \sum_{P'} R_{F_{P'}/F_{W_{G}P}}(G_{pk,\delta^{(y)}\circ \varphi})|_{d_Pk}U_{d_P/p}.
\]
Then $E \in M_{d_Pk}(\Gamma_{00}(\mathfrak{f}), \mathbb{C})$. Let $\alpha$ be a finite id\`{e}le of $F_{W_{G}P}$ whose image under the Artin symbol map coincides with $u$. Then by lemma \ref{constanttermofrestriction} and remark \ref{differenceofconstantterms}
\begin{align*}
& c(0,E)- N(\alpha_p)^{-d_Pk}c(0,\alpha,E) \\
& = 2^{-r_P} \Delta_P^{u^{d_P}}(\delta^{(y)},1-k) - \sum_{P'} 2^{-r_P/p}\Delta_{P'}^{u^{d_P/p}}(\delta^{(y)} \circ \varphi, 1-pk).
\end{align*}
As $2^{-r_P} \equiv 2^{-r_P/p} (\text{mod } r_P)$ and $r_P \geq |(W_GP)_y|$,
\begin{align*}
& 2^{-r_P} \Delta_P^{u^{d_P}}(\delta^{(y)},1-k) - \sum_{P'} 2^{-r_P/p}\Delta_{P'}^{u^{d_P/p}}(\delta^{(y)} \circ \varphi, 1-pk) \\
\equiv & 2^{-r_P} \Big( \Delta_P^{u^{d_P}}(\delta^{(y)},1-k) - \sum_{P'} \Delta_{P'}^{u^{d_P/p}}(\delta^{(y)} \circ \varphi, 1-pk) \Big) (\text{mod } |(W_GP)_y|\mathbb{Z}_p).
\end{align*}
Hence using the $q$-expansion principle it is enough to prove that the non-constant terms of the standard $q$-expansion of $E$ all lie in $|(W_GP)_y|\mathbb{Z}_{(p)}$ i.e. for all totally positive $\mu$ in $O_{F_{W_GP}}$, we have
\begin{align*}
& c(\mu,E) = c(d_P\mu, R_{F_P/F_{W_{G}P}}(G_{k,\delta^{(y)}})) - \sum_{P'} c(d_P \mu/p, R_{F_{P'}/F_{W_{G}P}}(G_{pk,\delta^{(y)}\circ \varphi})) \\
&=  \sum_{(\mathfrak{b},\eta)} \delta^{(y)}(g_{\mathfrak{b}})N(\mathfrak{b})^{k-1} - \sum_{P'} \sum_{(\mathfrak{a},\nu)} \delta^{(y)}(g_{\mathfrak{a}O_{F_P}})N(\mathfrak{a})^{pk-1} \in |(W_GP)_y|\mathbb{Z}_{(p)}.
\end{align*}
Here the pair $(\mathfrak{b},\eta)$ runs through all integral ideals $\mathfrak{b}$ of $F_P$ which are prime to $\Sigma_{F_P}$ and $\eta \in \mathfrak{b}$ is a totally positive element with $tr_{F_P/F_{W_GP}}(\eta) = d_P\mu$. The pair $(\mathfrak{a},\nu)$ runs through all integral ideals $\mathfrak{a}$ of $F_{P'}$ which are prime to $\Sigma_{F_{P'}}$ and $\nu \in \mathfrak{a}$ is a totally positive element with $tr_{F_{P'}/F_{W_GP}}(\nu) = d_P\mu/p$. The group $W_GP$ acts on the pairs $(\mathfrak{b},\eta)$ and $(\mathfrak{a},\nu)$. Let $V \subset (W_GP)_y$ be the stabiliser of $(\mathfrak{b},\eta)$. Then there is an integral ideal $\mathfrak{c}$ of $F_V := F_P^V$ and a totally positive element $\gamma$ in $O_{F_V}$ such that $\mathfrak{c}O_{F_P} = \mathfrak{b}$ and $\eta = \gamma$. Then the $(W_GP)_y$ orbit of $(\mathfrak{b},\eta)$ in the above sum is 
\begin{align*} 
& \sum_{g \in (W_GP)_y/V} \Big( \delta^{(y)}(gg_{\mathfrak{b}}g^{-1}) N(\mathfrak{b}^g)^{k-1} - \sum_{P'} \delta^{(y)}(gg_{\mathfrak{b}}g^{-1})N(\mathfrak{b}^g)^{pk-1} \Big) \\
= & |(W_GP)_y/V| \delta^{(y)}(g_{\mathfrak{b}}) \Big(N(\mathfrak{b})^{k-1} - N(\mathfrak{b})^{pk-1} \Big) \\
= & |(W_GP)_y/V| \delta^{(y)}(g_{\mathfrak{b}}) \Big( N(\mathfrak{c})^{|V|(k-1)} - N(\mathfrak{c})^{|V|(pk-1)/p} \Big) \\
\in & |(W_GP)_y| \mathbb{Z}_{(p)}.
\end{align*} 
Here in the first line $P'$ runs through all $P' \in C(G)$ with $[P':P]=p$ such that $P'/P \subset V$. The second sum if 0 if $V$ is trivial and in that case inclusion in the last line is trivial. The first equality uses $\delta^{(y)}(gg_{\mathfrak{b}}g^{-1}) = \delta^{(y)}(g_{\mathfrak{b}})$ as $g \in (W_GP)_y$. The last inclusion is because $N(\mathfrak{c})^{|V|} \equiv N(\mathfrak{c})^{|V|/p} (\text{mod }|V|)$. This proves the proposition. \qed

\begin{proposition} The sufficient conditions in proposition (\ref{propsuffcongruence4}) for proving proposition (\ref{propcongruence4}) hold.
\end{proposition}

\noindent{\bf Proof:} We just prove the sufficient condition for congruence (\ref{congruence4}). Proof of the other sufficient condition in proposition (\ref{propsuffcongruence4}) is similar. We must show that for any $P \in C(G)$ and any $P' \in C(G)$ such that $[P':P]=p$, for any $j \geq 0$, any coset $y$ of $Z^{p^j}$ in $\Delta \times U_P$ and any $u$ in $Z$
\[
\Delta_P^{u^{pd_{P'}}}(\delta^{(y)},1-k) \equiv \Delta_{P'}^{u^{d_{P'}}}(\delta^{(y)} \circ \varphi, 1-pk) (\text{mod } |(N_GP'/P)_y| \mathbb{Z}_p),
\]
for any positive integer $k$ divisible by $[F_{\infty}(\mu_p):F_{\infty}]$. 

Choose an integral ideal $\mathfrak{f}$ of $F_{W_GP'}$ such that the Hilbert Eisenstein series $G_{k,\delta^{(y)}}$ and $G_{pk, \delta^{(y)}\circ \varphi}$, given by proposition (\ref{propeseries}), on $\mathfrak{H}_{F_P}$ and $\mathfrak{H}_{F_{P'}}$ respectively are defined over $\Gamma_{00}(\mathfrak{f}O_{F_P})$ and $\Gamma_{00}(\mathfrak{f}O_{F_{P'}})$. Moreover, we may assume that all prime factors of $\mathfrak{f}$ are in $\Sigma_{F_{W_GP'}}$ and $\mathfrak{f} \subset pd_{P'} O_{F_{W_GP'}}$. Define 
\[
E= R_{F_P/F_{W_GP'}}(G_{k,\delta^{(y)}})|_{pd_{P'}k}U_{pd_{P'}} - R_{F_{P'}/F_{W_GP'}}(G_{pk,\delta^{(y)}\circ \varphi})|_{pd_{P'}k}U_{d_{P'}}.
\]
Then $E \in M_{pd_{P'}k}(\Gamma_{00}(\mathfrak{f}), \mathbb{C})$. Let $\alpha$ be a finite id\`{e}le of $F_{W_GP'}$ whose image under the Artin symbol map coincides with $u$. Then by lemma \ref{constanttermofrestriction} and remark \ref{differenceofconstantterms} 
\[
c(0,E) - N(\alpha_p)^{-pd_{P'}} c(0,\alpha,E)= 2^{-r_P}\Delta_P^{u^{pd_{P'}}}(\delta^{(y)},1-k) - 2^{-r_{P'}} \Delta_{P'}^{u^{d_{P'}}}(\delta^{(y)},1-pk).
\]
As $2^{-r_P} \equiv 2^{-r_{P'}} (\text{mod }r_P)$ and $r_P \geq |(N_GP'/P)_y|$,
\begin{align*}
& 2^{-r_P}\Delta_P^{u^{pd_{P'}}}(\delta^{(y)},1-k) - 2^{-r_{P'}} \Delta_{P'}^{u^{d_{P'}}}(\delta^{(y)},1-pk) \\
\equiv & 2^{-r_P} \Big( \Delta_P^{u^{pd_{P'}}}(\delta^{(y)},1-k) - \Delta_{P'}^{u^{d_{P'}}}(\delta^{(y)},1-pk) \Big) (\text{mod } |(N_GP'/P)_y| \mathbb{Z}_p).
\end{align*}
Hence using the $q$-expansion principle it is enough to prove that the non-constant terms of the standard $q$-expansion of $E$ all lie in $|(N_GP'/P)_y| \mathbb{Z}_{(p)}$ i.e. for all totally positive $\mu$ in $O_{F_{W_GP'}}$ we have
\begin{align*}
&c(\mu,E) = c(pd_{P'}\mu, R_{F_P/F_{W_GP'}}(G_{k,\delta^{(y)}})) - c(d_{P'}\mu, R_{F_{P'}/F_{W_GP'}}(G_{pk,\delta^{(y)}\circ \varphi})) \\
& =\sum_{(\mathfrak{b},\eta)} \delta^{(y)}(g_{\mathfrak{b}})N(\mathfrak{b})^{k-1} - \sum_{(\mathfrak{a},\nu)} \delta^{(y)}(g_{\mathfrak{a}O_{F_P}}) N(\mathfrak{a})^{pk-1} \in |(N_GP'/P)_y| \mathbb{Z}_{(p)}. 
\end{align*}
Here the pairs $(\mathfrak{b},\eta)$ runs through all integral ideals $\mathfrak{b}$ of $F_P$ which are prime to $\Sigma_{F_P}$ and $\eta \in \mathfrak{b}$ is a totally positive element with $tr_{F_P/F_{W_GP'}}(\eta) = pd_{P'} \mu$. The pairs $(\mathfrak{a},\nu)$ runs through all integral ideals $\mathfrak{a}$ of $F_{P'}$ which are prime to $\Sigma_{F_{P'}}$ and $\nu \in \mathfrak{a}$ is a totally positive element with $tr_{F_{P'}/F_{W_GP'}}(\nu) = d_{P'}\mu$. The group $N_GP'/P$ acts on the pairs $(\mathfrak{b},\eta)$ and $(\mathfrak{a},\nu)$. Let $V \subset (N_GP'/P)_y$ be the stabiliser of $(\mathfrak{b},\eta)$. Then there is an integral ideal $\mathfrak{c}$ of $F_V:= F_P^V$ and a totally positive element $\gamma$ of $O_{F_V}$ such that $\mathfrak{c}O_{F_P}=\mathfrak{b}$ and $\eta =\gamma$. First assume that $V$ is a non-trivial group. Then the $(N_GP'/P)_y$ orbit of $(\mathfrak{b},\eta)$ in the above sum is
\begin{align*}
& \sum_{g \in (N_GP'/P)_y/V} \Big( \delta^{(y)}(gg_{\mathfrak{b}}g^{-1}) N(\mathfrak{b}^g)^{k-1} - \delta^{(y)}(gg_{\mathfrak{b}}g^{-1})N(\mathfrak{b}^g)^{pk-1} \Big) \\
= & |(N_GP'/P)_y/V| \delta^{(y)}(g_{\mathfrak{b}})\Big( N(\mathfrak{b})^{k-1} - N(\mathfrak{b})^{pk-1} \Big) \\
= & |(N_GP'/P)_y/V| \delta^{(y)}(g_{\mathfrak{b}})\Big( N(\mathfrak{c})^{|V|(k-1)} - N(\mathfrak{c})^{|V|(pk-1)/p} \Big) \\
\in &  |(N_GP'/P)_y|\mathbb{Z}_{(p)}.
\end{align*}
On the other hand if $V$ is a trivial group then the $(N_GP'/P)_y$ orbit of the pair $(\mathfrak{b},\eta)$ in the above sum is
\[
 \sum_{g \in (N_GP'/P)_y}  \delta^{(y)}(gg_{\mathfrak{b}}g^{-1}) N(\mathfrak{b}^g)^{k-1} = |(N_GP'/P)_y| \delta^{(y)}(g_{\mathfrak{b}}) N(\mathfrak{b})^{k-1}.
\]
In both cases the first equality uses $\delta^{(y)}(gg_{\mathfrak{b}}g^{-1}) = \delta^{(y)}(g_{\mathfrak{b}})$ for $g \in (N_GP'/P)_y$. In the first case we also use the fact that $N(\mathfrak{c})^{|V|} \equiv N(\mathfrak{c})^{|V|/p} (\text{mod } |V|)$. This proves the proposition. \qed

\label{sectionofproofofcongruences}

\subsection{Proof of M4. from the basic congruences} \label{subsectionm3frombasiccongruences} We have proved the basic congruences in previous subsections. We want to deduce M4 from these congruences. However, we cannot do it directly for the extension $F_{\infty}/F$. We use the following trick: we extend our field slightly to $\tilde{F}_{\infty} \supset F_{\infty}$ such that $\tilde{F}_{\infty}/F$ is an admissible $p$-adic Lie extension satisfying $\mu=0$ hypothesis and $Gal(\tilde{F}_{\infty}/F) = \Delta \times \tilde{\mathcal{G}}$ with $\tilde{\mathcal{G}} \cong \tilde{H} \times \mathcal{G}$, where $\tilde{H}$ is a cyclic group of order $|G|$. We know the basic congruences for $\tilde{F}_{\infty}/F$ which we use to deduce the M4 for $\tilde{F}_{\infty}/F$. This proves the main conjecture for $\tilde{F}_{\infty}/F$ and hence implies the main conjecture for $F_{\infty}/F$.

\subsubsection{The field $\tilde{F}_{\infty}$} Choose a prime $l$ large enough such that $l \equiv 1 (\text{mod } |G|)$ and $\mathbb{Q}(\mu_l) \cap F_{\infty} = \mathbb{Q}$. Let $K$ be the extension of $\mathbb{Q}$ contained in $\mathbb{Q}(\mu_l)$ such that $[K:\mathbb{Q}] = |G|$. Define $\tilde{F} = KF$ and $\tilde{F}_{\infty} = \tilde{F}F_{\infty}$. Then 
\[
Gal(\tilde{F}_{\infty}/F) = Gal(\tilde{F}/F) \times Gal(F_{\infty}/F) =: \tilde{H} \times \Delta \times \mathcal{G} =: \tilde{\mathcal{G}}.
\]
Then define $\tilde{G} = \tilde{\mathcal{G}}/Z \cong \tilde{H} \times \Delta \times G$. 

\subsubsection{A key lemma} We extend the field $F_{\infty}$ to $\tilde{F}_{\infty}$ as we need the following key lemma. For any $P \in C(\tilde{G})$, define the integer $i_P$ by
\[
i_P = max_{P' \in C(\tilde{G})} \{[P':P] | P \subset P' \}
\]

\begin{lemma} Let $P \in C(\tilde{G})$. If $P \neq \{1\}$, then 
\[
T_P \subset pi_P^2 \Lambda( \Delta \times U_P)
\]
Clearly
\[
T_{\{1\}} = |\tilde{G}| \Lambda(\Delta \times Z).
\]
Similar statement holds for $T_{P,S}$ and $\widehat{T}_P$.
\label{keylemma}
\end{lemma}

\noindent{\bf Proof:} {\bf Case 1:} $P \subset \tilde{H}$. Then $i_P = [\tilde{H}:P]$ and $N_{\tilde{G}}P = \tilde{G}$ acts trivially on $\Lambda(\Delta \times U_P)$. Hence
\[
T_P = [\tilde{G}:P] \Lambda(\Delta \times U_P) = |G|[\tilde{H}:P]\Lambda(\Delta \times U_P).
\]
If $P \neq \{1\}$, then $|G| \geq pi_P$. Hence the claim. 

\noindent{\bf Case 2:} $P \nsubseteq \tilde{H}$. Let $P$ be generated by $(\tilde{h}, h)$, with $\tilde{h} \in \tilde{H}$ and $h \in G$. By assumption $h \neq 1$. Let $P' \in C(\tilde{G})$ such that $[P':P] = i_P$. Let $(\tilde{h}_0,h_0)$ be a generator of $P'$ such that $\tilde{h}_0^{i_P} = \tilde{h}$ and $h_0^{i_P} = h$. Now note that
\[
\tilde{H} \times \langle h_0 \rangle \subset N_{\tilde{G}}P
\]
acts trivially on $\Lambda(\Delta \times U_P)$. As $P \subset \tilde{H} \times \langle h_0 \rangle$ this implies that
\begin{align*}
T_P & \subset \frac{|\tilde{H} \times \langle h_0 \rangle|}{|P|} \Lambda(\Delta \times U_P) \\
& = \frac{|\tilde{H}||\langle h_0 \rangle|}{|P|} \Lambda(\Delta \times U_P) \\
& = |\tilde{H}|i_P \Lambda(\Delta \times U_P) \\
& \subset pi_P^2 \Lambda(\Delta \times U_P).
\end{align*}
The last containment holds because $|\tilde{H}| \geq pi_P$. The same proof goes through for the statement about $T_{P,S}$ and $\widehat{T}_P$. \qed

\subsubsection{Completion of the proof}
\begin{lemma} For any $P \in C(\tilde{G})$ and any $0 \leq k \leq p-1$, we have
\[
\zeta_P - \omega_P^k(\zeta_P) \in \mathfrak{p} \frac{T_{P,S}}{i_P}.
\] 
Hence $\zeta_P^p/\prod_{k=0}^{p-1}\omega_P^k(\zeta_P) \in 1+pT_{P,S}/i_P$.
\label{lastlemma}
\end{lemma}

\noindent{\bf Proof:} We use reverse induction on $|P|$. When $P$ is a maximal cyclic subgroup $i_P = 1$ and the required congruence is proven in proposition \ref{propcongruence2}. In general we use the congruence in proposition \ref{propcongruence3} so that 
\begin{align}
\zeta_P - \omega_P^k(\zeta_P) & \equiv \sum_{P'} \Big( \varphi(\zeta_{P'}) - \omega_P^k(\varphi(\zeta_{P'}))\Big) \\
& = \sum_{P'}  \varphi(\zeta_{P'} - \omega_{P'}^k(\zeta_{P'})) (\text{mod } \mathfrak{p}T_{P,S}), 
\label{lastcongruence}
\end{align}
for appropriately chosen $\omega_P$ and $\omega_{P'}$. But by induction hypothesis 
\[
\zeta_{P'} - \omega_{P'}^k(\zeta_{P'}) \in \mathfrak{p} \frac{T_{P',S}}{i_{P'}}.
\]
Now for any $P' \in C(\tilde{G})$ such that $[P':P]=p$, note that 
\[
\varphi(\sum_{ x \in N_{\tilde{G}}P/N_{\tilde{G}}P'} xT_{P',S}x^{-1}) \subset \frac{T_{P,S}}{p}.
\]
This finishes the proof of the first assertion noting that $i_P = pi_{P'}$. Hence 
\[
\zeta_P^p \equiv \prod_{k=0}^{p-1} \omega_P^k(\zeta_P) (\text{mod } \mathfrak{p}T_{P,S}/i_P).
\]
But since both sides of the congruence are invariant under the action of $Gal(\mathbb{Q}_p(\mu_p)/\mathbb{Q}_p)$, we get 
\[
\zeta_P^p \equiv \prod_{k=0}^{p-1} \omega_P^k(\zeta_P) (\text{mod }pT_{P,S}/i_P).
\]
\qed

Hence 
\[
log(\frac{\omega_P^k(\zeta_P)}{\zeta_P}) \equiv 1- \frac{\omega_P^k(\zeta_P)}{\zeta_P} \text{mod } (\mathfrak{p} \widehat{T}_P/i_P)^2,
\]
which implies
\[
log\Big( \frac{\prod_{k=0}^{p-1} \omega_P^k(\zeta_P)}{\zeta_P^p} \Big) \equiv \frac{p\zeta_P - \sum_{k=0}^{p-1} \omega_P^k(\zeta_P)}{\zeta_P} (\text{mod } (p\widehat{T}_P/i_P)^2).
\]
Then 
\begin{align*}
& log\Big( \frac{\prod_{P'} \varphi(\alpha_{P'}(\zeta_{P'}))}{\zeta_P^p/\prod_{k=0}^{p-1}\omega_P^k(\zeta_P)}\Big) \\ 
&\equiv \sum_{P'}\Big(\frac{p\varphi(\zeta_{P'}) - \sum_{k=0}^{p-1}\omega_P^k(\varphi(\zeta_{P'}))}{\varphi(\zeta_{P'})}\Big) - \Big(\frac{p\zeta_P-\sum_{k=0}^{p-1} \omega_P^k(\zeta_P)}{\zeta_P}\Big) \\
& \equiv \sum_{P'}\Big(\frac{p\varphi(\zeta_{P'}) - \sum_{k=0}^{p-1}\omega_P^k(\varphi(\zeta_{P'}))}{\varphi(\zeta_{P'})}\Big) - \sum_{P'} \Big(\frac{p\varphi(\zeta_{P'}) - \sum_{k=0}^{p-1}\omega_P^k(\zeta_{P'})}{\zeta_P}\Big) \\
& \equiv \sum_{P'} \frac{(p\varphi(\zeta_{P'}) - \sum_{k-0}^{p-1}\omega_P^k(\varphi(\zeta_{P'})))(\zeta_P-\varphi(\zeta_{P'}))}{\zeta_P\varphi(\zeta_{P'}} (\text{mod }p\widehat{T}_P).
\end{align*}
Here we use $(p\widehat{T}_P/i_P)^2 \subset p \widehat{T}_P$ as implied by lemma \ref{keylemma}. The second congruence above uses congruence \ref{lastcongruence}. Now note that 
\[
p \varphi(\zeta_{P'}) - \sum_{k=0}^{p-1} \omega_P^k(\varphi(\zeta_{P'})) \in p\varphi(\widehat{T}_{P'}/i_{P'}) \hspace{.5cm} \text{(by lemma \ref{lastlemma})}
\]
and
\[
\zeta_P - \varphi(\zeta_{P'}) \in \widehat{T}_{P,N_{\tilde{G}}P'} \hspace{.5cm} \text{(by congruence (\ref{congruence4}) and (\ref{congruence5}))}.
\]
Hence
\[
\big(p \varphi(\zeta_{P'}) - \sum_{k=0}^{p-1} \omega_P^k(\varphi(\zeta_{P'}))\big)\big(\zeta_P - \varphi(\zeta_{P'})\big) \in p\varphi(\widehat{T}_{P'}/i_{P'})\widehat{T}_{P,N_{\tilde{G}}P'} \subset p\widehat{T}_{P,N_{\tilde{G}}P'}.
\]
Which in turn implies that 
\[
\sum_{P'} \Big((p\varphi(\zeta_{P'}) - \sum_{k-0}^{p-1}\omega_P^k(\varphi(\zeta_{P'})))(\zeta_P-\varphi(\zeta_{P'}))\Big) \in p\widehat{T}_P.
\]
Hence 
\[
log\Big( \frac{\prod_{P'} \varphi(\alpha_{P'}(\zeta_{P'}))}{\zeta_P^p/\prod_{k=0}^{p-1}\omega_P^k(\zeta_P)}\Big) \in p\widehat{T}_P.
\]
As $log$ induces an isomorphism between $1+p\widehat{T}_P$ and $p\widehat{T}_P$, we have
\[
\frac{\prod_{P'} \varphi(\alpha_{P'}(\zeta_{P'}))}{\zeta_P^p/\prod_{k=0}^{p-1}\omega_P^k(\zeta_P)} \in 1+p\widehat{T}_P.
\]
But by lemma \ref{lastlemma} $\frac{\prod_{P'} \varphi(\alpha_{P'}(\zeta_{P'}))}{\zeta_P^p/\prod_{k=0}^{p-1}\omega_P^k(\zeta_P)} \in 1+pT_{P,S}/i_P$ and 
\[
1+p\widehat{T}_P \cap 1+pT_{P,S}/i_P = 1+pT_{P,S}.
\]
Hence 
\[
\frac{\prod_{P'} \varphi(\alpha_{P'}(\zeta_{P'}))}{\zeta_P^p/\prod_{k=0}^{p-1}\omega_P^k(\zeta_P)}  \in 1+pT_{P,S}.
\]
When $P \neq \{1\}$, this is the required congruence M4. When $P = \{1\}$, note that $\prod_{k=0}^{p-1}\omega_{\{1\}}^k(\zeta_{\{1\}}) = \zeta_0$. This can be seen either by interpolation properties of $\prod_{k=0}^{p-1}\omega_{\{1\}}^k(\zeta_{\{1\}})$ and $\zeta_0$ or by lemma \ref{normandtrace}. Hence we get 
\[
\frac{\prod_{P'} \varphi(\alpha_{P'}(\zeta_{P'}))}{\zeta_{\{1\}}^p/\zeta_0} \in 1+pT_{\{1\},S}.
\]
Now use the basic congruence (\ref{congruence5}) which says $\zeta_0 \equiv \varphi(\zeta_{\{1\}}) (\text{mod }p|\tilde{G}|)$. Note that $T_{\{1\},S} = |\tilde{G}|\Lambda(\Delta \times Z)_S$. Hence
\[
\frac{\prod_{P'} \varphi(\alpha_{P'}(\zeta_{P'}))}{\zeta_{\{1\}}^p/\varphi(\zeta_{\{1\}})} \in 1+pT_{\{1\},S}.
\]
This is M4 for $P=\{1\}$. This finishes proof of the main theorem \ref{mainconjecture}.

\bibliographystyle{plain}
\bibliography{mybib2}

\end{document}